\documentclass[10pt]{article}
\usepackage{amssymb,latexsym}
\usepackage{epsfig}
\usepackage{eufrak}
\usepackage{amsmath}
\usepackage{mathrsfs}
\usepackage{color}

\newtheorem{theorem}{Theorem}
\newtheorem{lemma}{Lemma}
\newtheorem{corollary}{Corollary}

\newcommand{\be}{\begin{equation}}
\newcommand{\ee}{\end{equation}}
\newcommand{\bee}{\begin{eqnarray*}}
\newcommand{\eee}{\end{eqnarray*}}
\newcommand{\bel}{\begin{eqnarray}}
\newcommand{\eel}{\end{eqnarray}}
\newcommand{\bec}{\begin{cases}}
\newcommand{\eec}{\end{cases}}
\newcommand{\bem}{\begin{bmatrix}}
\newcommand{\eem}{\end{bmatrix}}

\newcommand{\la}{\label}
\newcommand{\li}{\left}
\newcommand{\ri}{\right}

\newcommand{\ovl}{\overline}

\newcommand{\lf}{\lfloor}
\newcommand{\rf}{\rfloor}

\newcommand{\vep}{\varepsilon}
\newcommand{\lm}{\lambda}
\newcommand{\Lm}{\Lambda}

\newcommand{\si}{\sigma}
\newcommand{\Si}{\Sigma}

\newcommand{\de}{\delta}

\newcommand{\vDe}{\varDelta}

\newcommand{\ga}{\gamma}
\newcommand{\Ga}{\Gamma}
\newcommand{\vse}{\vartheta}
\newcommand{\se}{\theta}
\newcommand{\Se}{\Theta}

\newcommand{\ze}{\zeta}
\newcommand{\al}{\alpha}
\newcommand{\ba}{\beta}

\newcommand{\vro}{\varrho}
\newcommand{\ro}{\rho}
\newcommand{\ka}{\kappa}
\newcommand{\om}{\omega}
\newcommand{\Om}{\Omega}

\newcommand{\f}{\frac}
\newcommand{\sq}{\sqrt}
\newcommand{\cd}{\cdots}

\newcommand{\qu}{\quad}
\newcommand{\qqu}{\qquad}

\newcommand{\mscr}{\mathscr}
\newcommand{\mcal}{\mathcal}

\newcommand{\bb}{\mathbb}

\newcommand{\wh}{\widehat}

\newcommand{\mrm}{\mathrm}
\newcommand{\bs}{\boldsymbol}

\newcommand{\ap}{\approx}

\newcommand{\sh}{\slash}

\newcommand{\tx}{\text}

\newcommand{\iy}{\infty}

\newcommand{\pa}{\partial}

\newcommand{\bed}{\begin{description}}
\newcommand{\eed}{\end{description}}
\newcommand{\bei}{\begin{itemize}}
\newcommand{\eei}{\end{itemize}}
\newcommand{\ben}{\begin{enumerate}}
\newcommand{\een}{\end{enumerate}}
\newcommand{\bib}{\bibitem}
\newcommand{\beL}{\begin{lemma}}
\newcommand{\eeL}{\end{lemma}}
\newcommand{\beT}{\begin{theorem}}
\newcommand{\eeT}{\end{theorem}}
\newcommand{\sect}{\section}

\newcommand{\bpf}{\begin{pf}}
\newcommand{\epf}{\end{pf}}
\newcommand{\bsk}{\bigskip}
\newcommand{\bi}{\binom}

\newcommand{\pfbox}{\hfill\mbox{$\Box$}}

\newenvironment{pf}{\paragraph*{Proof{\rm.}}}{\pfbox\bigskip}

\begin{document}

\title{{\bf A Likelihood Ratio Approach for Probabilistic Inequalities}
\thanks{The author is working with the Department of Electrical Engineering and Computer Sciences at Louisiana
State University at Baton Rouge, LA 70803, USA, and the Department of Electrical Engineering at Southern University and A\&M College, Baton
Rouge, LA 70813, USA; Email: chenxinjia@gmail.com}}

\author{Xinjia Chen}

\date{August 18,  2013}

\maketitle

\begin{abstract}

We propose a new approach for deriving probabilistic inequalities based on bounding likelihood ratios.  We demonstrate that this approach is
more general and powerful than the classical method frequently used for deriving concentration inequalities such as Chernoff bounds. We discover
that the proposed approach is inherently related to statistical concepts such as monotone likelihood ratio,  maximum likelihood, and the method
of moments for parameter estimation.  A connection between the proposed approach and the large deviation theory is also established. We show
that, without using moment generating functions,  tightest possible concentration inequalities may be readily derived by the proposed approach.
We have derived new concentration inequalities using the proposed approach, which cannot be obtained by the classical approach based on moment
generating functions.

\end{abstract}


\sect{Introduction}

A problem of ubiquitous significance in engineering and sciences is to obtain bounds for probabilities of random events.  Formally, let $E$ be
an event defined in probability space $(\Om, \Pr , \mscr{F})$, where $\Om$ is the sample space, $\Pr$ denotes the probability measure, and
$\mscr{F}$ is the the $\si$-algebra.  In many situations, it is desirable to obtain bounds for $\Pr \{ E \}$.  For theoretical and practical
purposes, the bounds are expected to be simple, insightful and as tight as possible.   In general, the event $E$ can be expressed in terms of a
matrix-valued random variable $\bs{X}$.  In particular, $\bs{X}$ can be a random vector or scalar.   Clearly, the event $E$ can be represented
as $\{ \bs{X} \in \mscr{E} \}$, where $\mscr{E}$ is a certain set of deterministic matrices. In probability theory, one of the most frequently
used approach for deriving inequalities for  $\Pr \{ E \}$ is to bound the indicator function $\bb{I}_{ \{ \bs{X} \in \mscr{E} \} }$ by a family
of random variables having finite expectation and minimize the expectation.  The central idea of this approach is to seek a family of bounding
functions $w(\bs{X}, \bs{\vse})$ of $\bs{X}$, parameterized by $\bs{\vse} \in \varTheta$, such that \be \la{inviewvip}
 \bb{I}_{ \{  \bs{X} \in \mscr{E}   \} } \leq
w(\bs{X}, \bs{\vse}) \qu \tx{for all $\bs{\vse} \in \varTheta$}. \ee Here, the notion of inequality (\ref{inviewvip}) is that the inequality
$\bb{I}_{ \{  \bs{X} (\om) \in \mscr{E}   \} } \leq w(\bs{X} (\om) , \bs{\vse})$ holds for every $\om \in \Om$. As a consequence of the
monotonicity of the mathematical expectation $\bb{E} [. ]$, \be \la{inviewvip889b}
 \Pr \{  E \}
= \bb{E} [ \bb{I}_{ \{  \bs{X} \in \mscr{E}   \} }  ] \leq \bb{E} [ w(\bs{X}, \bs{\vse}) ] \qu \tx{for all $\bs{\vse} \in \varTheta$}. \ee
Minimizing the upper bound in (\ref{inviewvip889b}) with respect to $\bs{\vse} \in \varTheta$ yields  \be \la{inviewvip889c} \Pr \{  E \} \leq
\inf_{ \bs{\vse} \in \varTheta } \bb{E} [ w(\bs{X}, \bs{\vse}) ]. \ee Classical inequalities such as Chebyshev inequality, Markov inequality,
Cantelli's inequality, Bernstein's inequalities \cite{Bernstein}, Chernoff bounds \cite{Chernoff}, Benett's inequalities \cite{Bennett},
Hoeffding's inequalities \cite{Hoeffding} can be derived by this approach with various bounding functions $w(\bs{X}, \bs{\vse})$, where $\bs{X}$
is a scalar random variable. To illustrate, the derivation of some fundamental inequalities is presented follows.

\bed

\item [Chebyshev inequality:]  Let $X$ be a random variable with mean $\mu$ such that $\bb{E} [|X - \mu|^s ] < \iy$ for $s \in \Se$,
where $\Se$ is a subset of positive real numbers.  Let $\vep > 0$. Define $w( X, s) = \li ( \f{| X - \mu |}{\vep} \ri )^s$ for $s \in \Se$.
Then, $\bb{I}_{ \{ | X - \mu | \geq \vep \} } \leq w( X, s)$ for $s \in \Se$. It follows from (\ref{inviewvip889c}) that $\Pr \{ | X - \mu |
\geq \vep \} \leq \inf_{s \in \Se} \f{\bb{E}[ | X - \mu |^s ]}{\vep^s}$,  which is referred to as Chebyshev inequality in the particular case
that $\Se = \{ 2 \}$.

\item [Markov inequality:] Let $X$ be a non-negative random variable such that $\bb{E} [X^s ] < \iy$ for $s \in \Se$,
where $\Se$ is a subset of positive real numbers.  Let $\ga > 0$. Define $w(X, s) = \li ( \f{X}{\ga} \ri )^s$ for $s \in \Se$.  Then, $\bb{I}_{
\{ X \geq \ga \} } \leq w( X, s)$ for  $s \in \Se$. It follows from (\ref{inviewvip889c}) that $\Pr \{  X \geq \ga \} \leq \inf_{s \in \Se} \f{
\bb{E} [X^s ] }{\ga^s}$, which is referred to as Markov inequality in the particular case that $\Se = \{1\}$.

\item [Cantelli's inequality:] Let $X$ be a random variable with mean $\mu$ and finite variance $\si^2 > 0$.
Let $\vep > 0$.  Define $w(X, s) = \li ( \f{X + s}{\vep + \mu + s} \ri )^2$
for $s > - ( \vep + \mu )$. Then, $\bb{I}_{ \{ X - \mu \geq \vep \} } \leq w( X, s)$ and $\bb{E} [ w(X, s) ] =  \f{(\mu + s)^2 + \si^2}{(\vep +
\mu + s)^2}$  for all $s > - ( \vep + \mu )$.  It follows from (\ref{inviewvip889c}) that $\Pr \{  X - \mu \geq \vep \} \leq \inf_{s
> - ( \vep + \mu )} \f{(\mu + s)^2 + \si^2}{(\vep + \mu + s)^2} = \f{\si^2}{\si^2 + \vep^2}$.

\item [Chernoff bound:]  Let $X$ be a random variable such that the moment generating function $\phi(s) = \bb{E} [e^{s X}  ]$ is finite for $s \in \Se$,
where $\Se$ is a subset of positive real numbers.  Let $\ga$ be a real number. Define $w(X, s) = e^{s (X-\ga)} $ for $s \in \Se$. Then,
$\bb{I}_{ \{ X \geq \ga \} } \leq w( X, s)$ for $s \in \Se$. It follows from (\ref{inviewvip889c}) that $\Pr \{  X \geq \ga \} \leq \inf_{s \in
\Se} \bb{E} [ w(X, s) ] = \inf_{s \in \Se}  e^{- \ga s} \bb{E} [ e^{s X} ]$.

\eed

As can be seen from the above discussion, the central idea of this approach is to seek a family of bounding functions  $w (\bs{X}, \bs{\vse}),
\; \bs{\vse} \in \varTheta$ such that the mathematical expectation $\bb{E} [ w (\bs{X}, \bs{\vse}) ]$ is convenient for evaluation and
minimization.  We call this technique of deriving probabilistic inequalities as the {\it mathematical expectation} (ME) method, in view of the
crucial role played by the mathematical expectation of bounding functions.  The ME method is a very general approach. However, it has two
drawbacks.

First, in some situations, the mathematical expectation $\bb{E} [ w (\bs{X}, \bs{\vse}) ]$ may be intractable.  For example, if $X$ is a
hypermetric random variable and one wishes to derive the Chernoff bound for $\Pr \{ X \geq \ga \}$ with bounding function $w (X, s) = e^{s
(X-\ga)}$, it is not possible to obtain an tractable expression for $\bb{E} [ w (X, s) ]$ except some crude bounds. Similar difficulty is
encountered if $X$ is a generalized Poisson random variable. In the general case that $\bs{X}$ is matrix-valued or vector-valued, it may be more
cumbersome to evaluate and minimize $\bb{E} [ w (\bs{X}, \bs{\vse}) ]$.

Second, the ME method may not fully exploit the information of the underlying distribution, since the mathematical expectation is a only
quantity of summary for the distribution.  The issue may arise when the probability density or mass function of $\bs{X}$ is available.

In this paper,  we would like to propose a more general approach for deriving probabilistic inequalities, aiming at overcoming the drawbacks of
the ME method.  Let $f(.)$ denote the probability density function (pdf) or probability mass function (pmf) of $\bs{X}$. The primary idea of the
proposed approach is to seek a family of pdf or pmf $g(., \bs{\vse})$, parameterized by $\bs{\vse} \in \varTheta$,   and a deterministic
function $\Lm(\bs{\vse})$ of $\bs{\vse} \in \varTheta$  such that  for all $\bs{\vse} \in \varTheta$,  the indicator function $\bb{I}_{ \{
\bs{X} \in \mscr{E}   \} }$ is bounded from above by the product of $\Lm(\bs{\vse})$ and the likelihood ratio $\f{ g(\bs{X}, \bs{\vse}) }{
f(\bs{X}) }$. Then, the probability $\Pr \{ \bs{X} \in \mscr{E}   \}$ is bounded from above by the infimum of $\Lm(\bs{\vse})$ with respect to
$\bs{\vse} \in \varTheta$.  Due to the central role played by the likelihood ratio, this technique of deriving probabilistic inequalities is
referred to as the {\it likelihood ratio} (LR) method.  It can be demonstrated that the ME method is actually a special technique of the LR
method. From the perspective of the LR method, the idea of the ME method is to construct a family of pdf or pmf $g(., \bs{\vse})$ via
multiplying the pdf or pmf $f (. )$ of $\bs{X}$ by the bounding function $w(. , \bs{\vse})$ and performing a normalization.  In other words, the
bounding function $w(. , \bs{\vse})$ is used as a {\it weight function} to modify the distribution of $\bs{X}$.  The deterministic function
$\Lm(\bs{\vse})$ can be taken as $\bb{E} [ w (\bs{X}, \bs{\vse}) ]$. This approach for constructing the parameterized distribution $g(.,
\bs{\vse})$ and the function $\Lm(\bs{\vse})$ is useful provided that the mathematical expectation $\bb{E} [ w (\bs{X}, \bs{\vse}) ]$ is
tractable. In additional to this approach of constructing the parameterized distribution $g(., \bs{\vse})$ and the deterministic function
$\Lm(\bs{\vse})$, the LR method opens up another possible avenue. In many situations, since the pdf or pmf $f(.)$ of $\bs{X}$ belongs to a
family of distributions parameterized by $\bs{\se} \in \Se$, the desired parameterized distribution $g(., \bs{\vse})$ can be directly obtained
by restricting the parameter space $\Se$ to its subset $\varTheta$. In many cases,  by appropriately choosing the subset $\varTheta$, the
desired deterministic function $\Lm(\bs{\vse})$ can obtained by the monotonicity or concavity of the likelihood ratio. In this way, the hassle
of using a weight function is eliminated.  Since this approach does not involve the weight function, the derivation of bounds for the
probability $\Pr \{ \bs{X} \in \mscr{E} \}$ can be extremely simple. Moreover, this approach may also be powerful in situations that the ME
method is ineffective.

In the proposed LR method, an important step is to minimize the deterministic function $\Lm (\bs{\vse})$ to obtain the tightest bound. We
realize that this minimization process is  relevant to the maximization of a likelihood function.  This implies an inherent connection between
the LR method for probabilistic inequalities and Fisher's concept of {\it maximum likelihood estimation} of parameters. In situations that the
maximum likelihood estimation is not mathematically convenient,  we propose to use Pearson's {\it method of moments} as an alternative technique
to obtain simple bounds for $\Pr \{ \bs{X} \in \mscr{E} \}$. Of course, the bounds obtained from the method of moments is relatively more
conservative than their counterparts obtained from the maximum likelihood estimation.

We have established a link between the LR method and the large deviation theory. Under mild conditions,  we have obtained explicit bounds for
rate functions, which are tight for a wide class of distributions.  We have applied the LR method to derive new concentration inequalities. In
particular, we have derived bounds for exponential families. We show that the LR approach and the Berry-Essen inequality can be used together to
derive tighter concentration inequalities.  We have established inequalities for multivariate generalized hypergeometric and inverse
hypergeometric distributions,  Dirichlet distribution and matrix gamma distribution.

The remainder of the paper is organized as follows.  In Section 2, we introduce the fundamentals of the LR method.  Specially, we propose the
general principle and establish the connection between the LR method and statistical parameter estimation.  In Section 3, we apply the LR method
to the development of concentration inequalities. In Section 4, we explore the connection between the LR method and the large deviation theory.
In Section 5, we establish a unified approach for deriving concentration inequalities for exponential families.  In particular, we demonstrate
that tightest possible concentration inequalities can be derived without using moment generating functions. In Section 6, we apply the LR method
to derive new concentration inequalities for the binomial distribution, hypergeometric distribution and generalized Poisson distribution.   In
Section 7, we apply the LR method to establish concentration inequalities for multivariate generalized hypergeometric and inverse hypergeometric
distributions, Dirichlet distribution and matrix gamma distribution.  Section 8 is the conclusion.  Most proofs are given in Appendices.

Throughout this paper, we shall use the following notations.  Let $\bb{R}$ denote the set of all real numbers.  Let $\bb{Z}$ denote the set of
all integers.  Let $\bb{Z}^+$ denote the set of all nonnegative integers. Let $\bb{N}$ denote the set of all positive integers. Let $\bb{I}_E$
denote the indicator function such that $\bb{I}_E = 1$ if $E$ is true and $\bb{I}_E = 0$ otherwise.   We use the notation $\bi{t}{k}$ to denote
a generalized combinatoric number in the sense that
\[
\bi{t}{k} = \f{\prod_{\ell = 1}^k (t -\ell + 1)}{k!} = \f{ \Ga(t + 1) }{\Ga(k + 1) \; \Ga(t - k + 1)},  \qqu  \bi{t}{0} = 1,
\]
where $t$ is a real number and $k$ is a non-negative integer.  We use $\ovl{X}_n$ to denote the average of random variables $X_1, \cd, X_n$,
that is, $\ovl{X}_n = \f{\sum_{i=1}^n X_i}{n}$.  The notation $\top$ denotes the transpose of a matrix.  The trace of a matrix is denoted by
$\tx{tr}$.  We use pdf and pmf to represent probability density function and probability mass function, respectively.  The notation
``$\Leftrightarrow$'' means ``if and only if''.  The other notations will be made clear as we proceed.

\sect{Likelihood Ratio Method}

In this section, we shall propose a new method for deriving probabilistic inequalities, referred to as the likelihood ratio (LR) method.

Let $E$ be an event which can be expressed in terms of matrix-valued random variable $\bs{X}$, where $\bs{X}$ is defined on the sample space
$\Om$ and $\si$-algebra $\mscr{F}$ such that the true probability measure is one of two measures $\Pr$ and $\bb{P}_{\bs{\vse} }$. Here, the
measure $\Pr$ is determined by pdf or pmf $f(.)$. The  measure $\bb{P}_{\bs{\vse} }$ is determined by pdf or pmf $g(., \bs{\vse})$, which is
parameterized by $\bs{\vse} \in \varTheta$.  The subscript in $\bb{P}_{\bs{\vse}}$ is used to indicate the dependence on the parameter
$\bs{\vse}$. Clearly, there exists a set, $\mscr{E}$, of deterministic matrices of the same size as $\bs{X}$ such that $E = \{ \bs{X} \in
\mscr{E} \}$. In many applications, it is desirable to obtain an upper bound for the probability $\Pr \{ E \}$.

\subsection{General Principle}

 The LR method is based on the
following general result. \beT \la{ThM888} Assume that there exists a function $\Lm(\bs{\vse})$ of $\bs{\vse} \in \varTheta$ such that \be
\la{LRBVIP} f (\bs{X}) \; \bb{I}_{ \{ \bs{X} \in \mscr{E} \} } \leq \Lm(\bs{\vse}) \; g (\bs{X}, \bs{\vse}) \qu \tx{for all $\bs{\vse} \in
\varTheta$}. \ee Then, \be \la{mainVIP}
 \Pr \{  E   \} \leq \inf_{ \bs{\vse} \in \varTheta }  \Lm(\bs{\vse})  \; \bb{P}_{\bs{\vse} } \{ E \}  \leq \inf_{ \bs{\vse} \in \varTheta } \Lm(\bs{\vse}). \ee  \eeT

The notion of the inequality in (\ref{LRBVIP}) is that $f (\bs{X} (\om) )  \; \bb{I}_{ \{ \bs{X} (\om)  \in \mscr{E} \} } \leq \Lm(\bs{\vse}) \;
g (\bs{X} (\om), \bs{\vse})$ for every $\om \in \Om$.   We call the function $\Lm(\bs{\vse})$ in (\ref{LRBVIP}) as {\it likelihood-ratio
bounding function}. Theorem \ref{ThM888} asserts that the probability of event $E$ is no greater than the likelihood ratio bounding function. In
a comparison of (\ref{mainVIP}) with (\ref{inviewvip}), it can be seen that the distribution of $\bs{X}$ is directly involved in
(\ref{mainVIP}), which indicates that the LR method allows for more direct and complete use of the distribution of $\bs{X}$.  As can be seen
from Theorem \ref{ThM888}, the key of the LR method is to find the likelihood ratio bounding function which is tight and amenable for
minimization as well.

To establish Theorem \ref{ThM888}, note that if $f(.)$ is a pdf, then it follows from (\ref{LRBVIP}) that  \bee \Pr \{ E \}   =   \int_{ \bs{x}
\in \mscr{E} } f (\bs{x}) \; d \bs{x}   \leq  \int_{ \bs{x} \in \mscr{E} }  \; \Lm(\bs{\vse})   \; g (\bs{x}, \bs{\vse} ) \; d \bs{x}   =
\Lm(\bs{\vse})   \; \int_{ \bs{x} \in \mscr{E} } \; g (\bs{x}, \bs{\vse} ) \; d \bs{x} \qu \tx{for all $\bs{\vse} \in \varTheta$}.  \eee
Similarly, if $f(.)$ is a pmf, then \bee \Pr \{ E \} =  \sum_{ \bs{x} \in \mscr{E} } f (\bs{x}) \leq \sum_{ \bs{x} \in \mscr{E} }  \;
\Lm(\bs{\vse}) \; g (\bs{x}, \bs{\vse} ) = \Lm(\bs{\vse}) \; \sum_{ \bs{x} \in \mscr{E} } \; g (\bs{x}, \bs{\vse} ) \; \qu \tx{for all
$\bs{\vse} \in \varTheta$}. \eee This implies that \be \la{VIPB} \Pr \{ E \} \leq \Lm(\bs{\vse}) \; \bb{P}_{\bs{\vse} } \{ E \} \leq
\Lm(\bs{\vse}) \qu \tx{for all $\bs{\vse} \in \varTheta$}. \ee   Since the inequalities in (\ref{VIPB}) hold for any $\bs{\vse} \in \varTheta$,
minimizing the bounds with respect to $\bs{\vse} \in \varTheta$ yields (\ref{mainVIP}).

\subsection{Construction of Parameterized Distributions}

In the sequel, we will discuss two approaches for constructing parameterized distributions $g(., \bs{\vse})$ that are required for the
application of the LR method.

\subsubsection{Weight Function} \la{weightpara}

A natural approach to construct parameterized distribution $g(., \bs{\vse})$ is to modify the pdf or pmf $f(.)$ by multiplying it with a
parameterized function and performing a normalization. Specifically, let $w (.,  \bs{\vse} )$ be a non-negative function with parameter
$\bs{\vse} \in \varTheta$ such that $\bb{E} [ w (\bs{X}, \bs{\se}) ]  < \iy$ for all $\bs{\vse} \in \varTheta$,  where the expectation is taken
under the probability measure $\Pr$ determined by $f(.)$. Define a family of distributions as
\[ g(\bs{x}, \bs{\vse}) =  \f{ w (\bs{x}, \bs{\vse}) \; f( \bs{x} )  }{\bb{E} [ w (\bs{X}, \bs{\vse}) ] }
\]
for  $\bs{\vse} \in \varTheta$ and $\bs{x}$ in the range of $\bs{X}$.  In view of its role in the modification of $f(.)$ as $g(., \bs{\vse})$,
the function $w (. , \bs{\vse})$ is called a {\it weight function}. Note that \be \la{nowgood} f(\bs{X}) \; w (\bs{X}, \bs{\vse})  = \bb{E} [ w
(\bs{X}, \bs{\vse}) ] \;  g(\bs{X}, \bs{\vse}) \qu \tx{for all $\bs{\vse} \in \varTheta$}. \ee For simplicity, we choose the weight function
such that the condition (\ref{inviewvip}) is satisfied. Combining (\ref{inviewvip}) and (\ref{nowgood}) yields  \[ f(\bs{X}) \; \bb{I}_{ \{
\bs{X} \in \mscr{E} \} }  \leq  f(\bs{X}) \; w (\bs{X}, \bs{\vse}) = \bb{E} [ w (\bs{X}, \bs{\vse}) ] \; g(\bs{X}, \bs{\vse}) \qu \tx{for all
$\bs{\vse} \in \varTheta$}. \]   Thus, the likelihood ratio bounding function can be taken as
\[
\Lm (\bs{\vse}) = \bb{E} [ w (\bs{X}, \bs{\vse})  ] \qu \tx{for $\bs{\vse} \in \varTheta$}.
\]
It follows from Theorem \ref{ThM888} that  \[ \Pr \{  E  \} \leq \inf_{ \bs{\vse} \in \varTheta }  \Lm (\bs{\vse}) \; \bb{P}_{\bs{\vse}} \{ E \}
\leq \inf_{ \bs{\vse} \in \varTheta } \Lm (\bs{\vse}).
\]
Thus, we have demonstrated that the ME method  is actually a special technique of the LR method.

\subsubsection{Parameter Restriction}

In Section \ref{weightpara}, we have discussed an approach for constructing the parameterized distribution $g(., \bs{\vse})$  based on weight
function. If the mathematical expectation of $\bb{E} [ w (\bs{X}, \bs{\se}) ]$ is not tractable, then this approach is not effective for
deriving tight probabilistic inequalities. To overcome this problem, we propose another approach, which is based on the idea of restricting
parameter space.

In many situations, the pdf or pmf $f(.)$ of $\bs{X}$ comes from a family of distributions parameterized by $\bs{\se} \in \Se$.  If so, then the
parameterized distribution $g(., \bs{\vse})$ can be taken as the subset of pdf or pmf with parameter $\bs{\vse}$ contained in a subset
$\varTheta$ of parameter space $\Se$.  By appropriately choosing the subset $\varTheta$, the deterministic function $\Lm(\bs{\vse})$ may be
readily obtained. To illustrate the idea, consider $\Pr \{ \ovl{X}_n \geq z \}$, where $z \in [\f{1}{2}, 1)$ and $\ovl{X}_n = \f{\sum_{i=1}^n
X_i}{n}$ with $X_1, \cd, X_n$ being i.i.d. random variables uniformly distributed over interval $[0, 1]$. Let $\bs{X} = [X_1, X_2, \cd,
X_n]^\top$ and $\bs{x} = [x_1, x_2, \cd, x_n]^\top$.  The pdf of $\bs{X}$ is easily seen to be
\[
f(\bs{x}) = \bec 1  & \tx{if} \; 0 \leq x_i \leq 1 \; \tx{for} \; i = 1, \cd, n,\\
0 & \tx{otherwise} \eec
\]
Note that the pdf $f(\bs{x})$ is contained in the family
\[
g (\bs{x}, \se) = \bec [ C(\se) ]^n \exp ( \se \sum_{i=1}^n x_i  )   & \tx{if} \; 0 \leq x_i \leq 1 \; \tx{for} \; i = 1, \cd, n,\\
0 & \tx{otherwise} \eec
\]
where  \[ C(\se) = \f{1}{ \int_0^1 e^{\se x} dx } = \bec \f{\se}{e^\se - 1}  & \tx{for} \; \se \neq 0,\\
1 & \tx{for} \; \se = 0 \eec
\]
with parameter $\se \in \bb{R}$, that is, the parameter space $\Se$ is the set of all real numbers.  Letting $\se = 0$ yields  $f(\bs{x}) = g
(\bs{x}, 0)$. If we restrict the parameter space as $\varTheta = (0, \iy)$, then it can be readily seen that \be \la{unifollow89}
 f(\bs{X}) \;
\bb{I}_{ \{ \ovl{X}_n \geq z \} } \leq  \Lm (\vse) \;  g (\bs{X}, \vse) \qqu \tx{for all $\vse \in \varTheta$}, \ee where
\[
\Lm (\vse) = \f{ [ C(0) ]^n \exp ( n  z  \times 0  )}{ [ C(\vse) ]^n \exp ( n \vse z  ) } = \li (  \f{ e^\vse - 1 }{ \vse e^{z \vse}  }  \ri
)^n.
\]
It follows from (\ref{unifollow89}) and Theorem \ref{ThM888} that
\[
\Pr \{ \ovl{X}_n \geq z \} \leq \inf_{\vse \in (0, \iy)} \li (  \f{ e^\vse - 1 }{ \vse e^{z \vse}  }  \ri )^n = \li (  \f{ e^\nu - 1 }{ \nu e^{z
\nu}  }  \ri )^n ,
\]
where $\nu > 0$ is the unique number satisfying $z = 1 + \f{1}{e^\nu - 1} - \f{1}{\nu}$.

\subsection{Bounding Probability of Order Relation}

Let $\bs{Z}$ be a random matrix which can be expressed as a  function of random matrix $\bs{X}$.  Let $\bs{z}$ be a deterministic matrix in the
range of $\bs{Z}$. Let ``$\bs{\prec}$'' denote a partial order relation on the range of $\bs{Z}$.  It is a frequent problem to obtain an upper
bound for the probability $\Pr \{ \bs{Z} \bs{\prec} \bs{z} \}$.  For this purpose,  we seek a nonnegative  function $\mscr{M} (. , \bs{\vse})$,
parameterized by $\bs{\vse} \in \varTheta$, which satisfies the following requirements:

(i) For arbitrary parametric value $\bs{\vse} \in \varTheta$, $\mscr{M}(\bs{u}, \bs{\vse}) \leq \mscr{M}(\bs{z}, \bs{\vse})$ if $\bs{u}
\bs{\prec}  \bs{z}$, where $\bs{u}$ is in the range of $\bs{Z}$.

(ii)  \be \la{aspineq}
 f  ( \bs{X} )  \; \bb{I}_{\{ \bs{Z} \bs{\prec} \bs{z} \}} \leq \mscr{M}( \bs{Z}, \bs{\vse} ) \;
 g ( \bs{X}, \bs{\vse}) \qu \tx{for all $\bs{\vse} \in \varTheta$}. \ee   As a consequence of (i) and (ii), we can write (\ref{aspineq}) as
\[
f  ( \bs{X} ) \; \bb{I}_{\{ \bs{Z} \bs{\prec} \bs{z} \}} \leq \mscr{M}( \bs{z}, \bs{\vse} ) \; g ( \bs{X}, \bs{\vse} ) \qu \tx{for all
$\bs{\vse} \in \varTheta$}.
\]
Applying the general principle described by Theorem \ref{ThM888}, we have
\[
\Pr \{ \bs{Z} \bs{\prec} \bs{z} \} \leq \inf_{ \bs{\vse} \in \varTheta}  \mscr{M} ( \bs{z}, \bs{\vse} ) \; \bb{P}_{\bs{\vse}} \{ \bs{Z}
\bs{\prec} \bs{z} \}  \leq \inf_{ \bs{\vse} \in \varTheta} \mscr{M} ( \bs{z}, \bs{\vse} ).
\]

It should be noted that in many situations, especially in scenarios where the underlying distribution belongs to an exponential family,  the
likelihood ratio $\f{ f ( \bs{X} ) }{ g ( \bs{X}, \bs{\vse} ) }$ can be expressed as a function $\mscr{M}( \bs{Z}, \bs{\vse} )$ which is
monotone in $\bs{Z}$. Such property is called {\it monotone likelihood ratio} property, which has been extensively explored in our paper
\cite{ChenLR} for deriving probabilistic inequalities.

\subsubsection{Optimal Bound with Maximum-Likelihood Estimate} \la{secMLE}

Under some conditions, a connection between the LR method and Fisher's concept of maximum-likelihood estimation can be established as follows.
Assume that \be \la{bydef} f  (\bs{X})  = \mscr{M}(\bs{Z}, \bs{\vse}) \; g (\bs{X}, \bs{\vse}) \qu \tx{for all $\bs{\vse} \in \varTheta$} \ee
and that the  maximum-likelihood estimator of $\bs{\vse} \in \varTheta$ exists  and can be expressed as a function $\varphi (\bs{Z})$ of
$\bs{Z}$. That is, \be \la{eqdef} \max_{\bs{\vse} \in \varTheta} g (\bs{X}, \bs{\vse}) =  g (\bs{X}, \varphi (\bs{Z}) ). \ee Making use of
(\ref{bydef}) and (\ref{eqdef}), we have \bee \min_{\bs{\vse} \in \varTheta} \mscr{M}(\bs{Z}, \bs{\vse})  = \f{ f ( \bs{X}) } { \max_{
\bs{\bs{\vse}} \in \varTheta} g ( \bs{X}, \bs{\vse} )  }   =  \f{ f  ( \bs{X} ) } { g ( \bs{X}, \varphi (\bs{Z})) } =
 \mscr{M}(\bs{Z}, \varphi(\bs{Z}) ). \eee Without loss of generality, assume that $\bs{z}$ is contained in the range of $\bs{Z}$. Hence, \be
\la{MLEVIP}
 \min_{\bs{\vse} \in \varTheta} \mscr{M} ( \bs{z}, \bs{\vse}) = \mscr{M} ( \bs{z}, \varphi(\bs{z}) )
\ee and it follows from Theorem \ref{ThM888} that
\[
\Pr \{ \bs{Z}  \bs{\prec} \bs{z} \} \leq \mscr{M} ( \bs{z}, \varphi( \bs{z} ) ) \; \bb{P}_{\varphi( \bs{z} )}  \{ \bs{Z} \bs{\prec} \bs{z} \}
\leq \mscr{M} ( \bs{z}, \varphi( \bs{z} ) ).
\]

\subsubsection{Approximately Optimal Bound with  Method of Moments}

As can be seen from Section \ref{secMLE}, to make the bound as tight as possible, the minimization problem becomes finding the
maximum-likelihood estimate for $\bs{\vse} \in \varTheta$.  In situations that the computation of the maximum-likelihood estimate is
inconvenient, other estimation techniques may be used.  In this direction, we propose to use Pearson's {\it method of moments}, which may offer
an estimate which is close to the maximum-likelihood estimate. Specifically, for $\bs{z}$ in the range of $\bs{Z}$, let $\bs{\xi}$ be a number
in $\varTheta$ such that $\bs{z} = \bb{E}_{\bs{\xi}} [ \bs{Z} ]$,  where the expectation is taken with the underlying pdf or pmf $g(.,
\bs{\xi})$. Hence, $\bs{\xi}$ is actually a function of $\bs{z}$. So, we write $\bs{\xi} = \xi (\bs{z})$.   In many cases, $\varphi(\bs{z}) =
\xi (\bs{z})$ or $\varphi(\bs{z}) \ap \xi (\bs{z})$.  If $\xi(\bs{z}) \in \varTheta$, then it follows from Theorem \ref{ThM888} that
\[
\Pr \{ \bs{Z}  \bs{\prec} \bs{z} \} \leq \mscr{M} ( \bs{z}, \xi(\bs{z}) ) \; \bb{P}_{\xi(\bs{z})} \{ \bs{Z} \bs{\prec} \bs{z} \} \leq \mscr{M}
(\bs{z}, \xi(\bs{z}) ).
\]

\section{Concentration Inequalities}

Let $\bs{X}_1, \bs{X}_2, \cd, \bs{X}_n$  and $\bs{X}$ be matrix-valued random variables defined on the same sample space $\Om$ and $\si$-algebra
$\mscr{F}$ such that the true probability measure is one of two measures $\Pr$ and $\bb{P}_{\bs{\vse}}$.  Here, the measure $\Pr$ is determined
by pdf or pmf $f(.)$. The  measure $\bb{P}_{\bs{\vse}}$ is determined by pdf or pmf $g(., \bs{\vse})$, which is parameterized by $\bs{\vse} \in
\varTheta$. Under each of these probability measures, the $\bs{X}$'s are independent  and identically distributed.  Define $\ovl{\bs{X}}_n = \f{
\sum_{i = 1}^n \bs{X}_i }{n}$ for $n \in \bb{N}$.  Let $\bs{\prec}$ denote a partial order relation defined on the range of $\ovl{\bs{X}}_n$.
Let $\mscr{X}$ denote the union of the ranges of $\ovl{\bs{X}}_n$ for all $n \in \bb{N}$.  In this setting,  we have the following results.

\beT

\la{con89666}

Assume that $\ro (\bs{x}, \bs{\vse} )$ is a concave function of $\bs{x} \in \mscr{X}$, with parameter $\bs{\vse} \in \varTheta$,  such that
$f(\bs{X}) \leq \exp ( \ro ( \bs{X}, \bs{\vse}) ) \; g(\bs{X}, \bs{\vse})$ for all $\bs{\vse} \in \mscr{\varTheta}$.   Let $\bs{z} \in \mscr{X}$
be a matrix such that for $\bs{x} \in \mscr{X}$ satisfying $\bs{x}  \bs{\prec} \bs{z}$, the inequality $\ro (\bs{x}, \bs{\vse}) \leq \ro
(\bs{z}, \bs{\vse})$ holds for all $\bs{\vse} \in \mscr{\varTheta}$. Then, \be \la{concentra8996} \Pr \{ \ovl{\bs{X}}_n \bs{\prec} \bs{z} \}
\leq \inf_{\bs{\vse} \in \varTheta } \exp ( n \ro (\bs{z}, \bs{\vse}) ) \; \bb{P}_{\bs{\vse}} \{ \ovl{\bs{X}}_n \bs{\prec} \bs{z} \} \leq
\inf_{\bs{\vse} \in \varTheta } \exp ( n \ro (\bs{z}, \bs{\vse}) ). \ee
  \eeT

See Appendix \ref{con89666app} for a proof.

In the particular case that $f(\bs{X}) = \exp ( \ro ( \bs{X}, \bs{\vse}) ) \; g(\bs{X}, \bs{\vse})$ for all $\bs{\vse} \in \mscr{\varTheta}$,
the connection between the LR method and statistical parameter estimation can be established as follows.

To achieve the optimal bound, consider the minimization of
\[
\exp ( \ro (\bs{z}, \bs{\vse}) ) = \f{ f(\bs{z}) }{ g(\bs{z}, \bs{\vse}) } \qqu \tx{subject to $\bs{\vse} \in \varTheta$}.
\]
Let $\varphi(\bs{X})$ be the maximum-likelihood estimator of $\bs{\vse} \in \varTheta$ such that
\[
\max_{\bs{\vse} \in \varTheta} g(\bs{X}, \bs{\vse}) = g(\bs{X}, \varphi(\bs{X}))
\]
and that the function $\varphi(.)$ is defined on the range of $\bs{X}$.  Assume that $g$ and $\varphi$ can be extended to the range of
$\ovl{\bs{X}}_n$ such that
\[
\max_{\bs{\vse} \in \varTheta} g( \ovl{\bs{X}}_n, \bs{\vse}) = g(\ovl{\bs{X}}_n, \varphi(\ovl{\bs{X}}_n)).
\]
Note that
\[
\min_{ \bs{\vse} \in \varTheta } \f{ f(\ovl{\bs{X}}_n) }{ g(\ovl{\bs{X}}_n, \bs{\vse}) } = \f{ f(\ovl{\bs{X}}_n) }{ \max_{ \bs{\vse} \in
\varTheta } g(\ovl{\bs{X}}_n, \bs{\vse}) } = \f{ f(\ovl{\bs{X}}_n) }{ g(\ovl{\bs{X}}_n, \varphi(\ovl{\bs{X}}_n)) }.
\]
Assume that $\bs{z}$ is contained in the range of $\ovl{\bs{X}}_n$.  Then,
\[
\min_{ \bs{\vse} \in \varTheta } \f{ f(\bs{z}) }{ g(\bs{z}, \bs{\vse}) } = \f{ f(\bs{z}) }{ g(\bs{z}, \varphi(\bs{z})) }
\]
and it follows from (\ref{concentra8996}) that
\[
\Pr \{ \ovl{\bs{X}}_n \bs{\prec} \bs{z} \} \leq \li [ \f{ f(\bs{z}) }{ g(\bs{z}, \varphi(\bs{z})) } \ri ]^n \; \bb{P}_{\varphi(\bs{z})} \{
\ovl{\bs{X}}_n \bs{\prec} \bs{z} \} \leq \li [ \f{ f(\bs{z}) }{ g(\bs{z}, \varphi(\bs{z})) } \ri ]^n.
\]
In situations that the computation of the maximum-likelihood estimate is inconvenient, the {\it method of moments} may offer an estimate which
is close to the maximum-likelihood estimate. Since there is only one parameter to be estimated, it suffices to use the first moment.
Specifically, let $\bs{\xi}$ be a number in $\varTheta$ such that $\bs{z} = \bb{E}_{\bs{\xi}} [ \bs{X} ]$, where the expectation is taken with
the underlying pdf or pmf  $g(., \bs{\xi})$.  Hence, $\bs{\xi}$ is actually a function of $\bs{z}$. So, we write $\bs{\xi} = \xi(\bs{z})$.   In
many cases,  $\varphi(\bs{z}) \ap \xi (\bs{z})$. If $\xi(\bs{z}) \in \varTheta$, then it follows from (\ref{concentra8996}) that
\[
\Pr \{ \ovl{\bs{X}}_n \bs{\prec} \bs{z} \} \leq \li [ \f{ f(\bs{z}) }{ g(\bs{z}, \xi(\bs{z})) } \ri ]^n \; \bb{P}_{\xi(\bs{z})} \{
\ovl{\bs{X}}_n \bs{\prec} \bs{z} \} \leq \li [ \f{ f(\bs{z}) }{ g(\bs{z}, \xi(\bs{z})) } \ri ]^n.
\]

\section{A Link between LR Method and Large Deviation Theory}

In this section, we will show that there exists a fundamental connection between the LR method and the theory of large deviations.

Let $X_1, X_2, \cd, X_n$ and $X$ be random variables defined on the same sample space $\Om$ and $\si$-algebra $\mscr{F}$ such that the true
probability measure is one of two measures $\Pr$ and $\bb{P}_g$.  Here, the measure $\Pr$ is determined by pdf or pmf $f(.)$. The  measure
$\bb{P}_g$ is determined by pdf or pmf $g(.)$.  Under each of these probability measures, the $X$'s are independent  and identically
distributed. Let $\bb{E}[.]$ and $\bb{E}_g [. ]$ denote the mathematical expectations taken under probability measures $\Pr$ and $\bb{P}_g$,
respectively.   Let $o(n)$ be a function of $n$ such that $\f{o(n)}{n} \to 0$ as $n \to \iy$.  In this setting,  we have the following results.

\beT \la{convip} Let $\mu = \bb{E} [X]$ and $z = \bb{E}_g [ X ]$.  Assume that $\bb{E}_g [ | X - z |^2 ] < \iy$ and  that $\bb{E} [ e^{s X} ] <
\iy$ for $s$ in a neighborhood of $0$. Let $\ro(.)$ be a continuous function.  The following assertions hold.

(I): If $z \geq \mu, \; \Pr \{X > z \} > 0$ and  there exists $\de
> 0$ such that
\be \la{cona88}
 \li [ \prod_{i=1}^n f ( X_i ) \ri ] \bb{I}_{ \{ z \leq \ovl{X}_n < z + \de \} } \leq \exp \li ( n \ro ( \ovl{X}_n ) + o(n) \ri
) \prod_{i=1}^n g ( X_i ) \ee  for sufficiently large $n > 0$,  then $\Pr \{ \ovl{X}_n \geq z \} \leq  \exp (n \ro ( z ) )$ and $\lim_{n \to
\iy} \f{1}{n} \ln \Pr \{ \ovl{X}_n \geq z \}  \leq \ro ( z )$.

(II): If $z \geq \mu, \; \Pr \{X > z \} > 0$ and there exists $\de > 0$  such that (\ref{cona88}) and the following inequality \be \la{conb88}
 \li [ \prod_{i=1}^n g ( X_i ) \ri ] \bb{I}_{ \{ z \leq \ovl{X}_n < z + \de \} } \leq \exp \li ( - n \ro (
\ovl{X}_n ) + o(n) \ri ) \prod_{i=1}^n f ( X_i ) \ee hold for sufficiently large  $n > 0$, then  $\lim_{n \to \iy} \f{1}{n} \ln \Pr \{ \ovl{X}_n
\geq z \} = \ro ( z )$.

(III): If $z \leq \mu, \; \Pr \{X < z \} > 0$ and  there exists $\de > 0$ such that \be \la{cona88rev}
 \li [ \prod_{i=1}^n f ( X_i ) \ri ] \bb{I}_{ \{ z \geq \ovl{X}_n > z - \de \} } \leq \exp \li ( n \ro ( \ovl{X}_n ) + o(n) \ri
) \prod_{i=1}^n g ( X_i ) \ee  for sufficiently large $n > 0$,  then $\Pr \{ \ovl{X}_n \leq z \} \leq  \exp (n \ro ( z ) )$ and $\lim_{n \to
\iy} \f{1}{n} \ln \Pr \{ \ovl{X}_n \leq z \}  \leq \ro ( z )$.

(IV): If $z \leq \mu, \; \Pr \{X < z \} > 0$ and there exist $\de > 0$ such that (\ref{cona88rev}) and the following inequality \be
\la{conb88rev}
 \li [ \prod_{i=1}^n g ( X_i ) \ri ] \bb{I}_{ \{ z \geq \ovl{X}_n > z - \de \} } \leq \exp \li ( - n \ro (
\ovl{X}_n ) + o(n) \ri ) \prod_{i=1}^n f ( X_i ) \ee hold for sufficiently large $n > 0$,  then,  $\lim_{n \to \iy} \f{1}{n} \ln \Pr \{
\ovl{X}_n \leq z \} = \ro ( z )$.

(V): If $\ro (.)$ is a concave continuous function such that $ f(X) \leq \exp (\ro (X) ) \; g(X)$,  then \bee &  & \Pr \{ \ovl{X}_n \geq z \}
\leq \exp (n  \ro (z) )
\qqu \tx{provided that $z \geq \mu$ and $\Pr \{ X > z \} > 0$}, \\
&  & \Pr \{ \ovl{X}_n \leq z \} \leq  \exp ( n \ro (z) ) \qqu \tx{provided that $z \leq \mu$ and $\Pr \{ X < z \} > 0$}. \eee

(VI): If $\ro (.)$ is a linear continuous function such that $ f(X) = \exp (\ro (X) ) \; g(X)$, then \bee & & \lim_{n \to \iy} \f{1}{n} \ln \Pr
\{ \ovl{X}_n \geq z \}   = \ro ( z )
\qqu \tx{provided that $z \geq \mu$ and $\Pr \{X > z \} > 0$}, \\
&  & \lim_{n \to \iy} \f{1}{n} \ln \Pr \{ \ovl{X}_n \leq z \}   = \ro ( z ) \qqu \tx{provided that $z \leq \mu$ and $\Pr \{X < z \}
> 0$}. \eee

 \eeT

See Appendix \ref{convipapp} for a proof.  The assertions (V) and (VI) indicate that, under certain conditions, the probabilistic inequalities
obtained from the LR method coincide with the Chernoff bounds derived from moment generating functions.

\section{A Unified Theory for Exponential Families}

Our main objective for this section is to develop a unified theory for bounding the tail probabilities of exponential families.   A multivariate
exponential family has pdf or pmf of the form \be \la{vecexp}
 f (\bs{x}, \bs{\se}) = v (\bs{x}) \exp \li ( \bs{\eta}
(\bs{\se}) *  \bs{u} ( \bs{x} )  - \ze (\bs{\se}) \ri ), \qu \bs{\se} \in \Se, \ee where $\bs{x}$ and $\bs{\se}$ are column vectors, $\bs{\eta}
(\bs{\se})$ is column vector function of $\bs{\se}$, $\bs{u} ( \bs{x} )$ is a column vector function of $\bs{x}$, $\ze (\bs{\se})$ is a scalar
function of  $\bs{\se}$, and $v ( \bs{x} )$ is a scalar function of  $\bs{x}$.  The symbol ``$*$'' denotes the inner product. Assume that
$\bs{\eta} (\bs{\se})$ and $\ze (\bs{\se})$ are differentiable functions of $\bs{\se} \in \Se$. Let $\bs{X}$ be a random vector having a pdf or
pmf in the exponential family (\ref{vecexp}). It can be shown  that
\[
 \bs{\eta}^\prime (\bs{\se}  ) \; \bb{E}_{\bs{\se}} [ u ( \bs{X} ) ] = \ze^\prime (\bs{\se}  ),
\]
where $\bs{\eta}^\prime (\bs{\se}  )$ is the derivative matrix such that the element in the $i$th row and $j$th column is equal to $\f{\pa
\bs{\eta}_j (\bs{\se} )}{\pa \se_i}$.  The subscript in $\bb{E}_{\bs{\se}}$ indicates that the expectation is taken under the probability
measure determined by a pdf or pmf with parameter $\bs{\se}$.   Let $\Pr \{ E \mid \bs{\se} \}$ denote the probability of event $E$, where the
probability measure is determined by a pdf or pmf with parameter $\bs{\se}$.  For two vectors $\bs{V}_1$ and $\bs{V}_2$, we write $\bs{V}_1 \leq
\bs{V}_2$ if each element of $\bs{V}_1$ is no greater than the corresponding element of $\bs{V}_2$. Similarly, we write $\bs{V}_1 \geq \bs{V}_2$
if each element of $\bs{V}_1$ is no less than the corresponding element of $\bs{V}_2$.  In this setting, we have the following results.

\beT

\la{unifiedexp}

Let $\bs{X}_1, \bs{X}_2, \cd, \bs{X}_n$ be i.i.d. random vectors with common pdf or pmf parameterized by $\bs{\se} \in \Se$ as (\ref{vecexp}).
Define $\bs{Z} = \f{ 1 }{n} \sum_{i=1}^n \bs{u} (\bs{X}_i)$ and $\ro (\bs{z}, \bs{\vse}) =  [ \bs{\eta} (\bs{\se}) - \bs{\eta} (\bs{\vse}) ] *
\bs{z} - \ze (\bs{\se}) + \ze (\bs{\vse})$, where $\bs{\se}, \bs{\vse} \in \Se$ and $\bs{z}$ is a vector of the same size as $\bs{Z}$. Define
$\mscr{A}_{\bs{\se}} = \{ \bs{\vse} \in \Se: \bs{\eta} (\bs{\vse}) \leq \bs{\eta} (\bs{\se}) \}$ and $\mscr{B}_{\bs{\se}} = \{ \bs{\vse} \in
\Se: \bs{\eta} (\bs{\vse}) \geq \bs{\eta} (\bs{\se}) \}$. Then, \bel  &  & \Pr \{ \bs{Z} \leq \bs{z} \mid \bs{\se} \} \leq \inf_{ \bs{\vse} \in
\mscr{A}_{\bs{\se}} } \exp( n \ro (\bs{z}, \bs{\vse}) )  \; \Pr \{ \bs{Z} \leq \bs{z} \mid \bs{\vse} \} \leq \inf_{ \bs{\vse} \in
\mscr{A}_{\bs{\se}} } \exp( n \ro (\bs{z}, \bs{\vse}) ), \la{886633a}\\
&  & \Pr \{ \bs{Z} \geq \bs{z} \mid \bs{\se} \} \leq \inf_{ \bs{\vse} \in  \mscr{B}_{\bs{\se}} } \exp( n \ro (\bs{z}, \bs{\vse}) ) \; \Pr \{
\bs{Z} \geq \bs{z} \mid \bs{\vse} \}  \leq \inf_{ \bs{\vse} \in  \mscr{B}_{\bs{\se}} } \exp( n \ro (\bs{z}, \bs{\vse})). \la{886633b} \eel
 \eeT

It should be noted that Theorem \ref{unifiedexp} involves no moment generating function in the derivation of probabilistic inequalities.  To
find explicit expressions for the bounds, the method of moments may be useful. Specifically, if a vector $\bs{\al} \in \mscr{A}_{\bs{\se}}$ can
be obtained such that  $\bb{E}_{\bs{\al}} [ u( \bs{X} ) ] = \bs{z}$ or equivalently, $\bs{\eta}^\prime (\bs{\al} ) \; \bs{z} = \ze^\prime
(\bs{\al} )$, then it follows from (\ref{886633a}) that
\[
\Pr \{ \bs{Z} \leq \bs{z} \mid \bs{\se} \} \leq  \exp( n \ro (\bs{z}, \bs{\al}) ) \;  \Pr \{ \bs{Z} \leq \bs{z} \mid \bs{\al} \} \leq \exp( n
\ro (\bs{z}, \bs{\al}) ),
\]
Similarly, if a vector $\bs{\al} \in \mscr{B}_{\bs{\se}}$ can be obtained such that $\bs{\eta}^\prime (\bs{\al} ) \; \bs{z} = \ze^\prime
(\bs{\al} )$, then it follows from (\ref{886633b}) that
\[ \Pr \{ \bs{Z} \geq \bs{z} \mid \bs{\se} \} \leq  \exp( n \ro (\bs{z}, \bs{\al}) ) \; \Pr \{ \bs{Z} \geq \bs{z} \mid \bs{\al} \}
\leq \exp( n \ro (\bs{z}, \bs{\al}) ).
\]
To prove Theorem \ref{unifiedexp}, it suffices to make use of Theorem \ref{ThM888} and the observation that \bee &  &  \li [ \prod_{i=1}^n f(
\bs{X}_i, \bs{\se} ) \ri ] \bb{I}_{ \{ \bs{Z} \leq \bs{z} \} } \leq \exp( n \ro (\bs{z}, \bs{\vse}) ) \li [ \prod_{i=1}^n f ( \bs{X}_i,
\bs{\vse} ) \ri ] \qqu \tx{for} \; \bs{\vse} \in \mscr{A}_{\bs{\se}};\\
&  & \li [ \prod_{i=1}^n f( \bs{X}_i, \bs{\se}  ) \ri ] \bb{I}_{ \{ \bs{Z} \geq \bs{z} \} } \leq \exp( n \ro (\bs{z}, \bs{\vse}) ) \li [
\prod_{i=1}^n f( \bs{X}_i, \bs{\vse} ) \ri ] \qqu \tx{for} \; \bs{\vse} \in \mscr{B}_{\bs{\se}}. \eee

\bsk

In the above discussion, we have developed a general method for deriving probabilistic inequalities for the multivariate case.   In the
remainder of this section, we will discuss the univariate case. A univariate exponential family is a set of probability distributions whose pdf
or pmf  can be expressed in the form \be \la{expdef}
 f (x, \se) = v (x) \exp ( \eta (\se) u ( x ) - \ze (\se) ), \qqu \se
 \in \Se
\ee where $u (x),  v(x),  \eta(\se)$, and $\ze(\se)$ are known functions.  Assume that $\eta (\se)$ and $\ze (\se)$ are differentiable functions
of $\se \in \Se$.  Let $X$ be a random variable having a pdf or pmf in the exponential family (\ref{expdef}). Then,
\[
 \eta^\prime (\se) \; \bb{E}_\se [  u (X) ] =  \ze^\prime (\se), \qqu \se \in \Se.
\]
Clearly, the expectation of $u (X)$ is a function of $\se$.  Let this function be denoted by $\mu (\se)$.

For the exponential family described above, we have the following results.

\beT

\la{univexp}
 Let $X$ be a random variable with pdf or pmf parameterized by $\se \in \Theta$ as (\ref{expdef}).
Define $Z = \f{ 1 } { n } \sum_{i = 1}^n u (X_i)$, where $X_1, \cd, X_n$ are i.i.d. samples of $X$.  Let $\al \in \Se$ and $z = \mu ( \al )$.
Define $\ro (z, \vse) = [ \eta (\se) - \eta (\vse) ] z - \ze (\se) + \ze (\vse)$ for $\se, \vse \in \Se$.  Let $\mscr{C}$ denote the absolute
constant in the Berry-Essen inequality.  Define
\[ \vDe = \min \li \{ \f{1}{2}, \; \f{\mscr{C}}{\sq{n}} \f{ \bb{E}_\al [ |u (X) - z |^3 ] } { \bb{E}_\al^{\f{3}{2}} [ |u (X) - z |^2 ] } \ri \},
\]
where the expectation is taken under the probability measure determined by the pdf or pmf $f(x, \al)$. The following assertions hold.

(I): If $\mu (\se) \leq z$ and $\Pr \{ u(X) < z \mid \se \} > 0$, then \bee \Pr \li \{ Z \geq z \mid \se \ri \} \leq  \exp ( n \ro (z, \al) ),
\qqu \qqu \lim_{n \to \iy} \f{1}{n} \Pr \li \{ Z \geq z \mid \se \ri \} = \ro (z, \al). \eee

(II):  If $\mu (\se) \geq z$ and $\Pr \{ u(X) > z \mid \se \} > 0$, then \bee  \Pr \li \{ Z \leq z \mid \se \ri \} \leq \exp ( n \ro (z, \al) ),
\qqu \qqu \lim_{n \to \iy} \f{1}{n} \Pr \li \{ Z \leq z \mid \se \ri \} = \ro (z, \al). \eee

(III): If $\eta (\se) \leq \eta (\al)$, then $\Pr \li \{ Z \geq z \mid \se \ri \} \leq  \exp ( n \ro (z, \al) ) \times \Pr \li \{ Z \geq z \mid
\al \ri \} \leq \li (   \f{1}{2} + \vDe \ri ) \exp ( n \ro (z, \al) )$.

(IV): If $\eta (\se) \geq \eta (\al)$, then $\Pr \li \{ Z \leq z \mid \se \ri \} \leq  \exp ( n \ro (z, \al ) ) \times \Pr \li \{ Z \leq z \mid
\al \ri \} \leq \li (   \f{1}{2} + \vDe \ri ) \exp ( n \ro (z, \al) )$.

\eeT

The assertions (I) and (II) of Theorem \ref{univexp} follow immediately from the assertions (V) and (VI) of Theorem \ref{convip}.  The
assertions (III) and (IV) of Theorem \ref{univexp} are direct consequences of Theorem \ref{unifiedexp} and the Berry-Essen inequality.  The
famous Berry-Essen inequality \cite{Berry, Essen} asserts the following:

    Let $Y_1, Y_2, ...$ be i.i.d. samples of random variable $Y$ such that
    $\bb{E}[Y] = 0, \; \bb{E}[Y^2] > 0$, and $\bb{E}[|Y|^3] < \iy$. Also, let
$F_n$ be the cumulative distribution function (cdf) of ${\sum_{i = 1}^n Y_i \over {\sqrt{n \bb{E}[Y^2]} }}$, and $\Phi$ the cdf of the standard
normal distribution. Then, there exists a positive constant $\mscr{C}$ such that for all $y$ and $n$,
\[
\left|F_n(y) - \Phi(y)\right| \le \f{\mscr{C}}{\sq{n}} {\bb{E}[|Y|^3] \over \bb{E}^{3\sh 2} [Y^2]}.
\]
Recently, Tyurin \cite{Tyurin} has shown that $\mscr{C} < 0.4785 < \f{1}{2}$.

\section{Univariate Probabilistic Inequalities}

In this section, we shall apply the LR method to derive inequalities for some important univariate random variables.

\subsection{Bernoulli Distribution}

A Bernoulli random variable, $X$, of mean value $p \in (0, 1)$ has a pmf
\[
f(x, p) = \Pr\{ X = x \} = p^x (1 - p)^{1 - x}, \qqu x \in \{0, 1\},
\]
which belongs to the exponential family $f(x, p) = v (x) \exp \left( \eta (p) u (x) - \ze (p) \right), \; x \in \{0, 1\}, \; p \in (0, 1)$,
where
\[
u (x) = x,  \qqu v (x) = 1, \qqu  \eta (p) = \ln {p \over 1-p}, \qqu \ze (p) = \ln \f{1}{1-p}.
\]
Note that
\[
\bb{E} [X] = p, \qqu \bb{E} [ | X - p |^2 ] = p (1- p), \qqu \bb{E} [ |X - p |^3 ]  = [ p^2 + (1 - p)^2 ] p(1 - p)
\]
and
\[
\f{ \bb{E} [ |X - p |^3 ] } { \bb{E} ^{\f{3}{2}} [ | X - p |^2 ] } = \f{p^2 + (1 - p)^2 }{ \sq{ p (1 - p) } }.
\]
Making use of these facts,  the assertions (III) and (IV) of Theorem \ref{univexp},  we have the following results.

\beT \la{Bernui}
  Let $X_1, \cd, X_n$ be i.i.d. samples of Bernoulli random variable $X$ of mean value $p
\in (0, 1)$.  Define $\mscr{M} ( z, p) = \li ( \f{p}{z} \ri )^z  \li ( \f{ 1 - p }{ 1 - z} \ri )^{ 1 - z}$ for $z \in (0, 1)$ and $p \in (0,
1)$. Then, \bee &  & \Pr \li \{ \ovl{X}_n \leq z \ri \} \leq \li ( \f{1}{2} + \vDe \ri ) [ \mscr{M} ( z, p) ]^n \qqu \tx{for $z \in (0, p)$},\\
& & \Pr \li \{ \ovl{X}_n \geq z \ri \} \leq \li ( \f{1}{2} + \vDe \ri ) [ \mscr{M} ( z, p) ]^n \qqu \tx{for $z \in (p, 1)$},
 \eee where $\vDe =
\min \li \{ \f{1}{2}, \f{\mscr{C} [z^2 + (1 - z)^2] }{ \sq{ n z (1 - z) } } \ri \}$ with $\mscr{C}$ being the Berry-Essen constant.

\eeT

Clearly,  $\vDe \to 0$ as $n \to \iy$.  For moderate to large $n$, the factor $\f{1}{2} + \vDe$ is close to $\f{1}{2}$. This implies that
Theorem \ref{Bernui} has reduced the existing Chernoff bounds for the tail probabilities of a binomial distribution by a factor of $\f{1}{2}$.
This improvement is due to the application of the LR method.  Moreover, the LR method offers an extremely proof for the existing Chernoff bounds
in the case of a binomial distribution. Specifically, making use of Theorem \ref{ThM888} and the observation that
\[ \li [ \prod_{i=1}^n f( X_i, p ) \ri ] \bb{I}_{ \{ \ovl{X}_n \leq z \} } \leq [ \mscr{M} ( z, p) ]^n  \li [
\prod_{i=1}^n f ( X_i, z ) \ri ] \qqu \tx{for} \; z \in (0, p),
\]
we have $\Pr \li \{ \ovl{X}_n \leq z \ri \} \leq [\mscr{M} ( z, p) ]^n$ for $z \in (0, p)$.  On the other hand, making use of Theorem
\ref{ThM888} and the observation that
\[
\li [ \prod_{i=1}^n f( X_i, p  ) \ri ] \bb{I}_{ \{ \ovl{X}_n \geq z \} } \leq [ \mscr{M} ( z, p) ]^n  \li [ \prod_{i=1}^n f( X_i, z ) \ri ] \qqu
\tx{for} \; z \in (p, 1),
\]
we have $\Pr \li \{ \ovl{X}_n \geq z \ri \} \leq [ \mscr{M} ( z, p) ]^n$ for $z \in (p, 1)$.  In sharp contrast to the classical method, such
arguments involve no moment generating function.

\subsection{Hypergeometric Distribution}

The hypergeometric distribution can be described by the following model.  Consider a finite population of $N$ units, of which there are $R$
units having a certain attribute.  Draw $n$ units from the whole population by sampling without replacement. Let $X$ denote the number of units
having the attribute found in the $n$ draws.  Then, $X$ is a random variable possessing a hypergeometric distribution such that \be
\la{hyperCDF}
 \Pr \{ X = x \} = \bec \f{ \bi{R}{x} \bi{B}{n - x} } { \bi{N}{n}  } & \tx{for} \; x \in \{0, 1, \cd, n \} \; \tx{such that} \; x \leq R \; \tx{and} \; n - x \leq
 B,\\
0  & \tx{for} \; x \in \{0, 1, \cd, n \} \; \tx{such that} \; x > R \; \tx{or} \; n - x > B, \eec \ee where $B = N - R$.  The mean of $X$ is
$\mu = \bb{E} [ X ] = \f{n R}{N}$.   Since the moment generating function of the hypergeometric distribution is intractable, it is difficult to
derive tight inequalities for the tail distributions based on Chernoff bounds.   However, such difficulties can be overcame by the LR method. In
this direction, we have the following results.

\beT

\la{Hyperfirst}

Let $X$ be a random variable possessing a hypergeometric distribution  defined by (\ref{hyperCDF}).  Let $r$ and $b$ be nonnegative integers
such that $r \leq R, \; b \leq B$ and $r + b = n$.  Define $\wh{R} = \min \{N, \lf (N + 1) \f{r}{n} \rf \}$ and $\wh{B} = N - \wh{R}$.   Then,
\be \la{IneqChena} \Pr \{ X \leq r \} \leq \f{ \bi{ R }{r} \bi{B } {b} } { \bi{ \wh{R} }{r} \bi{ \wh{B} } {b}  } \qu \tx{for $r \leq \mu$}, \ee
and \be \la{IneqChenb} \Pr \{ X \geq r \} \leq \f{ \bi{ R }{r} \bi{ B } {b} } { \bi{ \wh{R} }{r} \bi{ \wh{B} } {b}  }   \qu \tx{for $r \geq
\mu$}. \ee

 \eeT

 See Appendix \ref{Hyperfirstapp} for proof.  Actually, in August 2010, Chen had derived inequalities (\ref{IneqChena}) and (\ref{IneqChenb})
  based on the LR method and applied them  in \cite[page 54, inequalities (35), (36)]{Chenestimation} and \cite[page 20, inequalities (15), (16) ]{Chentest}
   for developing multistage sampling schemes for statistical inference of the population proportion $p$.

\subsection{Generalized Poisson Distribution}

A random variable $X$  is said to possess a generalized Poisson distribution if \be \la{genPos8996388} \Pr \{ X = x \} = \f{ \lm (\lm + x
\al)^{x - 1} e^{-\lm - x \al}  }{ x! }, \qqu x = 0, 1, 2, \cd, \ee where $\lm > 0$ and $0 \leq \al < 1$.   This distribution is sometimes called
Consul's generalized Poisson distribution (see, \cite{Consul} and the references therein).  It can be shown that the expectation of $X$ is
$\bb{E} [ X ] = \f{\lm}{1 - \al}$. For $\al \in (0, 1)$, the moment generating function of $X$ is not tractable for deriving Chernoff bounds for
the distribution. However, the LR method can be applied to derive tight bounds as follows.

\beT

\la{ThmPos} Let $X_1, \cd, X_n$ be i.i.d. random variables possessing the generalized Poisson distribution described by (\ref{genPos8996388}).
Then, \bel & & \Pr \{ \ovl{X}_n \leq z \} \leq  \f{\lm }{ (1 - \al) (\lm + z \al )} \li [ \li (  \f{ \lm }{ z } + \al \ri )^z
\f{ e^{(1 -\al) z}  }{e^\lm }  \ri ]^{n}  \qu \tx{provided that $0 < z \leq \f{\lm}{1 - \al }$}, \la{Pos8899a}\\
&  & \Pr \{ \ovl{X}_n \geq z \} \leq  \f{\lm }{ (1 - \al) (\lm + z \al )} \li [ \li (  \f{ \lm }{ z } + \al \ri )^z \f{ e^{(1 -\al) z}  }{e^\lm
}  \ri ]^{n} \qu \tx{provided that $z \geq \f{\lm}{1 - \al }$}. \la{Pos8899b} \eel Moreover, \bel &  & \Pr \{ \ovl{X}_n \leq z \} \leq \f{\lm
(\nu + z \al)}{\nu (\lm + z \al )} \li [ \li (  \f{ \lm + z \al }{ \nu + z \al } \ri )^z
\f{ e^\nu  }{e^\lm }  \ri ]^{n}  \qu \tx{provided that $0 < z \leq \f{\lm}{1 - \al + \f{\al}{n \lm} }$}, \la{Pos8899c}\\
&  & \Pr \{ \ovl{X}_n \geq z \} \leq  \f{\lm (\nu + z \al)}{\nu (\lm + z \al )} \li [ \li (  \f{ \lm + z \al }{ \nu + z \al } \ri )^z \f{ e^\nu
}{ e^\lm }  \ri ]^{n}  \qu \tx{provided that $z \geq \f{\lm}{1 - \al + \f{\al}{n \lm} }$}, \la{Pos8899d} \eel  where $\nu = \f{1}{2} [ (1 - \al)
z + \sq{ (1 - \al)^2 z^2 + 4 z \al \sh n}  ]$.  Furthermore, if $\al = 0$, then \bel &  & \Pr \li \{ \ovl{X}_n \geq z \ri \} \leq \li ( \f{1}{2}
+ \vDe \ri ) \li ( \f{ \lm^z e^{z} } {
z^z e^{\lm} } \ri )^n \qqu \tx{for $z \geq \lm$}, \la{Pos88a}\\
&  & \Pr \li \{ \ovl{X}_n \leq z \ri \} \leq \li ( \f{1}{2} + \vDe \ri ) \li ( \f{ \lm^z e^{z} } { z^z e^{\lm} } \ri )^n \qqu \tx{for $0 < z
\leq \lm$}, \la{Pos88b} \eel where $\vDe = \min \li \{ \f{1}{2}, \f{ \mscr{C} } { \sq{n}  }  \li ( 3  +  \f{1}{z} \ri )^{ \f{3}{4}} \ri \}$ with
$\mscr{C}$ being the Berry-Essen constant.

\eeT

See Appendix \ref{ThmPosapp} for a proof.

\subsection{Gamma Distribution}

In probability theory and statistics, a random variable $X$ is said to have a gamma distribution if its pdf is of the form \be \la{gapdf}
 f(x) = \frac{x^{k - 1}} { \Gamma(k) \se ^{ k }    }  \exp \li ( - \frac{x}{\se} \ri ) \;\;\; {\rm for} \;\;\; 0 < x < \infty,
\ee where $\se > 0, \; k > 0$ are referred to as the scale parameter and shape parameter respectively.  By letting \[ u (x) = \f{x}{k}, \qqu v
(x) = \f{ x^{k - 1} } { \Ga(k)  }, \qqu \eta (\se) = - \f{k}{\se}, \qqu \ze ( \se ) = k \ln \se,
\]
we can write $f(x) = v (x) \exp \left( \eta (\se) u (x) - \ze (\se) \right)$, which can be seen to be an exponential family. The moment
generating function of $X$ is $\phi (s) = \bb{E} [ e^{s X} ] = (1 - \se s)^{-k}$ for $s < \f{1}{\se}$. It can be shown by induction that
\[
\f{ d^{\ell + 1} \phi (s)  } {  d s^{\ell + 1} }  = \f{( k + \ell) \se}{1 - \se s} \f{ d^\ell \phi (s) } {  d s^\ell }, \qqu \bb{E} [ X^{\ell +
1} ] = (k + \ell) \se \li. \f{ d^\ell \phi (s) } {  d s^\ell } \ri |_{s = 0} = \se^{\ell + 1} \prod_{i = 0}^\ell (k + i)
\]
for $\ell = 0, 1, 2, \cd$. Therefore, \[ \bb{E} [ | X - k \se |^2 ] = k \se^2, \qqu  \bb{E} [ | X - k \se |^4 ] = \sum_{i = 0}^4 \bi{4}{i}  (- k
\se)^i \bb{E} [ X^{4 - i}] = 3 k (k + 2) \se^4, \] and by Lyapunov's inequality,  \[ \f{  \bb{E} [ | X - k \se |^3 ]  } { \bb{E}^{ \f{3}{2} } [
| X - k \se |^2 ] } < \f{ \bb{E}^{ \f{3}{4} } [ | X - k \se |^4 ]  } { \bb{E}^{ \f{3}{2} } [ | X - k \se |^2 ]  } = \li ( 3 + \f{6}{k} \ri )^{
\f{3}{4} }.
\] Making use of these facts,  the assertions (III) and (IV) of Theorem \ref{univexp},  we have the following results.

\beT

\la{gauniinq}

Let $X_1, \cd, X_n$ be i.i.d. samples of Gamma random variable $X$ possessing a probability density function (\ref{gapdf}). Then, \bee & & \Pr
\li \{ \ovl{X}_n \leq \vro k \se \ri \} \leq \li ( \f{1}{2} + \vDe \ri ) \li [ \vro  \exp \li ( 1 - \vro \ri ) \ri ]^{k n} \qqu \tx{for $0 <
\vro \leq
1$},\\
& & \Pr \li \{ \ovl{X}_n \geq \vro k \se \ri \} \leq \li ( \f{1}{2} + \vDe \ri ) \li [ \vro \exp \li ( 1 - \vro \ri ) \ri ]^{k n} \qqu \tx{for
$\vro \geq 1$},
 \eee where $\vDe = \min \li \{ \f{1}{2},  \f{ \mscr{C} } { \sq{n} } \li ( 3 + \f{6}{k} \ri )^{ \f{3}{4} }  \ri \}$ with $\mscr{C}$ being the Berry-Essen constant.

\eeT

\section{Multivariate Probabilistic Inequalities}

In this section, we shall apply the LR method to derive inequalities for multivariate random variables.  We need to define two partial order
relations $\bs{\prec}$ and $\bs{\succ}$ for vectors as follows.

Let $\ka \in \bb{N}$.  For vectors  $\bs{x} = [x_0, x_1, \cd, x_\ka]^\top$ and $\bs{y} = [y_0, y_1, \cd, y_\ka]^\top$, we write $\bs{x}
\bs{\prec} \bs{y}$ if
\[
x_i \leq y_i \qu \tx{for} \; i = 1, \cd, \ka.
\]
Similarly, we write $\bs{x} \bs{\succ} \bs{y}$ if
\[
x_i \geq y_i \qu \tx{for} \; i = 1, \cd, \ka.
\]
These partial order relations will be used in Sections 7.1, 7.2 and 7.3.

\subsection{Multivariate Generalized Hypergeometric Distribution}

 In probability theory, random variables $X_0, X_1, \cd, X_\ka$ are said to possess a multivariate hypergeometric
distribution if  \[
\Pr \{ X_i = x_i, \; i = 0, 1, \cd, \ka \}  = \bec 0 & \tx{if $x_i \geq 1 + C_i$ for some $i \in \{ 0, 1, \cd, \ka \}$},\\
\f{ \prod_{i=0}^k \bi{  C_i }{x_i}   }{  \bi{N} {n} } & \tx{otherwise} \eec
\]
where $C_i \in \bb{N}, \; x_i \in \bb{Z}^+$ for $i = 0, 1, \cd, \ka$ and
\[
N = \sum_{i = 0}^\ka C_i  \geq \sum_{i = 0}^\ka x_i = n.
\]
Actually, under mild restrictions, the multivariate hypergeometric distribution can be generalized by allowing $N$ and $C_i$ to be real numbers.
More formally, random variables $X_0, X_1, \cd, X_n$ are said to possess a {\it multivariate generalized  hypergeometric distribution} if  \be
\la{GINHY}
\Pr \{ X_i = x_i, \; i = 0, 1, \cd, \ka \} = \bec 0 & \tx{if $x_i \geq 1 + \mcal{C}_i > 1$ for some $i \in \{ 0, 1, \cd, k \}$},\\
\f{ \prod_{i=0}^k \bi{  \mcal{C}_i }{x_i}   }{  \bi{\mcal{N}} {n} } & \tx{otherwise} \eec \ee where $\mcal{C}_i \in \bb{R}, \; x_i \in \bb{Z}^+$
for $i = 0, 1, \cd, \ka$ and
\[
\mcal{N} = \sum_{i = 0}^k \mcal{C}_i \neq 0,  \qqu \sum_{i = 0}^\ka x_i = n, \qqu \f{n - 1}{\mcal{N}} < 1, \qqu
 \f{\mcal{C}_i}{\mcal{N}} > 0 \qu \tx{for} \;  i = 0, 1, \cd, k.
\]
By an urn model approach, we have shown in \cite{ChenCDF} that (\ref{GINHY}) indeed defines a distribution.

The multivariate generalized  hypergeometric distribution includes many important distributions as special cases.  Clearly, the multivariate
hypergeometric distribution is obtained from the multivariate generalized hypergeometric distribution by restricting $\mcal{C}_i, \; i = 0, 1,
\cd, k$ and $\mcal{N}$ as positive integers. The multivariate negative hypergeometric distribution is obtained from the multivariate generalized
hypergeometric distribution by taking $\mcal{C}_i, \; i = 0, 1, \cd, k$ and $\mcal{N}$ as negative integers.  The multinomial distribution is
obtained from the multivariate generalized hypergeometric distribution by letting $\mcal{N} \to \iy$ under the constraint that
$\f{\mcal{C}_i}{\mcal{N}}, \; i = 0, 1, \cd, k$ converge to positive numbers sum to $1$.  Moreover, the multivariate P\'{o}lya-Eggenberger
distribution \cite{Steyn} can be accommodated  as a special case of the multivariate generalized hypergeometric distribution.

Let $X_0, X_1, \cd, X_\ka$ be random variables possessing a multivariate generalized hypergeometric distribution defined by (\ref{GINHY}). It
can be shown that the means of these variables are given as \[ \mu_i = \bb{E} [ X_i ] = \f{n \mcal{C}_\ell}{\mcal{N}}, \qqu i = 0, 1, \cd, \ka.
\]
 Define vectors
\[
\bs{x} = [x_0,  x_1, \cd x_\ka ]^\top, \qqu \bs{X} = [X_0,  X_1, \cd X_\ka ]^\top, \qqu \bs{\mu} = [\mu_0, \mu_1, \cd, \mu_\ka]^\top.
\]
Clearly, $\bb{E} [ \bs{X} ] = \bs{\mu}$.  We use $f(\bs{x})$ to denote the probability mass function defined by (\ref{GINHY}), that is,
\[
f(\bs{x}) = \Pr \{ X_i = x_i, \; i = 0, 1, \cd, \ka \}.
\]
 In this setting, we have established the following results.

\beT

\la{GenmulHyper}

Let $z_i, \; i = 0, 1, \cd, k$ be nonnegative integers such that {\small $\sum_{i = 0}^k z_i= n$}. Define {\small $\bs{z} = [z_0, z_1, \cd,
z_\ka]^\top$} and $\wh{\mcal{C}}_i = \mcal{N} \f{z_i}{n}$ for $i = 0, 1, \cd, k$.  Assume that $f(\bs{z}) > 0$.  The following assertions hold.

(I):

\be \la{inevipa}
 \Pr \{ \bs{X} \bs{\prec} \bs{z} \} \leq \prod_{i = 0}^k \f{ \bi{ \mcal{C}_i }{z_i} } { \bi{ \wh{\mcal{C}}_i } {z_i}  }
 \qqu \tx{provided that $\bs{z} \bs{\prec} \bs{\mu}$.}  \ee

(II):

\be \la{inevipb}
 \Pr \{ \bs{X} \bs{\succ} \bs{z} \} \leq \prod_{i = 0}^k \f{ \bi{ \mcal{C}_i }{z_i} } { \bi{ \wh{\mcal{C}}_i } {z_i}  }
 \qqu \tx{provided that $\bs{z} \bs{\succ} \bs{\mu}$.}  \ee

(III):  \be \la{holds889}
 \prod_{i = 0}^k \f{ \bi{ \mcal{C}_i }{z_i} } { \bi{ \wh{\mcal{C}}_i } {z_i}  }  < \li [ \li (
\f{n}{\mcal{N}} \ri )^n \prod_{i = 0}^k \li (  \f{\mcal{C}_i}{z_i} \ri )^{z_i} \ri ] \li ( \f{e^2}{ 2 \pi} \ri )^{k +1} \li ( \f{\mcal{N} -
n}{\mcal{N}} \ri )^{\mcal{N} - n + \f{k+1}{2}} \prod_{i = 0}^k \li ( \f{\mcal{C}_i}{\mcal{C}_i - z_i} \ri )^{ \mcal{C}_i - z_i + \f{1}{2}} \ee
provided that $\mcal{C}_i - z_i \geq 1$ and $\wh{\mcal{C}}_i - z_i \geq 1$ for $i = 0, 1, \cd, k$.

(IV):
\[
\li ( \f{\mcal{N} - n}{\mcal{N}} \ri )^{\mcal{N} - n + \f{k+1}{2}} \prod_{i = 0}^k \li ( \f{\mcal{C}_i}{\mcal{C}_i - z_i} \ri )^{ \mcal{C}_i -
z_i + \f{1}{2}} \to 0 \] as $\mcal{N} \to \iy$ under the constraint that $\f{n}{\mcal{N}} \to \al$ and $\f{\mcal{C}_i}{\mcal{N}} \to \se_i, \;
\f{z_i}{n} \to \ba_i$ for $i = 0, 1, \cd, k$, where $\se_i, \ba_i$ and $\al$ are positive numbers less than $1$ such that $\se_i
> \al \ba_i$ for $i = 0, 1, \cd, k$ and that $\se_i \neq \ba_i$ for some $i \in \{0, 1, \cd, k \}$.

\eeT

See Appendix \ref{GenmulHyperapp} for a proof.

As $\mcal{N} \to \iy$ under the constraint that $\f{\mcal{C}_i}{\mcal{N}} \to p_i, \; i = 0, 1, \cd, k$, \bee \f{ \prod_{i=0}^k \bi{  \mcal{C}_i
}{x_i}   }{  \bi{\mcal{N}} {n} } = \bi{n}{\bs{x}} \li ( \prod_{i=0}^k \prod_{\ell = 1}^{x_i} \f{ \mcal{C}_i - \ell + 1 }{ \mcal{N} - \ell + 1 }
\ri ) \to \bi{n}{\bs{x}} \prod_{i=0}^k {p_i}^{x_i}, \eee where
\[
\bi{n}{\bs{x}} = \f{n!}{ \prod_{i=0}^k x_i! }
\]
is the multinomial coefficient.  This implies that the multivariate generalized hypergeometric distribution converges to the multinomial
distribution \be \la{multnomial} \Pr \{ \bs{X} = \bs{x} \}  = \bi{n}{\bs{x}} \prod_{i=0}^k {p_i}^{x_i}. \ee Accordingly, \bee  \prod_{i = 0}^k
\f{ \bi{ \mcal{C}_i }{z_i} } { \bi{ \wh{\mcal{C}}_i } {z_i}  }  =   \prod_{i=0}^k \prod_{\ell = 1}^{z_i} \f{  \mcal{C}_i - \ell  + 1   }{
\wh{\mcal{C}}_i - \ell + 1 } =
 \prod_{i=0}^k \prod_{\ell = 1}^{z_i} \f{ \f{n \mcal{C}_i}{\mcal{N}}  - \f{n(\ell - 1)}{\mcal{N}} }{ \f{n \wh{\mcal{C}}_i}{\mcal{N}}  - \f{n(\ell - 1)}{\mcal{N}} }
  =  \prod_{i=0}^k \prod_{\ell = 1}^{z_i} \f{ \mu_i  - \f{n(\ell - 1)}{\mcal{N}} }{ z_i  - \f{n(\ell - 1)}{\mcal{N}} }
   \to  \prod_{i = 0}^k \li ( \f{ \mu_i } { z_i } \ri)^{z_i}, \eee where \[ \mu_i = \bb{E} [ X_i ] = n p_i
\]
with $p_i, \; i = 0, 1, \cd, \ka$ being positive real numbers sum to $1$.  As before, let $\bs{X} = [X_0,  X_1, \cd,  X_\ka ]^\top$ and
$\bs{\mu} = [\mu_0, \mu_1, \cd, \mu_\ka]^\top$.  As applications of (\ref{inevipa}) and (\ref{inevipb}) of Theorem \ref{GenmulHyper}, we have
the following results.

\begin{corollary}

\la{Genmulmultinomial}

Let $z_i, \; i = 0, 1, \cd, k$ be nonnegative integers such that $\sum_{i = 0}^k z_i= n$.  Let $\bs{z} = [z_0, z_1, \cd, z_\ka]^\top$. Then, \be
\la{inevipnomiala}
 \Pr \{ \bs{X} \bs{\prec} \bs{z} \} \leq \prod_{i = 0}^k \li ( \f{ \mu_i } { z_i } \ri)^{z_i}  \qu \tx{provided that
$\bs{z} \bs{\prec}  \bs{\mu}$}, \ee and \be \la{inevipnomialb}
 \Pr \{ \bs{X} \bs{\succ} \bs{z} \} \leq \prod_{i = 0}^k \li ( \f{ \mu_i } { z_i } \ri)^{z_i} \qu \tx{provided that $\bs{z} \bs{\succ} \bs{\mu}$}.  \ee

\end{corollary}

It should be noted that Corollary \ref{Genmulmultinomial} is a multivariate generalization of the Chernoff bounds for binomial distributions.

\subsection{Multivariate Generalized Inverse Hypergeometric Distribution}

Let $\ka$ be a positive integer.  In probability theory, random variables $X_1, \cd, X_\ka$ are said to possess a {\it multivariate inverse
hypergeometric distribution} if  \[
\Pr \{ X_i = x_i, \; i = 1, \cd, \ka \} = \bec 0 & \tx{if $n \geq 1 + N$ or $x_i \geq 1 + C_i$ for some $i \in \{ 1, \cd, k \}$},\\
\f{\ga}{n} \f{ \prod_{i=0}^k \bi{ C_i }{x_i}   }{  \bi{N} {n} } & \tx{otherwise} \eec
\]
where $C_i \in \bb{N}, \; x_i \in \bb{Z}^+$ for $i = 1, \cd, k$ and
\[
C_0 \in \bb{N}, \qqu \ga \in \bb{N}, \qqu x_0 = \ga \leq C_0, \qqu n = \sum_{i=0}^k x_i, \qqu N = \sum_{i = 0}^k C_i. \]

Actually,  under mild restrictions, the multivariate inverse hypergeometric distribution can be generalized by allowing $N$ and $C_i$ to be real
numbers. More formally, random variables $X_1, \cd, X_\ka$ are said to possess a {\it multivariate generalized inverse hypergeometric
distribution} if \be \la{geninvhhyper89}
\Pr \{ X_i = x_i, \; i = 1, \cd, \ka \} = \bec 0 & \tx{if $\f{n - 1}{\mcal{N}} \geq 1$ or $x_i \geq 1 + \mcal{C}_i > 1$ for some $i \in \{ 1, \cd, k \}$},\\
\f{\ga}{n} \f{ \prod_{i=0}^k \bi{ \mcal{C}_i }{x_i}   }{  \bi{\mcal{N}} {n} } & \tx{otherwise} \eec \ee where $\mcal{C}_i \in \bb{R}, \; x_i \in
\bb{Z}^+$ for $i = 1, \cd, k$,
\[
\mcal{C}_0 \in \bb{R}, \qqu \ga \in \bb{N}, \qqu x_0 = \ga, \qqu n = \sum_{i = 0}^k x_i, \qqu \mcal{N} = \sum_{i = 0}^k \mcal{C}_i \neq 0, \qqu
\f{\mcal{C}_0} {\mcal{N}}
> \f{\ga - 1}{\mcal{N}},
\]
and $\f{\mcal{C}_i}{\mcal{N}} > 0$ for $i = 1, \cd, k$.  By an urn model approach, we have shown in \cite{ChenCDF} that (\ref{geninvhhyper89})
indeed defines a distribution.

Let $X_1, \cd, X_\ka$ be random variables possessing a multivariate generalized inverse hypergeometric distribution defined by
(\ref{geninvhhyper89}). Let $X_0 =\ga$.  It can be shown that the means of these variables are given as \[ \mu_i = \bb{E} [ X_i ] = \f{\ga
\mcal{C}_\ell}{\mcal{C}_0}, \qqu i = 0, 1, \cd, \ka.
\]
 Define vector
\[
\bs{X} = [ X_0, X_1, \cd,  X_\ka ]^\top.
\]
Let the probability mass function defined by (\ref{geninvhhyper89}) be denoted by $f(\bs{x})$, that is,
\[ f(\bs{x}) = \Pr \{ X_i = x_i, \; i = 0, 1, \cd, \ka \},
\]
where $\bs{x} = [ x_0, x_1, \cd, x_\ka]^\top$.   In this setting, we have established the following results.

\beT

\la{GenmulInvHyper}

Let $z_0 = \ga$ and $z_i, \; i = 1, \cd, k$ be nonnegative integers. Define $n = \sum_{i = 0}^k z_i$ and $\bs{z} = [z_0, z_1, \cd, z_\ka]^\top$.
Define $\wh{\mu}_i = \f{n \mcal{C}_i}{\mcal{N}}$ and $\wh{\mcal{C}}_i = \mcal{N} \f{z_i}{n}$ for $i = 0, 1, \cd, k$.  Define $\wh{\bs{\mu}} =
[\wh{\mu}_0, \wh{\mu}_1, \cd, \wh{\mu}_\ka]^\top$.  Assume that $\f{n - 1}{\mcal{N}} < 1$ and $f(\bs{z}) > 0$. Then, the following assertions
hold.

(I):

\be \la{inevipainv}
 \Pr \{ \bs{X} \bs{\prec} \bs{z} \} \leq \prod_{i = 0}^k \f{ \bi{ \mcal{C}_i }{z_i} } { \bi{ \wh{\mcal{C}}_i } {z_i}  }  \qu \tx{provided that
$\bs{z} \bs{\prec} \wh{\bs{\mu}}$}. \ee

(II):

\be \la{inevipbinv}
 \Pr \{ \bs{X} \bs{\succ} \bs{z} \} \leq \prod_{i = 0}^k \f{ \bi{ \mcal{C}_i }{z_i} } { \bi{ \wh{\mcal{C}}_i } {z_i}  } \qu \tx{provided that
$\bs{z} \bs{\succ} \wh{\bs{\mu}}$}.  \ee

(III):  \[ \prod_{i = 0}^k \f{ \bi{ \mcal{C}_i }{z_i} } { \bi{ \wh{\mcal{C}}_i } {z_i}  }  < \li [ \li ( \f{n}{\mcal{N}} \ri )^n \prod_{i = 0}^k
\li (  \f{\mcal{C}_i}{z_i} \ri )^{z_i} \ri ] \li ( \f{e^2}{ 2 \pi} \ri )^{k +1} \li ( \f{\mcal{N} - n}{\mcal{N}} \ri )^{\mcal{N} - n +
\f{k+1}{2}} \prod_{i = 0}^k \li ( \f{\mcal{C}_i}{\mcal{C}_i - z_i} \ri )^{ \mcal{C}_i - z_i + \f{1}{2}}
\]
provided that $\mcal{C}_i - z_i \geq 1$ and $\wh{\mcal{C}}_i - z_i \geq 1$ for $i = 0, 1, \cd, k$.

(IV):
\[
\li ( \f{\mcal{N} - n}{\mcal{N}} \ri )^{\mcal{N} - n + \f{k+1}{2}} \prod_{i = 0}^k \li ( \f{\mcal{C}_i}{\mcal{C}_i - z_i} \ri )^{ \mcal{C}_i -
z_i + \f{1}{2}} \to 0 \] as $\mcal{N} \to \iy$ under the constraint that $\f{n}{\mcal{N}} \to \al$ and $\f{\mcal{C}_i}{\mcal{N}} \to \se_i, \;
\f{z_i}{n} \to \ba_i$ for $i = 0, 1, \cd, k$, where $\se_i, \ba_i$ and $\al$ are positive numbers less than $1$ such that $\se_i
> \al \ba_i$ for $i = 0, 1, \cd, k$ and that $\se_i \neq \ba_i$ for some $i \in \{0, 1, \cd, k \}$.

\eeT

See Appendix \ref{GenmulInvHyperapp} for a proof.

As $\mcal{N} \to \iy$ under the constraint that $\f{\mcal{C}_i}{\mcal{N}} \to p_i$ for $i = 0, 1, \cd, k$, \bee \f{\ga}{n} \f{ \prod_{i=0}^k
\bi{  \mcal{C}_i }{x_i}   }{  \bi{\mcal{N}} {n} } = \f{\ga}{n} \bi{n}{\bs{x}} \li ( \prod_{i=0}^k \prod_{\ell = 1}^{x_i} \f{ \mcal{C}_i - \ell +
1 }{ \mcal{N} - \ell + 1  } \ri ) \to \f{\ga }{n} \bi{n}{\bs{x}}  \prod_{i=0}^k {p_i}^{x_i}. \eee This implies that the multivariate generalized
inverse hypergeometric distribution converges to the negative multinomial distribution \be \la{multnomialinv} \Pr \{ X_i = x_i, \; i = 1, \cd,
\ka \} = \f{\ga}{n}  \bi{n}{\bs{x}} \prod_{i=0}^k {p_i}^{x_i}. \ee Accordingly, \bee  \prod_{i = 0}^k \f{ \bi{ \mcal{C}_i }{z_i} } { \bi{
\wh{\mcal{C}}_i } {z_i}  }  =  \prod_{i=0}^k \prod_{\ell = 1}^{z_i} \f{  \mcal{C}_i - \ell  + 1   }{ \wh{\mcal{C}}_i - \ell + 1 } =
 \prod_{i=0}^k \prod_{\ell = 1}^{z_i} \f{ \f{n \mcal{C}_i}{\mcal{N}}  - \f{n(\ell - 1)}{\mcal{N}} }{ \f{n \wh{\mcal{C}}_i}{\mcal{N}}  - \f{n(\ell - 1)}{\mcal{N}} }
  =  \prod_{i=0}^k \prod_{\ell = 1}^{z_i} \f{ n p_i  - \f{n(\ell - 1)}{\mcal{N}} }{ z_i  - \f{n(\ell - 1)}{\mcal{N}} }
   \to  \prod_{i = 0}^k \li ( \f{ n p_i } { z_i } \ri)^{z_i}. \eee Define
   \[  \bs{X} = [X_0, X_1, \cd, X_\ka]^\top
\]
with $X_1, \cd, X_\ka$ possessing a negative multinomial distribution defined by (\ref{multnomialinv}). As applications of (\ref{inevipainv})
and (\ref{inevipbinv}) of Theorem \ref{GenmulInvHyper}, we have the following results.

\begin{corollary}

\la{Genmulmultinomialinv}

Let $z_0 = \ga$ and $z_i, \; i = 1, \cd, k$ be nonnegative integers. Define $n = \sum_{i = 0}^k z_i$ and $\bs{z}  =[ z_0, z_1, \cd,
z_\ka]^\top$.  Define $\wh{\bs{\mu}} = [\wh{\mu}_0, \wh{\mu}_1, \cd, \wh{\mu}_\ka]^\top$, where $\wh{\mu}_i = n p_i, \; i = 0, 1, \cd, \ka$.
Then, \be \la{inevipnomialainv}
 \Pr \{ \bs{X} \bs{\prec} \bs{z} \} \leq \prod_{i = 0}^k \li ( \f{ \wh{\mu}_i } { z_i } \ri)^{z_i}  \qu \tx{provided that
$\bs{z} \bs{\prec} \wh{\bs{\mu}}$}, \ee and \be \la{inevipnomialbinv}
 \Pr \{ \bs{X} \bs{\succ} \bs{z} \} \leq \prod_{i = 0}^k \li ( \f{ \wh{\mu}_i } { z_i } \ri)^{z_i} \qu \tx{provided that $\bs{z} \bs{\succ} \wh{\bs{\mu}}$}.  \ee

\end{corollary}

\subsection{Dirichlet Distribution}

Let $\ka$ be a positive integer.   Random variables $X_0, X_1, \cd, X_\ka$ are said to possess a Dirichlet distribution if they have a
probability density function \[ f (\bs{x}, \bs{\al})  = \bec \f{1}{ \mcal{B} (\bs{\al})} \prod_{i=0}^\ka
x_i^{\al_i - 1} & \tx{for $x_i \in \bb{R}^+, \; i = 0, 1\cd, \ka$ such that  $\sum_{i=0}^\ka x_i = 1$},\\
0 & \tx{else} \eec
\]
where $\bs{x} = [x_0, x_1, \cd, x_\ka]^\top$, $\al_0, \al_1, \cd, \al_\ka$ are positive real numbers, $\bs{\al} = [ \al_0, \al_1, \cd, \al_\ka
]^\top$, and
\[
\mcal{B} (\bs{\al}) =  \f{ \prod_{i=0}^\ka \Ga (\al_i) }{ \Ga ( \sum_{i=0}^\ka \al_i)  }.
\]
The means of $X_0, X_1, \cd, X_\ka$ are
\[
\mu_i = \bb{E} [X_i] = \f{\al_i}{\sum_{\ell = 0}^\ka \al_\ell} \qqu \tx{for $i = 0, 1, \cd, \ka$}.
\]
Define vectors
\[
\bs{X} = [X_0, X_1, \cd, X_\ka]^\top, \qqu \bs{\mu} = [\mu_0, \mu_1, \cd, \mu_\ka]^\top.
\]
Let $\bs{\mcal{X}}_1, \bs{\mcal{X}}_2, \cd, \bs{\mcal{X}}_n$ be independent random vectors possessing the same distribution as $\bs{X}$.  Define
$\ovl{\bs{\mcal{X}}}_n = \f{ \sum_{\ell = 1}^n \bs{\mcal{X}}_i }{n}$.  We have the following result.

\beT

\la{mulDrich} Let $z_\ell, \; \ell = 0, 1, \cd, \ka$ be positive real numbers such that $\sum_{\ell = 0}^k z_\ell= 1$. Define $\bs{z} = [z_0,
z_1, \cd, z_\ka]^\top$.  Then, \be \la{conthm1388} \Pr \{ \ovl{\bs{\mcal{X}}}_n  \bs{\prec} \bs{z} \} \leq \li [ \f{ \mcal{B} (\wh{\bs{\al}}) }{
\mcal{B} (\bs{\al}) } \prod_{i=1}^\ka  \f{z_i^{\al_i}}{z_i^{\wh{\al}_i}}  \ri ]^n \qqu \tx{provided that $\bs{z} \bs{\prec} \bs{\mu}$}, \ee
where $\wh{\bs{\al}} = [\wh{\al}_0, \wh{\al}_1, \cd, \wh{\al}_\ka]^\top$ with $\wh{\al}_i = \f{ \al_0  z_i}{ z_0}$ for $i = 0, 1, \cd, \ka$.

 \eeT

 See Appendix \ref{mulDrichapp} for a proof.

\subsection{Matrix Gamma Distribution}

In statistics, a matrix gamma distribution is a generalization of the gamma distribution to positive-definite matrices (see, \cite{Gupta} and
the references therein). It is a more general version of the Wishart distribution.  The pdf of a matrix gamma distribution is \be \la{CDFmga}
 f(\bs{x}) = \f{| \bs{\Si} |^{-\al}}{ \ba^{p \al}
\Ga_p (\al ) } | \bs{x} |^{\al -( p + 1) \sh 2} \exp \li (  - \f{1}{\ba} \tx{tr} ( \bs{\Si}^{-1} \bs{x} )\ri ), \qqu \al > 0, \;  \ba > 0 \ee
where $\bs{x}$ and  $\bs{\Si}$ are positive definite matrices of size $p \times p$. Here, the notation ``$|.|$'' denotes the determinant
function. Let $\bs{\mcal{X}}_1, \bs{\mcal{X}}_2, \cd, \bs{\mcal{X}}_n$ be i.i.d. positive definite random matrices possessing the same pdf
defined by (\ref{CDFmga}). Define
\[
\ovl{\bs{\mcal{X}}}_n = \f{ \sum_{\ell = 1}^n \bs{\mcal{X}}_i }{n}, \qqu \bs{\mu} = \al \ba \bs{\Si}.
\]
Clearly, $\bb{E} [ \ovl{\bs{\mcal{X}}}_n ] = \bs{\mu}$.  We write $\bs{x} \prec \bs{y}$ if $\bs{y} - \bs{x}$ is positive definite.  Similarly,
we write $\bs{x} \succ \bs{y}$ if $\bs{x} - \bs{y}$ is positive definite.   In this setting, we have the following results.

\beT

\la{matrixgga} \bee &  &  \Pr \{ \ovl{\bs{\mcal{X}}}_n \bs{\succ} \bs{z} \} \leq  \li [ \li (  \f{e}{\al \ba} \ri )^{ p \al } \li (  \f{ |
\bs{z} | } { | \bs{\Si} | } \ri )^{\al} \exp \li (   - \f{1}{\ba} \mrm{tr} ( \bs{\Si}^{-1} \bs{z} ) \ri ) \ri ]^n \qu \tx{for $\bs{z} \bs{\succ}
\bs{\mu}$},\\
&  & \Pr \{ \ovl{\bs{\mcal{X}}}_n \bs{\prec} \bs{z} \} \leq  \li [ \li (  \f{e}{\al \ba} \ri )^{ p \al } \li (  \f{ | \bs{z} | } { | \bs{\Si} |
} \ri )^{\al} \exp \li (   - \f{1}{\ba} \mrm{tr} ( \bs{\Si}^{-1} \bs{z} ) \ri ) \ri ]^n \qu \tx{for $\bs{z} \bs{\prec} \bs{\mu}$}. \eee In
particular, $\Pr \{ \ovl{\bs{\mcal{X}}}_n \bs{\succ} \vro \bs{\mu} \} \leq  \li [  \vro \exp (1 - \vro ) \ri ]^{n p \al}$ for $\vro  \geq 1$ and
$\Pr \{ \ovl{\bs{\mcal{X}}}_n \bs{\prec} \vro \bs{\mu} \} \leq  \li [  \vro \exp (1 - \vro ) \ri ]^{n p \al}$ for $0 \leq \vro \leq 1$.

\eeT

See Appendix \ref{matrixggaapp} for a proof.

\section{Conclusion}

We have developed the LR method for deriving probabilistic inequalities.  We have established fundamental connections among the LR method, the
theory of large deviations and  statistical concepts such as likelihood ratio, maximum likelihood, and the method of moments for parameter
estimation. The LR method overcomes the limitations of the classical approach based on mathematical expectation. The LR method may provide easy
derivation of concentration inequalities in situations that moment generating functions are not readily tractable. We have applied the LR method
to obtain a wide spectrum of new concentration inequalities.

\appendix

\sect{Proof of Theorem \ref{con89666} } \la{con89666app}

By the assumption that  $f(\bs{X}) \leq \exp ( \ro ( \bs{X}, \bs{\vse}) ) \; g(\bs{X}, \bs{\vse})$ for all $\bs{\vse} \in \varTheta$, we have
\be \la{use86639a}
 \prod_{i=1}^n f(\bs{X}_i) \leq \exp \li ( \sum_{i=1}^n \ro ( \bs{X}_i, \bs{\vse}) \ri ) \; \prod_{i=1}^n
g(\bs{X}_i, \bs{\vse}) \qqu \tx{for all $\bs{\vse} \in \varTheta$}. \ee By the assumption that $\ro (\bs{x}, \bs{\vse} )$ is a concave function
of $\bs{x}$, we have \be \la{use86639b} \exp \li ( \sum_{i=1}^n \ro ( \bs{X}_i, \bs{\vse}) \ri ) \leq \exp \li ( n \ro ( \ovl{\bs{X}}_n,
\bs{\vse}) \ri ) \qqu \tx{for all $\bs{\vse} \in \varTheta$}. \ee Combining (\ref{use86639a}) and (\ref{use86639b}) yields \be \la{mulp8996}
 \prod_{i=1}^n f(\bs{X}_i) \leq \exp ( n
\ro ( \ovl{\bs{X}}_n, \bs{\vse}) ) \; \prod_{i=1}^n  g(\bs{X}_i, \bs{\vse}) \ee for all $\bs{\vse} \in \varTheta$. Multiplying (\ref{mulp8996})
by the indicator function $\bb{I}_{ \{ \ovl{\bs{X}}_n  \bs{\prec} \bs{z} \} }$ yields \be \la{foll899}
 \li [ \prod_{i=1}^n f(\bs{X}_i) \ri ] \bb{I}_{ \{ \ovl{\bs{X}}_n  \bs{\prec} \bs{z} \} } \leq \exp ( n \ro
( \ovl{\bs{X}}_n, \bs{\vse}) ) \; \li [ \prod_{i=1}^n  g(\bs{X}_i, \bs{\vse}) \ri ] \bb{I}_{ \{ \ovl{\bs{X}}_n \bs{\prec} \bs{z} \} } \qqu
\tx{for all $\bs{\vse} \in \varTheta$}. \ee Since $\bs{z}$ is a matrix such that $\ro (\bs{x}, \bs{\vse}) \leq \ro (\bs{z}, \bs{\vse})$ for all
$\bs{\vse} \in \varTheta$ provided that $\bs{x}  \bs{\prec} \bs{z}$, it follows from (\ref{foll899}) that
\[
\li [ \prod_{i=1}^n f(\bs{X}_i) \ri ] \bb{I}_{ \{ \ovl{\bs{X}}_n  \bs{\prec} \bs{z} \} } \leq \exp ( n \ro ( \bs{z}, \bs{\vse}) ) \; \li [
\prod_{i=1}^n  g(\bs{X}_i, \bs{\vse}) \ri ] \bb{I}_{ \{ \ovl{\bs{X}}_n  \bs{\prec} \bs{z} \} } \qqu \tx{for all $\bs{\vse} \in \varTheta$},
\]
which implies that \be \la{fin86696}
 \li [ \prod_{i=1}^n f(\bs{X}_i) \ri ] \bb{I}_{ \{ \ovl{\bs{X}}_n  \bs{\prec} \bs{z} \} } \leq
\exp ( n \ro ( \bs{z}, \bs{\vse}) ) \; \prod_{i=1}^n  g(\bs{X}_i, \bs{\vse}) \qqu \tx{for all $\bs{\vse} \in \varTheta$}. \ee Applying
(\ref{fin86696}) and Theorem \ref{ThM888} yields (\ref{concentra8996}).  This completes the proof of the theorem.

\section{Proof of Theorem \ref{convip} } \la{convipapp}

We need some preliminary results.  The results in the following lemma are due to Petrov \cite{Petrov}.

\beL \la{petrov} Let $X$ be a non-degenerate random variable with mean $\mu = \bb{E} [ X ]$ and  moment generating function $\phi(s) = \bb{E} [
e^{s X}]$ for $s \in (-a, b)$,  where $a$ and $b$ are $\iy$ or positive real numbers.  Let $\psi (s) = \ln \phi (s)$.   Let $X_1, X_2, \cd$ be
i.i.d. samples of $X$. Let $\ovl{X}_n = \f{\sum_{i=1}^n X_i}{n}$ for $n \in \bb{N}$.  Let $\al > \mu$ be a real number such that there exists
$\tau_\al \in (0, b)$ satisfying $\psi^\prime (\tau_\al) = \al$.  Let  $\ga (\al) = \al \tau_\al - \psi (\tau_\al )$. The following assertions
hold.

(I)  If $X$ is non-lattice valued, then
\[
\Pr \{ \ovl{X}_n \geq \al \} =  \f{\exp ( - n  \ga (\al) )}{ \tau_\al \sq{  2 \pi n \psi^{\prime \prime} (\tau_\al) }  } [ 1 + o(1) ].
\]

(II) If $X$ is a lattice-valued random variable with span $\nu$, then
\[
\Pr \{ \ovl{X}_n \geq \al \} = \f{\nu}{  \sq{  2 \pi n \psi^{\prime \prime} (\tau_\al ) }  } \f{ \exp ( - n  \ga(\al) ) } { 1 - \exp ( - \nu
\tau_\al ) } [ 1 + o(1) ].
\]

\eeL

\beL  \la{lemuseful} For $z \geq \mu$ and $\Pr \{ X > z \} > 0$, \be \la{lemcon89}
 \lim_{n \to \iy} \f{1}{n} \ln \f{ \Pr \{ z \leq \ovl{X}_n < z
+ \vep \} }{ \Pr \{ \ovl{X}_n \geq z \}  } = 0 \ee  provided that $\vep > 0$ is sufficiently small.  \eeL

\bpf

To show the lemma, we need to show that  $\Pr \{ \ovl{X}_n \geq z \} \geq \Pr \{ z \leq \ovl{X}_n < z + \vep \} > 0$ for large enough $n$ and
that (\ref{lemcon89}) holds for small enough $\vep > 0$. Since
\[
\ln \f{ \Pr \{ z \leq \ovl{X}_n < z + \vep \} }{  \Pr \{ \ovl{X}_n \geq z \}  } = \ln \li (  1 -  \f{ \Pr \{ \ovl{X}_n \geq z + \vep \} }{ \Pr
\{ \ovl{X}_n \geq z \}  } \ri ),
\]
to show (\ref{lemcon89}), it suffices to show that \be \la{goallem}
 \lim_{n \to \iy} \f{ \Pr \{ \ovl{X}_n \geq z + \vep \} }{ \Pr \{ \ovl{X}_n
\geq z \} } = 0 \ee holds for small enough $\vep > 0$.

 If $z = \mu$, then \be \la{special8996}
 \f{ \Pr \{
\ovl{X}_n \geq z + \vep \} }{ \Pr \{ \ovl{X}_n \geq z \} }  =  \f{ \Pr \{ \ovl{X}_n \geq \mu + \vep \} }{ \Pr \{ \ovl{X}_n \geq \mu \} } \leq
\f{ \Pr \{ | \ovl{X}_n - \mu | \geq \vep \} }{  \Pr \{ \f{ \sq{n} (\ovl{X}_n - \mu)}{\si} \geq 0 \} },  \ee where $\si > 0$ is the standard
deviation of $X$.  Making use of (\ref{special8996}) and the observation that $\lim_{n \to \iy} \Pr \{ | \ovl{X}_n - \mu | \geq \vep \} = 0$ due
to the weak law of large numbers and that $\lim_{n \to \iy} \Pr \{ \f{ \sq{n} (\ovl{X}_n - \mu)}{\si} \geq 0 \} = \f{1}{2}$ as a result of the
central limit theorem,  we have that $\Pr \{ \ovl{X}_n \geq z \} \geq \Pr \{ z \leq \ovl{X}_n < z + \vep \} > 0$ for large enough $n$
 and that (\ref{goallem}) holds.   In the sequel, we will show (\ref{goallem}) by restricting $z$ to be greater than $\mu$.

Define $\varphi(s) = \bb{E} [ e^{s (X - z)} ]$ and \be \la{defb88}
 b = \sup \{ c \geq 0: \varphi(s) < \iy \; \tx{for any} \; s \in [0, c] \}.
\ee Invoking the assumption that $\phi(s) = \bb{E} [ e^{s X} ] < \iy$ for $s$ in a neighborhood of $0$, we have that either $b = \iy$ or $b$ is
a positive number. Hence, $\varphi(s) < \iy$ and $\phi(s) < \iy$ for $0 < s < b$.  We claim that $\lim_{s \uparrow b} \varphi(s) = \iy$.   To
show the claim, note that, as a consequence of the assumption that $\Pr \{ X > z \} > 0$,  there exists a number $y$ greater than $z$ such that
$\Pr \{ X \geq y \}
> 0$.  This implies that $\varphi(s) \geq e^{s (y - z)}  \Pr \{ X \geq y \}$ for $s > 0$ such that $\phi(s) < \iy$.
Therefore, it must be true that either $\lim_{s \to \iy} \varphi(s) = \iy$ or there exists a positive number $\ga > 0$ such that $\lim_{s
\uparrow \ga} \varphi(s) = \iy$.  This fact together with the definition (\ref{defb88}) for $b$ imply that the claim is true.

Note that $\varphi(s) = \psi(s) - s z$, where $\psi (s) = \ln \phi(s)$.  Since $\varphi(s)$ is a convex function of $s \in (0, b)$ and
\[
\varphi(0) = 0, \qqu \varphi^\prime (0) = \psi^\prime (0) - z = \mu - z < 0, \qqu \lim_{s \uparrow b} \varphi(s) = \iy,
\]
there must exist a unique number $\tau_z \in (0, b)$ such that $\varphi^\prime (\tau_z) = \psi^\prime (\tau_z) -  z = 0$.  Let $\ze = \f{\tau_z
+ b}{2}$.  Then, $\psi^\prime (\ze) > z$, since $\ze > \tau_z$ and $\psi(.)$ is convex.  Let $\al$ be a number such that $z \leq \al <
\psi^\prime (\ze)$. Since $\mu < z = \psi^\prime (\tau_z) \leq \al < \psi^\prime (\ze)$, there must exist a unique number $\tau_{\al} \in (0,
b)$ such that $\psi^\prime (\tau_{\al}) =  \al$.  Clearly, $\tau_\al$ is a function of $\al$. Define $\ga (\al) = \al \tau_\al - \psi (\tau_\al
)$ for $\al$ such that $z \leq \al < \psi^\prime (\ze)$.  The derivative of $\ga (\al)$ with respect to $\al$ is \bee \ga^\prime (\al)  = \f{ d
}{d \al} \li [  \al \tau_\al - \psi (\tau_\al ) \ri ]  =  \tau_\al + \al \f{d \tau_\al}{d \al} - \psi^\prime (\tau_\al) \f{d \tau_\al}{d \al}
 =  \tau_\al + \al \f{d \tau_\al}{d \al} - \al \f{d \tau_\al}{d \al}
  =  \tau_\al > 0, \eee which implies that $\ga(\al)$ is strictly increasing with respect to $\al$ such that $z \leq \al < \psi^\prime (\ze)$.
From now on, we restrict $\vep$ to be positive and less than $\psi^\prime (\ze) - z$.  Since $0 < \vep < \psi^\prime (\ze) - z$, we have that
$\tau_{z+ \vep}$ exists and that \be \la{use8996a}
 0 < \tau_z < \tau_{z+ \vep} < b, \qqu 0 < \ga(z) < \ga(z + \vep).
\ee If $X$ is non-lattice valued, then according to assertion (I) of Lemma \ref{petrov}, we have \bel &  & \Pr \{ \ovl{X}_n \geq z \} =  \f{\exp
( -
n \ga (z) )}{ \tau_z \sq{  2 \pi n \psi^{\prime \prime} (\tau_z) }  } [ 1 + o(1) ], \la{use8996b}\\
&  & \Pr \{ \ovl{X}_n \geq z + \vep \} =  \f{\exp ( - n  \ga (z + \vep) )}{ \tau_{z + \vep} \sq{  2 \pi n \psi^{\prime \prime} (\tau_{z+ \vep})
} } [ 1 + o(1) ]. \la{use8996c} \eel  It follows from (\ref{use8996a}), (\ref{use8996b}) and (\ref{use8996c}) that  $\Pr \{ \ovl{X}_n \geq z \}
\geq \Pr \{ z \leq \ovl{X}_n < z + \vep \} > 0$ for large enough $n$ and that \bee  \f{ \Pr \{ \ovl{X}_n \geq z + \vep \} }{ \Pr \{ \ovl{X}_n
\geq z \} }  =  \f{ \tau_z }{ \tau_{z + \vep} } \sq{ \f{ \psi^{\prime \prime} (\tau_z) }{  \psi^{\prime \prime} (\tau_{z+ \vep}) } } \li [ \f{
\exp( \ga(z)  )  }{ \exp (\ga (z + \vep) ) } \ri ]^n [ 1 + o(1) ]  \to 0 \eee as $n \to \iy$.  This shows (\ref{goallem}) in the case that $X$
is non-lattice valued.

If $X$ is a lattice-valued random variable with span $\nu$, then according to assertion (II) of Lemma \ref{petrov}, we have \bel &  & \Pr \{
\ovl{X}_n \geq z \} =  \f{\nu}{ 1 - \exp ( - \nu \tau_z ) } \f{\exp ( -
n \ga (z) )}{  \sq{  2 \pi n \psi^{\prime \prime} (\tau_z) }  } [ 1 + o(1) ], \la{use8996b88}\\
&  & \Pr \{ \ovl{X}_n \geq z + \vep \} =  \f{\nu}{ 1 - \exp ( - \nu \tau_{z + \vep} ) }  \f{\exp ( - n  \ga (z + \vep) )}{ \sq{ 2 \pi n
\psi^{\prime \prime} (\tau_{z+ \vep}) } } [ 1 + o(1) ]. \la{use8996c88} \eel  It follows from (\ref{use8996a}), (\ref{use8996b88}) and
(\ref{use8996c88}) that  $\Pr \{ \ovl{X}_n \geq z \} \geq \Pr \{ z \leq \ovl{X}_n < z + \vep \} > 0$ for large enough $n$ and that  \bee \f{ \Pr
\{ \ovl{X}_n \geq z + \vep \} }{ \Pr \{ \ovl{X}_n \geq z \} } & = & \f{  1 - \exp ( - \nu \tau_{z} ) }{ 1 - \exp ( - \nu \tau_{z + \vep} ) }
\sq{ \f{ \psi^{\prime \prime} (\tau_z) }{ \psi^{\prime \prime} (\tau_{z+ \vep}) } } \li [ \f{ \exp( \ga(z) ) }{ \exp (\ga (z + \vep) ) } \ri ]^n
[ 1 + o(1) ]  \to 0 \eee as $n \to \iy$. This shows (\ref{goallem}) in the case that $X$ is a lattice-valued random variable.  The proof of the
lemma is thus completed.

\epf

\beL \la{CLTCHB} \[ \lim_{n \to \iy} \f{1}{n} \ln \bb{P}_g \{ z \leq \ovl{X}_n < z + \vep \} = 0 \qqu \tx{for any $\vep > 0$}.
\]

\eeL

\bpf

For notational simplicity, let $\si =  \sq{ \bb{E}_g [ | X - z |^2 ] }$.  If $\si = 0$, then $\bb{P}_g \{ z \leq \ovl{X}_n < z + \vep \} = 1$
and the lemma is obviously true.  So, it remains to show the lemma for the case that $0 < \si < \iy$.

Since the random variable $X$ has mean $z$ and finite variance $\si
> 0$ associated with pdf or pmf $g(.)$, it follows from the central limit theorem that $\f{ \sq{n} ( \ovl{X}_n - z ) }{\si}$ converges to a Gaussian random
variable with zero mean and unit variance.  Therefore, there exists an integer $m > 0$ such that
\[
\f{1}{3} < \bb{P}_g \li \{ \f{ \sq{n} ( \ovl{X}_n - z ) }{\si} \geq 0  \ri \} < \f{2}{3}
\]
for any integer $n$ greater than $m$.  Since the variance $\si^2$ is finite, it follows from the Chebyshev inequality that
\[
\bb{P}_g \{ \ovl{X}_n \geq z + \vep \} \leq \bb{P}_g \{ | \ovl{X}_n  - z| \geq  \vep \} \leq \f{ \si^2}{n \vep^2}.
\]
Hence, for $n > \max \{ m, \f{ 3 \si^2 }{ \vep^2 } \}$, we have
\[
\bb{P}_g \{ z \leq \ovl{X}_n < z + \vep  \} \leq \bb{P}_g \{ \ovl{X}_n \geq z \} = \bb{P}_g \li \{ \f{ \sq{n} ( \ovl{X}_n - z ) }{\si} \geq 0
\ri \} < \f{2}{3}
\]
and
\bee  \bb{P}_g \{ z \leq \ovl{X}_n < z + \vep  \} & = & \bb{P}_g \{ \ovl{X}_n \geq z \} - \bb{P}_g \{ \ovl{X}_n \geq z + \vep \}\\
& = & \bb{P}_g \li \{ \f{ \sq{n} ( \ovl{X}_n - z ) }{\si} \geq 0  \ri \} - \bb{P}_g \{ \ovl{X}_n \geq z + \vep \} \\
&  > & \f{1}{3} -  \f{ \si^2}{n \vep^2} > 0. \eee It follows that
\[
\f{1}{n} \ln \li ( \f{1}{3} -  \f{ \si^2}{n \vep^2} \ri ) < \f{1}{n} \ln \bb{P}_g \{ z \leq \ovl{X}_n < z + \vep \} < \f{1}{n} \ln \f{2}{3}
\]
for $n > \max \{ m, \f{ 3 \si^2 }{ \vep^2 } \}$.  Thus,  $\lim_{n \to \iy} \f{1}{n} \ln \bb{P}_g \{ z \leq \ovl{X}_n < z + \vep \} = 0$ for any
$\vep
> 0$.

\epf

\beL \la{lem3388} \bel &  & \limsup_{n \to \iy} \f{1}{n} \ln \Pr \{
z \leq \ovl{X}_n < z + \vep  \} \leq \sup_{y \in [z, z + \vep) } \ro ( y ), \la{limineqa}\\
&  & \liminf_{n \to \iy} \f{1}{n} \ln \Pr \{ z \leq \ovl{X}_n < z + \vep \} \geq \inf_{y \in [z, z + \vep) } \ro ( y ) \la{limineqb} \eel
provided that $\vep
> 0$ is sufficiently small. \eeL

\bpf

We only prove the lemma for the case that $f(.)$ is a pmf, since the proof for the case that $f(.)$ is a pdf is similar. For simplicity of
notations, let $\ovl{x}_n$ denote a realization of $\ovl{X}_n$.  Let $0 < \vep < \de$.  From (\ref{cona88}), we have  \bel \Pr \{ z \leq
\ovl{X}_n < z + \vep \} & =
& \sum_{ z \leq \ovl{x}_n < z + \vep } \; \prod_{i=1}^n f ( x_i ) \nonumber\\
& \leq & \sum_{ z \leq \ovl{x}_n < z + \vep  } \exp \li ( n \ro ( \ovl{x}_n ) + o(n) \ri ) \times \prod_{i=1}^n g ( x_i ) \nonumber\\
&  \leq & \sup_{y \in [z, z + \vep) } \exp \li ( n \ro ( y )  + o (n) \ri ) \sum_{ z \leq \ovl{x}_n < z + \vep  } \; \prod_{i=1}^n g ( x_i ) \nonumber \\
& = & \sup_{y \in [z, z + \vep) } \exp \li ( n \ro ( y )  + o (n)  \ri ) \; \bb{P}_g \{ z \leq \ovl{X}_n < z + \vep \}. \la{CH86336} \eel From
the assumption that $\bb{E}_g [ X ] = z$ and $\bb{E}_g [ | X - z |^2 ] < \iy$, according to Lemma \ref{CLTCHB}, we have that \be \la{CLTCH8996}
\lim_{n \to \iy} \f{1}{n} \ln \bb{P}_g \{ z \leq \ovl{X}_n < z + \vep \} = 0. \ee  Combining (\ref{CH86336}) and (\ref{CLTCH8996}) yields
(\ref{limineqa}).

From (\ref{conb88}), we have  \bel \Pr \{ z \leq \ovl{X}_n < z + \vep \} & =
& \sum_{ z \leq \ovl{x}_n < z + \vep } \; \prod_{i=1}^n f ( x_i ) \nonumber\\
& \geq & \sum_{ z \leq \ovl{x}_n < z + \vep  } \exp \li ( n \ro ( \ovl{x}_n ) - o(n) \ri ) \times \prod_{i=1}^n g ( x_i ) \nonumber\\
&  \geq & \inf_{y \in [z, z + \vep) } \exp \li ( n \ro ( y )  - o (n) \ri ) \sum_{ z \leq \ovl{x}_n < z + \vep  } \; \prod_{i=1}^n g ( x_i ) \nonumber \\
& = & \inf_{y \in [z, z + \vep) } \exp \li ( n \ro ( y )  - o (n)  \ri ) \; \bb{P}_g \{ z \leq \ovl{X}_n < z + \vep \}. \la{CHb86336} \eel
Combining (\ref{CHb86336}) and (\ref{CLTCH8996}) yields (\ref{limineqb}).

\epf

\bsk

We are now in a position to prove the theorem.  From the proof of Lemma \ref{lemuseful}, we know that \be \la{rest8899} \Pr \{ \ovl{X}_n \geq z
\} \geq \Pr \{ z \leq \ovl{X}_n < z + \vep \} > 0 \ee  if $\vep > 0$ is sufficiently small and $n$ is sufficiently large.  Therefore, in the
sequel, we can restrict $\vep > 0$ to be small enough and $n$ to be large enough so that (\ref{rest8899}) holds.

To show assertion (I) of Theorem \ref{convip}, note that \be \la{id89963}
 \ln \Pr \{ \ovl{X}_n \geq z \} = \ln \Pr \{ z \leq \ovl{X}_n < z +
\vep \} - \ln \f{ \Pr \{ z \leq \ovl{X}_n < z + \vep   \} }{  \Pr \{ \ovl{X}_n \geq z \} }. \ee By virtue of (\ref{id89963}), Lemma
\ref{lemuseful}, and (\ref{limineqa}) of Lemma \ref{lem3388}, we have
 \bee \limsup_{n \to \iy} \f{1}{n} \ln \Pr \{ \ovl{X}_n \geq z   \} & = & \limsup_{n \to \iy} \f{1}{n} \ln \Pr \{ z \leq \ovl{X}_n < z + \vep
\} - \lim_{n \to \iy} \f{1}{n} \ln \f{ \Pr \{ z \leq \ovl{X}_n < z + \vep   \} }{  \Pr \{ \ovl{X}_n
\geq z   \}  } \\
&  = & \limsup_{n \to \iy} \f{1}{n} \ln \Pr \{ z \leq \ovl{X}_n < z + \vep   \} \leq \sup_{y \in [z, z + \vep) } \ro ( y ) \eee for sufficiently
small $\vep > 0$.  By the assumption that $\ro ( y )$ is continuous in $y$, it must be true that \be \la{follw8996}
 \lim_{n \to \iy} \f{1}{n}
\ln \Pr \{ \ovl{X}_n \geq z   \} \leq \limsup_{n \to \iy} \f{1}{n} \ln \Pr \{ \ovl{X}_n \geq z   \} \leq \ro ( z ). \ee It follows from
(\ref{follw8996}) and Cramer-Chernoff theorem  that
\[
\Pr \{ \ovl{X}_n \geq z   \} \leq \exp \li ( n   \lim_{m \to \iy} \f{1}{m} \ln \Pr \{ \ovl{X}_m \geq z   \}  \ri ) \leq \exp ( n \ro ( z )).
\]
This completes the proof of the assertion (I) of Theorem \ref{convip}.

To show assertion (II) of Theorem \ref{convip},  making use of (\ref{id89963}), Lemma \ref{lemuseful}, and (\ref{limineqb}) of Lemma
\ref{lem3388}, we have \bee \liminf_{n \to \iy} \f{1}{n} \ln \Pr \{ \ovl{X}_n \geq z   \} & = & \liminf_{n \to \iy} \f{1}{n} \ln \Pr \{ z \leq
\ovl{X}_n < z + \vep \} - \lim_{n \to \iy} \f{1}{n} \ln \f{ \Pr \{ z \leq \ovl{X}_n < z + \vep   \} }{ \Pr \{ \ovl{X}_n
\geq z   \}  } \\
&  = & \liminf_{n \to \iy} \f{1}{n} \ln \Pr \{ z \leq \ovl{X}_n < z + \vep   \} \geq \inf_{y \in [z, z + \vep) } \ro ( y )  \eee for
sufficiently small $\vep > 0$.  By the assumption that $\ro ( y )$ is continuous in $y$, it must be true that \be \la{follw8996b}
 \lim_{n \to \iy} \f{1}{n}
\ln \Pr \{ \ovl{X}_n \geq z   \} \geq \liminf_{n \to \iy} \f{1}{n} \ln \Pr \{ \ovl{X}_n \geq z   \} \geq \ro ( z ). \ee Combining
(\ref{follw8996}) and (\ref{follw8996b}) yields $\lim_{n \to \iy} \f{1}{n} \ln \Pr \{ \ovl{X}_n \geq z   \} = \ro ( z )$. This completes the
proof of the assertion (II) of Theorem \ref{convip}.

The proof of assertion (III) is similar to that of assertion (I).  The proof of assertion (IV) is similar to that of assertion (II).

To show assertion (V), note that $\prod_{i=1}^n f ( X_i ) \leq \exp \li ( n \ro ( \ovl{X}_n ) \ri ) \prod_{i=1}^n g ( X_i )$,  as a consequence
of the assumption that $\ro (.)$ is a concave function such that $ f(X) \leq \exp (\ro (X) ) \; g(X)$. Hence, applying the established
assertions (I) and (III), we have that assertion (V) is true.

To show assertion (VI), note that $\prod_{i=1}^n f ( X_i ) = \exp \li ( n \ro ( \ovl{X}_n ) \ri ) \prod_{i=1}^n g ( X_i )$,  as a consequence of
the assumption that $\ro (.)$ is a linear function such that $ f(X) = \exp (\ro (X) ) \; g(X)$. Hence, applying the established assertions (II)
and (IV), we have that assertion (VI) is true.  This completes the proof of the theorem.

\section{Proof of Theorem \ref{Hyperfirst}} \la{Hyperfirstapp}

Let the pmf $\Pr \{X = x \}$ defined by (\ref{hyperCDF}) be denoted by $f(x)$.  Let $g(x, \vse)$ denote the pmf $\Pr \{X = x \}$ defined by
(\ref{hyperCDF}) with $R$ replaced by $\vse \in \{0, 1, \cd, N \}$.  Define
\[
\mscr{M} (x, \vse) = \f{ f(x) }{ g(x, \vse) }
\]
for $\vse \in \{0, 1, \cd, N \}$ and $x \in \{0, 1, \cd, n \}$ such that $\vse + n - N \leq x \leq \vse$.  Note that \be \la{defM8963}
 \mscr{M} (x, \vse)  = \f{ \bi{ R }{x} \bi{N - R } {n -
x} } { \bi{ \vse }{x} \bi{N - \vse } {n - x}  }. \ee Define \[ \mscr{A} = \{ \vse \in \bb{Z}^+: r \leq \vse \leq R \}, \qqu \mscr{B} = \{ \vse
\in \bb{Z}^+: R \leq \vse \leq r + N - n \}
\]
for $r \in \{0, 1, \cd, n \}$ and $R \in \{0, 1, \cd, N \}$ such that $R + n - N \leq r \leq R$.    Clearly, both $\mscr{A}$ and $\mscr{B}$ are
nonempty as a consequence of $R + n - N \leq r \leq R$.  We need to establish some preliminary results.

\beL \la{pre89663} Assume that $r \in \{0, 1, \cd, n \}, \; R \in \{0, 1, \cd, N \}$ and that $R + n - N \leq r \leq R$. Then, $f(X) \;
\bb{I}_{\{X \leq r \}} \leq \mscr{M} (r, \vse) g(X, \vse)$ holds for all $\vse \in \mscr{A}$.   Similarly,  $f(X) \; \bb{I}_{\{X \geq r \}} \leq
\mscr{M} (r, \vse) g(X, \vse)$ holds for all $\vse \in \mscr{B}$.

\eeL

\bpf

Note that $\mscr{M} (r, \vse)$ is well-defined for $\vse \in \mscr{A}$ and $r \in \{0, 1, \cd, n \}$ such that $R + n - N \leq r \leq R$. It is
easy to see that $f(x) \; \bb{I}_{\{x \leq r \}} \leq \mscr{M} (r, \vse) g(x, \vse)$ holds  for $x \in \{0, 1, \cd, r \}$ such that $f(x) = 0$.
To show the first assertion, it suffices to show that $f(x) \; \bb{I}_{\{x \leq r \}} \leq \mscr{M} (r, \vse) g(x, \vse)$ holds for $\vse \in
\mscr{A}$ and $x \in \{0, 1, \cd, r \}$ such that $f(x) > 0$. Clearly, $f(x) > 0$ implies $R + n - N \leq x \leq R$. Hence, for $\vse \in
\mscr{A}$ and $x \in \{0, 1, \cd, r \}$ such that $f(x) > 0$, it must be true that $R + n - N \leq x \leq r \leq \vse \leq R$ and $\mscr{M} (x,
\vse) > 0$. Making use of (\ref{defM8963}), we can verify that
\[
\f{ \mscr{M} (x + 1, \vse) }{ \mscr{M} (x, \vse) } =  \f{ (R - x) ( N - \vse - n + x + 1) } { (\vse - x) (N - R - n + x + 1) } \geq 1
\]
for $\vse \in \mscr{A}$ and $x \in \{0,1, \cd, r \}$ such that $R + n - N \leq x < \vse$.  This implies that for $\vse \in \mscr{A}$, $\mscr{M}
(x, \vse)$ is increasing with respect to $x \in \{0,1, \cd, r \}$ such that $R + n - N \leq x < \vse$. Consequently, $\f{  f(x) }{g(x, \vse) } =
\mscr{M} (x, \vse) \leq \mscr{M} (r, \vse)$ holds for $\vse \in \mscr{A}$ and $x \in \{0, 1, \cd, r \}$ such that $f(x) > 0$.  Thus, we have
shown that for $r \in \{0, 1, \cd, n \}$ such that $R + n - N \leq r \leq R$, the inequality $f(X) \; \bb{I}_{\{X \leq r\}} \leq \mscr{M} (r,
\vse) g(X, \vse)$ holds for all $\vse \in \mscr{A}$.

The second assertion can be shown in a similar manner. Note that $\mscr{M} (r, \vse)$ is well-defined for $\vse \in \mscr{B}$ and $r \in \{0, 1,
\cd, n \}$ such that $R + n - N \leq r \leq R$. It is easy to see that $f(x) \; \bb{I}_{\{x \geq r \}} \leq \mscr{M} (r, \vse) g(x, \vse)$ holds
for $x \in \{r, r+1, \cd, n \}$ such that $f(x) = 0$.  For $\vse \in \mscr{B}$ and $x \in \{r, r+1, \cd, n \}$ such that $f(x) > 0$, it must be
true that $R + n - N \leq \vse + n - N \leq r \leq x \leq R \leq \vse$ and $\mscr{M} (x, \vse)
> 0$.  Making use of (\ref{defM8963}), we can verify that $\f{ \mscr{M} (x + 1, \vse) }{ \mscr{M} (x, \vse) } \leq 1$ for $\vse \in \mscr{B}$
and $x \in \{r, r+1, \cd, n \}$ such that $x < R$.
 Consequently, $\f{  f(x) }{g(x, \vse) } = \mscr{M} (x, \vse) \leq \mscr{M} (r, \vse)$ holds for $\vse \in \mscr{B}$ and $x \in \{r, r+1, \cd, n
\}$ such that $f(x) > 0$. It follows that for $r \in \{0, 1, \cd, n \}$ such that $R + n - N \leq r \leq R$, the inequality
 $f(X) \; \bb{I}_{\{X \geq r \}} \leq \mscr{M} (r, \vse) g(X, \vse)$ holds for all $\vse \in \mscr{B}$.  This completes the proof of the lemma.

\epf

Applying  Theorem \ref{ThM888} and Lemma \ref{pre89663}, we have the following results.

\beL

\la{gen7896}  Assume that $r \in \{0, 1, \cd, n \}, \; R \in \{0, 1, \cd, N \}$ and that $R + n - N \leq r \leq R$.  Then,

 \bee &  & \Pr \{ X \leq r \} \leq \f{ \bi{ R }{r} \bi{N - R } {n - r} } { \bi{ \vse }{r} \bi{N - \vse } {n - r}  }
\qu \tx{for $\vse \in \mscr{A}$},\\
&   & \Pr \{ X \geq r \} \leq \f{ \bi{ R }{r} \bi{N - R } {n - r} } { \bi{ \vse }{r} \bi{N - \vse } {n - r}  } \qu \tx{for $\vse \in \mscr{B}$}.
\eee
 \eeL

\beL

\la{68896} Let $\wh{R} = \min \li \{ N, \li \lf (N + 1) \f{ r  }{ n } \ri \rf \ri \}$.  Then,  $\wh{R} \in \mscr{A}$ provided that $\f{r}{n}
\leq \f{R}{N}$.  Similarly,  $\wh{R} \in \mscr{B}$ provided  that $\f{r}{n} \geq \f{R}{N}$.

\eeL

\bpf

First, we shall show $\wh{R} \in \mscr{A}$ under the assumption that $\f{r}{n} \leq \f{R}{N}$.  Clearly, $r \leq \wh{R}$.  To show $\wh{R} \in
\mscr{A}$, it remains to show $\wh{R} \leq R$ by considering the cases that $r = n$ and $0 \leq r < n$.  In the case of $r = n$, we have $\wh{R}
= R = N$. In the case of $0 \leq r < n$, we have $\wh{R} = \li \lf (N + 1) \f{ r }{ n } \ri \rf$ and $N r + r < n R + n$. Hence, $(N + 1) \f{ r
}{ n } < R + 1$  and it follows that $\wh{R} \leq R$.  This proves that $\wh{R} \in \mscr{A}$ provided that $\f{r}{n} \leq \f{R}{N}$.

Next, we shall show  $\wh{R} \in \mscr{B}$ under the assumption  that $\f{r}{n} \geq \f{R}{N}$.  We need to consider the cases that $r = n$ and
$0 \leq r < n$.  In the case of $r = n$, we have $\wh{R} = r + N - n = N \geq R$. In the case of $0 \leq r < n$, we have $\wh{R} = \li \lf (N +
1) \f{ r }{ n } \ri \rf$ and $N r + r \geq n R$, which imply that $\wh{R} \geq R$.  Moreover, we observe that
\[
(N + 1) \f{ r  }{ n } - ( r + N - n + 1 ) = - \li (1 - \f{r}{n} \ri ) (N + 1 - n) < 0
\]
and hence $\wh{R} = \li \lf (N + 1) \f{ r  }{ n } \ri \rf \leq r + N - n$.

\epf

Finally,  Theorem \ref{Hyperfirst} follows from Lemmas \ref{gen7896} and \ref{68896}.

\section{Proof of Theorem  \ref{ThmPos}} \la{ThmPosapp}

Define $Y = \sum_{i = 1}^n X_i$ and $\se = n \lm$.  Using the property of parameter additivity  of the generalized Poisson distribution, we have
that $Y$ possesses the following distribution: \be \la{genPos89963}
 f(y) = \Pr
\{ Y = y \} = \f{ \se (\se + y \al)^{y - 1} e^{-\se - y \al}  }{ y! }, \qqu y = 0, 1, 2, \cd. \ee It can be shown that the mean of $Y$ is $\mu =
\bb{E} [ Y ] = \f{\se}{1 - \al}$. Let $\al \in (0, 1)$ be fixed.   Define
\[
 g(y, \vse) = \f{ \vse (\vse + y \al)^{y - 1} e^{-\vse - y \al}  }{ y! }, \qqu y = 0, 1, 2, \cd
\]
and
\[
\mscr{M} (y, \vse) = \f{ f(y) }{ g(y, \vse) } = \f{ \se e^{-\se}  }{ \vse e^{-\vse} }  \li ( \f{ \se + y \al }{ \vse + y \al } \ri )^{y - 1} =
\f{ \se e^{-\se}  }{ \vse e^{-\vse} }  \exp ( h (y, \vse) ),
\]
where $h(y, \vse) = (y-1) \ln \f{ \se + y \al }{ \vse + y \al }$.  Straightforward but tedious calculation shows that
\[
\f{ \pa \; h(y, \vse) }{\pa y} =  \f{\al (\vse - \se ) (y-1)}{ (\se + y \al)(\vse + y \al) } + \ln \f{ \se + y \al }{ \vse + y \al }
\]
 and
\[
\f{ \pa^2 \; h(y, \vse) }{\pa y^2} =   \f{\al (\vse - \se )}{ (\se + y \al)(\vse + y \al) } \li ( \f{\se +  \al}{ \se + y \al } + \f{\vse +
\al}{ \vse + y \al } \ri ).
\]
Let $0 < \vse \leq \se$. Observing that $\lim_{y \to \iy} \f{ \pa \; h(y, \vse) }{\pa y} = 0$ and $\f{ \pa^2 \; h(y, \vse) }{\pa y^2}  < 0$ for
$y \geq 0$, we have that $\f{ \pa \; h(y, \vse) }{\pa y} > 0$ for $y \geq 0$. Thus, $h(y, \vse)$ is increasing with respect to $y \geq 0$. This
shows that for $0 < \vse \leq \se$, $\mscr{M} (y, \vse)$ is increasing with respect to $y \geq 0$.  It follows that \[ f(Y) \; \bb{I}_{ \{ Y
\leq y \} } \leq \mscr{M}(y, \vse) g(Y, \vse) \qqu \tx{for $0 < \vse \leq \se$}.
\]
Therefore, invoking Theorem \ref{ThM888}, we have \be \la{Pos89638}
 \Pr \{ \ovl{X}_n \leq z \} = \Pr \{ Y \leq n z \} \leq \inf_{\vse \leq \se}
\f{ \se e^{-\se}  }{ \vse e^{-\vse} }  \li ( \f{ \se + n z \al }{ \vse + n z \al } \ri )^{n z - 1}. \ee To show (\ref{Pos8899a}), we use the
method of moments. Specifically, we seek $\vse \in (0, \se]$ such that $\bb{E}_\vse [Y] = n z$, that is,  $\f{\vse}{1 - \al} = n z$, from which
we obtain $\vse = n z (1 - \al)$. As a consequence of $z \leq \f{\lm}{1 - \al}$, such $\vse$ is no greater than $\se$.  It follows from
(\ref{Pos89638}) that
\[
\Pr \{ \ovl{X}_n \leq z \} \leq \f{ \se e^{-\se}  }{ \vse e^{-\vse} }  \li ( \f{ \se + n z \al }{ \vse + n z \al } \ri )^{n z - 1}
\]
with $\vse = n z (1 - \al)$. Substituting $\vse$ into the right side of this inequality yields (\ref{Pos8899a}).

To show (\ref{Pos8899c}), we use the method of maximum likelihood to find the tightest bound.  For this purpose, we define $L(z, \vse) = \vse -
\ln \vse - (n z - 1) \ln ( \vse + n z \al )$ and attempt to  minimize the function $L(z, \vse)$ subject to $0 < \vse \leq \se$. Note that
$\f{\pa L (z, \vse) }{\pa \vse} = 1 - \f{1}{\vse} - \f{n z - 1}{\vse + n z \al}$ and
\[
\f{\pa^2 L (z, \vse) }{\pa \vse^2} = \f{1}{\vse^2} - \f{1}{(\vse + n z \al)^2} + \f{n z}{(\vse + n z \al)^2} >  0.
\]
To make $\f{\pa L (z, \vse) }{\pa \vse} = 0$, it suffices to have $\vse^2 - (1 - \al) n z \vse - n z \al = 0$. This equation has a unique
nonnegative solution, which is
\[
\vse = \f{1}{2} \li [ (1 - \al) n z + \sq{ [(1 - \al) n z]^2 + 4 n z \al } \ri ] = n \nu.
\]
It can be checked that $\f{\pa L (z, \vse) }{\pa \vse} |_{\vse = \se} = 1 - \f{1}{\se} - \f{n z - 1}{\se + n z \al} \geq 0$ provided that $z
\leq \f{\lm}{1 - \al + \f{\al}{n \lm} }$.  Hence, $n \nu \leq \se$ and it follows from (\ref{Pos89638}) that \[ \Pr \{ \ovl{X}_n \leq z \}  \leq
\f{ \se e^{-\se} }{ n \nu e^{- n \nu} }  \li ( \f{ \se + n z \al }{ n \nu + n z \al } \ri )^{n z - 1} = \f{\lm (\nu + z \al)}{\nu (\lm + z \al
)} \li [ \li ( \f{ \lm + z \al }{ \nu + z \al } \ri )^z \f{ e^\nu  }{ e^\lm }  \ri ]^{n}
\]
for $0 < z \leq \f{\lm}{1 - \al + \f{\al}{n \lm} }$. This establishes (\ref{Pos8899c}).

To show (\ref{Pos8899b}) and (\ref{Pos8899d}),  let $\vse \geq \se$. Noting that $\lim_{y \to \iy} \f{ \pa \; h(y, \vse) }{\pa y} = 0$ and $\f{
\pa^2 \; h(y, \vse) }{\pa y^2}  > 0$ for $y \geq 0$, we have that $\f{ \pa \; h(y, \vse) }{\pa y} < 0$ for $y \geq 0$.  Thus, $h (y, \vse)$ is
decreasing for $y \geq 0$. This shows that for $\vse \geq \se$, $\mscr{M} (y, \vse)$ is decreasing with respect to $y \geq 0$.   It follows that
\[ f(Y) \; \bb{I}_{ \{ Y \geq y \} } \leq \mscr{M}(y, \vse) g(Y, \vse) \qqu \tx{for $\vse \geq \se$}.
\] Hence, invoking Theorem \ref{ThM888}, we
have \be \la{P889966}
 \Pr \{ \ovl{X}_n \geq z \} = \Pr \{ Y \geq n z \} \leq \inf_{\vse \geq \se} \f{ \se e^{-\se} }{ \vse e^{-\vse} }  \li (
\f{ \se + n z \al }{ \vse + n z \al } \ri )^{n z - 1}. \ee To show (\ref{Pos8899b}), we use the method of moments. Specifically, we seek $\vse
\in [\se, \iy)$ such that $\bb{E}_\vse [Y] = n z$, that is,  $\f{\vse}{1 - \al} = n z$, from which we obtain $\vse = n z (1 - \al)$. As a
consequence of $z \geq \f{\lm}{1 - \al}$, such $\vse$ is no less than $\se$.  It follows from (\ref{P889966}) that
\[
\Pr \{ \ovl{X}_n \geq z \} \leq \f{ \se e^{-\se}  }{ \vse e^{-\vse} }  \li ( \f{ \se + n z \al }{ \vse + n z \al } \ri )^{n z - 1}
\]
with $\vse = n z (1 - \al)$. Substituting $\vse$ into the right side of this inequality yields (\ref{Pos8899b}).

To show (\ref{Pos8899d}), we use the method of maximum likelihood to find the tightest bound.  For this purpose, it suffices to minimize the
function $L(z, \vse)$ subject to $\vse \geq \se$.  It can be verified that $\f{\pa L (z, \vse) }{\pa \vse} |_{\vse = \se} = 1 - \f{1}{\se} -
\f{n z - 1}{\se + n z \al} \leq 0$ provided that $z \geq \f{\lm}{1 - \al + \f{\al}{n \lm} }$. Hence, $n \nu \geq \se$ and it follows from
(\ref{P889966}) that \[ \Pr \{ \ovl{X}_n \geq z \} \leq \f{ \se e^{-\se}  }{ n \nu e^{- n \nu} }  \li ( \f{ \se + n z \al }{ n \nu + n z \al }
\ri )^{n z - 1} = \f{\lm (\nu + z \al)}{\nu (\lm + z \al )} \li [ \li (  \f{ \lm + z \al }{ \nu + z \al } \ri )^z \f{ e^\nu  }{ e^\lm } \ri
]^{n}
\]
for $z \geq \f{\lm}{1 - \al + \f{\al}{n \lm} }$. This establishes (\ref{Pos8899d}).

In the special case that $\al = 0$, the generalized Poisson distribution  reduces to the standard Poisson distribution in the exponential family
with pmf
\[
f(x) = \Pr \{ X = x \} = \f{ \lm^{x} e^{-\lm} } { x! } = v (x) \exp \left( \eta (\lm) u (x) - \ze (\lm) \right), \qqu x \in \{ 0, 1, 2, \cd \}
\]
where
\[
u (x) = x,  \qqu v (x) = \f{1}{x!}, \qqu  \eta (\lm) = \ln \lm, \qqu \ze (\lm) = \lm.
\]
The moment generating function is $\phi(s) = \bb{E} [ e^{s X} ] = e^{- \lm} \exp ( \lm e^s )$. Using $\bb{E} [ X^{\ell} ] =   \li. \f{ d^\ell
\phi (s) } { d s^\ell } \ri |_{s = 0}$ for $\ell = 1, 2, \cd$, we have
\[
\bb{E} [ X]  = \lm, \qu  \bb{E} [ X^2 ] = \lm (1 + \lm), \qu \bb{E} [ X^3 ]  = \lm^2 + (1 + \lm)^2 \lm, \qu \bb{E} [ X^4 ]  = \lm^2 + (1 + \lm)
[3 \lm^2 + (1 + \lm)^2 \lm ].
\]
Hence, $\bb{E} [ |X - \lm|^2 ] = \lm$ and $\bb{E} [ |X - \lm|^4 ]  =  \bb{E} [ X^4 - 4 X^3 \lm + 6 X^2 \lm^2 - 4 X \lm^3 + \lm^4 ] = \lm (3 \lm
+ 1)$.  Using these facts and Lyapunov's inequality, we have \[ \f { \bb{E} [ |X - \lm|^3 ] } { \bb{E}^{\f{3}{2}} [ |X - \lm|^2 ] } < \f {
\bb{E}^{\f{3}{4}} [ |X - \lm|^4 ] } { \bb{E}^{\f{3}{2}} [ |X - \lm|^2 ] } = \li ( 3 + \f{1}{\lm} \ri )^{3 \sh 4}. \] Making use of this
inequality, the assertions (III) and (IV) of Theorem \ref{univexp}, we can establish (\ref{Pos88a}) and (\ref{Pos88b}).

\section{Proof of Theorem \ref{GenmulHyper} } \la{GenmulHyperapp}

We need to define some quantities. Define
\[
\bs{\mcal{C}} = [\mcal{C}_0, \mcal{C}_1, \cd, \mcal{C}_\ka]^\top, \qqu \wh{\bs{\mcal{C}}} = [\wh{\mcal{C}}_0, \wh{\mcal{C}}_1, \cd,
\wh{\mcal{C}}_\ka]^\top,
\]
\[
\mscr{X} = \li \{ [x_0, x_1, \cd, x_\ka]^\top: x_i \in \bb{Z}^+, \; i = 0, 1, \cd, \ka \; \tx{and} \; \sum_{i=0}^\ka x_i = n \ri \}
\]
and
\[
\bs{\varTheta} = \li \{ [\vse_0, \vse_1, \cd, \vse_\ka]^\top:  \; \vse_i \in \bb{R}, \; \f{\vse_i}{\mcal{N}}
> 0 \; \tx{for} \; i = 0, 1, \cd, \ka \; \tx{and} \; \sum_{i=0}^\ka \vse_i = \mcal{N} \ri \}.
\]
Let $g (\bs{x}, \bs{\vse})$, where $\bs{x} = [x_0, x_1, \cd, x_\ka]^\top \in \mscr{X}$ and $\bs{\vse} = [\vse_0, \vse_1, \cd, \vse_\ka]^\top \in
\bs{\varTheta}$,  denote the pmf defined by (\ref{GINHY}) with $\mcal{C}_i$ replaced by $\vse_i$ for $i = 0, 1, \cd, \ka$. Define $\mscr{M}
(\bs{x}, \bs{\vse}) = \f{ f(\bs{x}) }{ g (\bs{x}, \bs{\vse}) }$ for $\bs{x} \in \mscr{X}$ and $\bs{\vse}  \in \bs{\varTheta}$  such that $g
(\bs{x}, \bs{\vse})
> 0$.  Define $\mscr{S} = \{ \bs{x} \in \mscr{X}: \; f(\bs{x}) > 0 \}$
and
\[
\mscr{A} = \{ \bs{\vse} \in \bs{\varTheta}: | \bs{\vse} | \bs{\prec} | \bs{\mcal{C}} |, \; \; g(\bs{z}, \bs{\vse}) > 0  \}, \qqu \mscr{B} = \{
\bs{\vse} \in \bs{\varTheta}: | \bs{\vse} | \bs{\succ} | \bs{\mcal{C}} |, \; \; g(\bs{z}, \bs{\vse}) > 0  \}.
\]

\beL \la{pre89663mul} Assume that $f(\bs{z}) > 0$.  Then,  $f(\bs{X}) \bb{I}_{\{\bs{X} \bs{\prec} \bs{z} \}} \leq \mscr{M} (\bs{z}, \bs{\vse})
g(\bs{X}, \bs{\vse})$ holds for all $\bs{\vse} \in \mscr{A}$. Similarly,  $f(\bs{X}) \bb{I}_{\{\bs{X} \bs{\succ} \bs{z} \}} \leq \mscr{M}
(\bs{z}, \bs{\vse}) g(\bs{X}, \bs{\vse})$ holds for all $\bs{\vse} \in \mscr{B}$.

\eeL

\bpf First, we shall show the first assertion that $f(\bs{X}) \bb{I}_{\{\bs{X} \bs{\prec} \bs{z} \}} \leq \mscr{M} (\bs{z}, \bs{\vse}) g(\bs{X},
\bs{\vse})$ holds for all $\bs{\vse} \in \mscr{A}$.   Since $f(\bs{z}) = g(\bs{z}, \bs{\mcal{C}})$, the set $\mscr{A}$ is nonempty  as a
consequence of the assumption that $f(\bs{z}) > 0$.  Hence, $\mscr{M} (\bs{z}, \bs{\vse})$ is well-defined for $\bs{z} \in \mscr{S}$ and
$\bs{\vse} \in \mscr{A}$.  Note that $f(\bs{x}) \bb{I}_{\{\bs{x} \bs{\prec} \bs{z} \}} \leq \mscr{M} (\bs{z}, \bs{\vse}) g(\bs{x}, \bs{\vse})$
holds for $\bs{\vse} \in \mscr{A}$ and $\bs{x} \in \mscr{X}$ such that $f(\bs{x}) = 0$. To show the first assertion,  it suffices to show that
$f(\bs{x}) \bb{I}_{\{\bs{x} \bs{\prec} \bs{z} \}} \leq \mscr{M} (\bs{z}, \bs{\vse}) g(\bs{x}, \bs{\vse})$ holds for $\bs{\vse} \in \mscr{A}$ and
$\bs{x} \in \mscr{S}$.  We need to consider two cases as follows.

Case (I): $\mcal{N} > 0$.

Case (II): $\mcal{N} < 0$.

In  Case (I), as a consequence of $\mcal{N} > 0$, we have $\vse_i > 0, \; \mcal{C}_i > 0$ for $i = 0, 1, \cd, \ka$.  Consider $\bs{x} \in
\mscr{S}$ such that $\bs{x} \bs{\prec} \bs{z}$.  Note that $\prod_{i=1}^k \prod_{\ell = 1}^{x_i} (\vse_i - \ell + 1) \geq \prod_{i=1}^k
\prod_{\ell = 1}^{z_i} (\vse_i - \ell + 1) > 0 $ because of $\bs{x} \bs{\prec} \bs{z}$ and $g(\bs{z}, \bs{\vse}) > 0$.  Moreover, $\prod_{\ell =
1}^{x_0} (\vse_0 - \ell + 1) \geq \prod_{\ell = 1}^{x_0} (\mcal{C}_0 - \ell + 1) > 0$ because of $\vse_0 \geq \mcal{C}_0$ and $f(\bs{x})
> 0$.  This implies that $g(\bs{x}, \bs{\vse}) > 0$ for $\bs{\vse} \in \mscr{A}$ and $\bs{x} \in \mscr{S}$ such that $\bs{x} \bs{\prec} \bs{z}$.  Since $f(\bs{x}) > 0$, we have that
\[
\mscr{M} (\bs{x}, \bs{\vse}) = \f{ f(\bs{x}) }{ g (\bs{x}, \bs{\vse}) } = \prod_{i=0}^\ka \f{ \bi{\mcal{C}_i}{x_i}  }{ \bi{\vse_i}{x_i} } =
\prod_{i=0}^k  \prod_{\ell = 1}^{x_i}  \f{ \mcal{C}_i - \ell  + 1  }{ \vse_i - \ell + 1 }
\]
is well-defined and positive for $\bs{\vse} \in \mscr{A}$ and $\bs{x} \in \mscr{S}$ such that $\bs{x} \bs{\prec} \bs{z}$. Consequently, we have
a meaningful ratio \be \la{case899}
 \f{ \mscr{M} (\bs{z}, \bs{\vse})  }{ \mscr{M} (\bs{x}, \bs{\vse})  } =  \li ( \prod_{\ell = z_0 + 1}^{x_0}  \f{ \vse_0 - \ell + 1 }{ \mcal{C}_0 - \ell + 1 }
\ri ) \li ( \prod_{i=1}^k  \prod_{\ell = x_i + 1}^{z_i}  \f{ \mcal{C}_i - \ell  + 1  }{ \vse_i - \ell + 1 } \ri ) \qqu  \tx{for $\bs{\vse} \in
\mscr{A}$ and $\bs{x} \in \mscr{S}$ such that} \; \bs{x} \bs{\prec} \bs{z}. \ee
Since $\vse_0 > \mcal{C}_0$ and $\vse_i < \mcal{C}_i$ for $i =
1, \cd, \ka$, it follows from (\ref{case899}) that $\mscr{M} (\bs{x}, \bs{\vse} ) \leq \mscr{M} (\bs{z}, \bs{\vse})$ holds for $\bs{\vse} \in
\mscr{A}$ and $\bs{x} \in \mscr{S}$ such that $\bs{x} \bs{\prec} \bs{z}$.

In  Case (II), as a consequence of $\mcal{N} < 0$, we have $\vse_i < 0, \; \mcal{C}_i < 0$ for $i = 0, 1, \cd, \ka$.  Consider $\bs{x} \in
\mscr{X}$ such that $\bs{x} \bs{\prec} \bs{z}$.  Note that for $\bs{\vse} \in \mscr{A}$ and $\bs{x} \in \mscr{X}$,
\[
\mscr{M} (\bs{x}, \bs{\vse}) =  \prod_{i=0}^k  \prod_{\ell = 1}^{x_i} \f{ \mcal{C}_i - \ell  + 1 }{ \vse_i - \ell + 1  } = \prod_{i=0}^k
\prod_{\ell = 1}^{x_i} \f{  |\mcal{C}_i| + \ell  - 1 }{ |\vse_i|  + \ell - 1 }
\]
is well-defined and positive.  Moreover, for $\bs{\vse} \in \mscr{A}$ and $\bs{x} \in \mscr{X}$ such that $\bs{x} \bs{\prec} \bs{z}$, we have a
meaningful ratio \be \la{case899b}
 \f{ \mscr{M}
(\bs{z}, \bs{\vse})  }{ \mscr{M} (\bs{x}, \bs{\vse})  } =  \li ( \prod_{\ell = z_0 + 1}^{x_0}  \f{ |\vse_0 | + \ell - 1 }{ |\mcal{C}_0| + \ell -
1 } \ri ) \li ( \prod_{i=1}^k  \prod_{\ell = x_i + 1}^{z_i}  \f{ | \mcal{C}_i | + \ell  - 1  }{ | \vse_i | + \ell - 1 } \ri ). \ee Since
$|\vse_0| > | \mcal{C}_0 |$ and $| \vse_i | < | \mcal{C}_i |$ for $i = 1, \cd, \ka$, it follows from (\ref{case899b}) that $\mscr{M} (\bs{x},
\bs{\vse} ) \leq \mscr{M} (\bs{z}, \bs{\vse})$ holds for $\bs{\vse} \in \mscr{A}$ and $\bs{x} \in \mscr{X}$ such that $\bs{x} \bs{\prec}
\bs{z}$. Therefore, in both cases, we have that
\[
\f{ f(\bs{x})}{g(\bs{x}, \bs{\vse})} = \mscr{M} (\bs{x}, \bs{\vse} )  \leq \mscr{M} (\bs{z}, \bs{\vse} )
\]
holds for $\bs{\vse} \in \mscr{A}$ and $\bs{x} \in \mscr{S}$ such that $\bs{x} \bs{\prec} \bs{z}$. Thus, we have shown that $f(\bs{x})
\bb{I}_{\{\bs{x} \bs{\prec} \bs{z} \}} \leq \mscr{M} (\bs{z}, \bs{\vse}) g(\bs{x}, \bs{\vse})$ holds for $\bs{\vse} \in \mscr{A}$ and  $\bs{x}
\in \mscr{S}$.  This completes the proof of the first assertion of the lemma.

Next, we shall show the second assertion that $f(\bs{X}) \bb{I}_{\{\bs{X} \bs{\succ} \bs{z} \}} \leq \mscr{M} (\bs{z}, \bs{\vse}) g(\bs{X},
\bs{\vse})$ holds for all $\bs{\vse} \in \mscr{B}$.  Since $f(\bs{z}) = g(\bs{z}, \bs{\mcal{C}})$, the set $\mscr{B}$ is nonempty  as a
consequence of the assumption that $f(\bs{z}) > 0$.  Hence, $\mscr{M} (\bs{z}, \bs{\vse})$ is well-defined for $\bs{z} \in \mscr{S}$ and
$\bs{\vse} \in \mscr{B}$.  Note that $f(\bs{x}) \bb{I}_{\{\bs{x} \bs{\succ} \bs{z} \}} \leq \mscr{M} (\bs{z}, \bs{\vse}) g(\bs{x}, \bs{\vse})$
holds for $\bs{\vse} \in \mscr{B}$ and $\bs{x} \in \mscr{X}$ such that $f(\bs{x}) = 0$. To show the second assertion,  it suffices to show that
$f(\bs{x}) \bb{I}_{\{\bs{x} \bs{\succ} \bs{z} \}} \leq \mscr{M} (\bs{z}, \bs{\vse}) g(\bs{x}, \bs{\vse})$ holds for $\bs{\vse} \in \mscr{B}$ and
$\bs{x} \in \mscr{S}$.  We need to consider two cases as follows.

Case (I): $\mcal{N} > 0$.

Case (II): $\mcal{N} < 0$.

In  Case (I), as a consequence of $\mcal{N} > 0$, we have $0 < \vse_0 \leq \mcal{C}_0$ and $\vse_i \geq \mcal{C}_i > 0$ for $i = 1, \cd, \ka$.
Consider $\bs{x} \in \mscr{S}$ such that $\bs{x} \bs{\succ} \bs{z}$.  Note that $\prod_{i=1}^k \prod_{\ell = 1}^{x_i} (\vse_i - \ell + 1) \geq
\prod_{i=1}^k \prod_{\ell = 1}^{x_i} (\mcal{C}_i - \ell + 1) > 0 $ because of $\bs{\vse} \bs{\succ} \bs{\mcal{C}}$ and $f (\bs{x}) > 0$.
Moreover, $\prod_{\ell = 1}^{x_0} (\vse_0 - \ell + 1) \geq \prod_{\ell = 1}^{z_0} (\vse_0 - \ell + 1) > 0$ because of $\bs{x} \bs{\succ} \bs{z}$
and $g(\bs{z}, \bs{\vse}) > 0$.  This implies that $g(\bs{x}, \bs{\vse}) > 0$ for $\bs{\vse} \in \mscr{B}$ and $\bs{x} \in \mscr{S}$ such that
$\bs{x} \bs{\succ} \bs{z}$. Since $f (\bs{x}) > 0$, we have that
\[
\mscr{M} (\bs{x}, \bs{\vse}) = \f{ f(\bs{x}) }{ g (\bs{x}, \bs{\vse}) } = \prod_{i=0}^\ka \f{ \bi{\mcal{C}_i}{x_i}  }{ \bi{\vse_i}{x_i} } =
\prod_{i=0}^k  \prod_{\ell = 1}^{x_i}  \f{ \mcal{C}_i - \ell  + 1  }{ \vse_i - \ell + 1 }
\]
is well-defined and positive for $\bs{\vse} \in \mscr{B}$ and $\bs{x} \in \mscr{S}$ such that $\bs{x} \bs{\succ} \bs{z}$. Consequently, we have
a meaningful ratio \be \la{case899BB}
 \f{ \mscr{M} (\bs{x}, \bs{\vse})  }{ \mscr{M} (\bs{z}, \bs{\vse})  } =  \li ( \prod_{\ell = x_0 + 1}^{z_0}  \f{ \vse_0 - \ell + 1 }{ \mcal{C}_0 - \ell + 1 }
\ri ) \li ( \prod_{i=1}^k  \prod_{\ell = z_i + 1}^{x_i}  \f{ \mcal{C}_i - \ell  + 1  }{ \vse_i - \ell + 1 } \ri ) \qqu  \tx{for $\bs{\vse} \in
\mscr{B}$ and $\bs{x} \in \mscr{S}$ such that} \; \bs{x} \bs{\succ} \bs{z}. \ee It can be seen from (\ref{case899BB}) that $\mscr{M} (\bs{x},
\bs{\vse} ) \leq \mscr{M} (\bs{z}, \bs{\vse})$ holds for $\bs{\vse} \in \mscr{B}$ and $\bs{x} \in \mscr{S}$ such that $\bs{x} \bs{\succ}
\bs{z}$.

In  Case (II), as a consequence of $\mcal{N} < 0$, we have $\vse_i < 0, \; \mcal{C}_i < 0$ for $i = 0, 1, \cd, \ka$.  Consider $\bs{x} \in
\mscr{X}$ such that $\bs{x} \bs{\succ} \bs{z}$.  Note that for $\bs{\vse} \in \mscr{B}$ and $\bs{x} \in \mscr{X}$,
\[
\mscr{M} (\bs{x}, \bs{\vse}) =  \prod_{i=0}^k  \prod_{\ell = 1}^{x_i} \f{ \mcal{C}_i - \ell  + 1 }{ \vse_i - \ell + 1  } = \prod_{i=0}^k
\prod_{\ell = 1}^{x_i} \f{  |\mcal{C}_i| + \ell  - 1 }{ |\vse_i|  + \ell - 1 }
\]
is well-defined and positive.  Moreover, for $\bs{\vse} \in \mscr{B}$ and $\bs{x} \in \mscr{X}$ such that $\bs{x} \bs{\succ} \bs{z}$, we have a
meaningful ratio \be \la{case899bBB}
 \f{ \mscr{M}
(\bs{x}, \bs{\vse})  }{ \mscr{M} (\bs{z}, \bs{\vse})  } =  \li ( \prod_{\ell = x_0 + 1}^{z_0}  \f{ |\vse_0 | + \ell - 1 }{ |\mcal{C}_0| + \ell -
1 } \ri ) \li ( \prod_{i=1}^k  \prod_{\ell = z_i + 1}^{x_i}  \f{ | \mcal{C}_i | + \ell  - 1  }{ | \vse_i | + \ell - 1 } \ri ). \ee Since
$|\vse_0| < |\mcal{C}_0|$ and $|\vse_i| > |\mcal{C}_i|$ for $i = 1, \cd, \ka$, it follows from (\ref{case899bBB}) that $\mscr{M} (\bs{x},
\bs{\vse} ) \leq \mscr{M} (\bs{z}, \bs{\vse})$ holds for $\bs{\vse} \in \mscr{B}$ and $\bs{x} \in \mscr{X}$ such that $\bs{x} \bs{\succ}
\bs{z}$. Therefore, in both cases, we have that
\[
\f{ f(\bs{x})}{g(\bs{x}, \bs{\vse})} = \mscr{M} (\bs{x}, \bs{\vse} )  \leq \mscr{M} (\bs{z}, \bs{\vse} )
\]
holds for $\bs{\vse} \in \mscr{B}$ and $\bs{x} \in \mscr{S}$ such that $\bs{x} \bs{\succ} \bs{z}$. Thus, we have shown that $f(\bs{x})
\bb{I}_{\{\bs{x} \bs{\succ} \bs{z} \}} \leq \mscr{M} (\bs{z}, \bs{\vse}) g(\bs{x}, \bs{\vse})$ holds for $\bs{\vse} \in \mscr{B}$ and  $\bs{x}
\in \mscr{S}$.  This completes the proof of the second assertion of the lemma.

\epf

Applying  Theorem \ref{ThM888} and Lemma \ref{pre89663mul}, we have the following results.

\beL

\la{generhyp86389}

Assume that $f(\bs{z}) > 0$.  Then,
\[
\Pr \{ \bs{X} \bs{\prec} \bs{z} \}  \leq \prod_{i=0}^\ka \f{ \bi{ \mcal{C}_i  }{ z_i } } { \bi{\vse_i}{z_i}  } \qu \tx{for $\bs{\vse} \in
\mscr{A}$},
\]
\[
\Pr \{ \bs{X} \bs{\succ} \bs{z} \}  \leq \prod_{i=0}^\ka \f{ \bi{ \mcal{C}_i  }{ z_i } } { \bi{\vse_i}{z_i}  } \qu \tx{for $\bs{\vse} \in
\mscr{B}$}.
\]

\eeL

\beL

\la{proper68833}

$g(\bs{z}, \wh{\bs{\mcal{C}}}) > 0$ for $\bs{z} \in \mscr{X}$.  \eeL

\bpf

Recall that $\wh{\bs{\mcal{C}}} = [\wh{\mcal{C}}_0, \wh{\mcal{C}}_1, \cd, \wh{\mcal{C}}_\ka]^\top$ with $\wh{\mcal{C}}_i = \mcal{N} \f{z_i}{n}$
for $i = 0, 1, \cd, \ka$.

In the case of $\mcal{N} < 0$, we have $\wh{\mcal{C}}_i \leq 0$ for $i = 0, 1, \cd, \ka$, which implies that $1 + \wh{\mcal{C}}_i \leq 1$ for $i
= 0, 1, \cd, \ka$. Hence, $g(\bs{z}, \wh{\bs{\mcal{C}}}) > 0$ for $\bs{z} \in \mscr{X}$ if $\mcal{N} < 0$.

In the case of $\mcal{N} > 0$, we have $0 < n < \mcal{N} + 1$ and $\wh{\mcal{C}}_i \geq 0$ for $i = 0, 1, \cd, \ka$. For $i = 0, 1, \cd, \ka$,
we have $z_i < \wh{\mcal{C}}_i + 1 = 1$ if $z_i = 0$; and $ z_i = \f{n}{\mcal{N}} \wh{\mcal{C}}_i < \wh{\mcal{C}}_i + \f{ \wh{\mcal{C}}_i
}{\mcal{N}} = \wh{\mcal{C}}_i + \f{z_i}{n} \leq \wh{\mcal{C}}_i + 1$ if $z_i > 0$. Hence, $g(\bs{z}, \wh{\bs{\mcal{C}}}) > 0$ for $\bs{z} \in
\mscr{X}$ if $\mcal{N} > 0$.

\epf

\beL

\la{more89663}

$|\wh{\bs{\mcal{C}}} | \bs{\prec} |\bs{\mcal{C}}| \Leftrightarrow \bs{z} \bs{\prec} \bs{\mu}$ and $|\wh{\bs{\mcal{C}}} | \bs{\succ}
|\bs{\mcal{C}}| \Leftrightarrow \bs{z} \bs{\succ} \bs{\mu}$.

\eeL

\bpf

In the case of $\mcal{N} > 0$, we have $\mcal{C}_i > 0$ and $\wh{\mcal{C}}_i \geq 0$ for $i = 0, 1, \cd, \ka$. Hence,
\[
|\wh{\bs{\mcal{C}}} | \bs{\prec} |\bs{\mcal{C}}| \Leftrightarrow \wh{\bs{\mcal{C}}} \bs{\prec} \bs{\mcal{C}} \Leftrightarrow \bs{z} \bs{\prec}
\bs{\mu}, \qqu |\wh{\bs{\mcal{C}}} | \bs{\succ} |\bs{\mcal{C}}| \Leftrightarrow \wh{\bs{\mcal{C}}} \bs{\succ} \bs{\mcal{C}} \Leftrightarrow
\bs{z} \bs{\succ} \bs{\mu}.
\]

In the case of $\mcal{N} < 0$, we have $\mcal{C}_i < 0$ and $\wh{\mcal{C}}_i \leq 0$ for $i = 0, 1, \cd, \ka$. Hence,
\[
|\wh{\bs{\mcal{C}}} | \bs{\prec} |\bs{\mcal{C}}| \Leftrightarrow \wh{\bs{\mcal{C}}} \bs{\succ} \bs{\mcal{C}} \Leftrightarrow \bs{z} \bs{\prec}
\bs{\mu}, \qqu |\wh{\bs{\mcal{C}}} | \bs{\succ} |\bs{\mcal{C}}| \Leftrightarrow \wh{\bs{\mcal{C}}} \bs{\prec} \bs{\mcal{C}} \Leftrightarrow
\bs{z} \bs{\succ} \bs{\mu}.
\]

\epf

\bsk

We are now in a position to prove the theorem. Clearly, the assertions (I) and (II) of Theorem \ref{GenmulHyper} follow from Lemmas
\ref{generhyp86389}, \ref{proper68833} and \ref{more89663}.

To show assertion (III) of Theorem \ref{GenmulHyper}, we need to use the following inequalities \be \la{chencina}
 \sq{2 \pi} x^{x - 1 \sh 2}
e^{-x} < \Ga (x) \leq e x^{x - 1 \sh 2} e^{-x}, \qqu x \geq 1. \ee This result is established by Mortici and Chen \cite[page 70, eq.
(1.3)]{Mortici}. Making use of (\ref{chencina}) and the fact that $\Ga(x+1) = x \Ga(x)$, we have \be \la{ChenM}
 \sq{2 \pi x} \; x^x e^{ - x } < \Ga(x+1) \leq e \sq{ x} \; x^x e^{ - x } \qu \tx{for all real number $x \geq 1$}.
\ee Using (\ref{ChenM}) and the assumption that $\mcal{C}_i - z_i \geq 1$ and $\wh{\mcal{C}}_i - z_i \geq 1$ for $i = 0, 1, \cd, k$,  we have
\bee \f{ \bi{\mcal{C}_i}{z_i}  }{ \bi{\wh{\mcal{C}}_i}{z_i} }
& = & \f{ \Ga(\mcal{C}_i + 1) }{  \Ga (\mcal{C}_i - z_i + 1) }  \f{  \Ga (\wh{\mcal{C}}_i - z_i + 1) }{ \Ga(\wh{\mcal{C}}_i + 1) } \\
& < & \f{ e^2  }{ 2 \pi } \; \f{ \mcal{C}_i^{\mcal{C}_i + \f{1}{2}} }{ ( \mcal{C}_i - z_i)^{\mcal{C}_i - z_i + \f{1}{2}} } \;
 \f{ ( \wh{\mcal{C}}_i - z_i)^{\wh{\mcal{C}}_i - z_i + \f{1}{2}} }{ \wh{\mcal{C}}_i^{\wh{\mcal{C}}_i + \f{1}{2}} } \\
& = &  \f{ e^2  }{ 2 \pi } \; \f{ \mcal{C}_i^{\mcal{C}_i + \f{1}{2}} }{ ( \mcal{C}_i - z_i)^{\mcal{C}_i - z_i + \f{1}{2}} } \;
 \f{ ( \mcal{N} \f{z_i}{n}  - z_i)^{\wh{\mcal{C}}_i - z_i + \f{1}{2}} }{ ( \mcal{N} \f{z_i}{n} )^{\wh{\mcal{C}}_i + \f{1}{2}} } \\
 & = &  \f{ e^2  }{ 2 \pi } \; \li (  \f{\mcal{C}_i}{z_i} \ri )^{z_i}  \li ( \f{\mcal{C}_i}{\mcal{C}_i - z_i} \ri )^{ \mcal{C}_i - z_i + \f{1}{2}}
  \li ( \f{n}{\mcal{N}} \ri )^{z_i}  \li ( \f{ \mcal{N} - n }{ \mcal{N} } \ri )^{\wh{\mcal{C}}_i - z_i + \f{1}{2}}   \eee
  for $i = 0, 1, \cd, k$.  Making use of this inequality and the fact that
  \[
\sum_{i = 0}^k z_i = n, \qqu  \sum_{i = 0}^k \li ( \wh{\mcal{C}}_i - z_i + \f{1}{2} \ri ) = \mcal{N} - n + \f{k + 1}{2},
  \]
we have that (\ref{holds889}) holds.

To show assertion (IV) of Theorem \ref{GenmulHyper}, we observe that \be \la{crt8896a}
 \f{1}{\mcal{N}} \ln \li [ \li ( \f{\mcal{N} - n}{\mcal{N}} \ri
)^{\mcal{N} - n} \prod_{i = 0}^k \li ( \f{\mcal{C}_i}{\mcal{C}_i - z_i} \ri )^{ \mcal{C}_i - z_i } \ri ] \to (1 - \al) \ln ( 1 - \al) + \sum_{i
= 0}^k ( \se_i - \al \ba_i ) \ln \f{ \se_i }{ \se_i - \al \ba _i } \ee as $\mcal{N} \to \iy$ under the constraint that $\f{n}{\mcal{N}} \to \al
\in (0, 1)$ and $\f{z_i}{n} \to \ba_i \in (0, 1)$ with  $\se_i > \al \ba_i$ for $i = 0, 1, \cd, k$.  Using the weighted AM-GM inequality, we
have \bel  \prod_{i = 0}^k  \li ( \f{ \se_i }{ \se_i - \al \ba _i } \ri
)^{ \f{ \se_i - \al \ba_i}{1 - \al} } & \leq & \sum_{i = 0}^k  \li ( \f{ \se_i }{ \se_i - \al \ba _i } \ri )  \f{ \se_i - \al \ba_i}{1 - \al} \la{ineq889966}\\
& = & \f{1}{1 - \al}. \nonumber \eel We claim that (\ref{ineq889966}) does not hold with equality.  To prove the claim, note that if
(\ref{ineq889966}) holds with equality,  then there exists a number $\nu$ such that
\[
 \f{ \se_i }{ \se_i - \al \ba _i } = \f{1}{\nu}, \qqu i = 0, 1, \cd, k.
\]
It follows that
\[
\f{1}{1 - \al} = \prod_{i = 0}^k  \li ( \f{ \se_i }{ \se_i - \al \ba _i } \ri )^{ \f{ \se_i - \al \ba_i}{1 - \al} } = \prod_{i = 0}^k  \li (
\f{1}{\nu} \ri )^{ \f{ \se_i - \al \ba_i}{1 - \al} } = \f{1}{\nu},
\]
which implies that $\nu = 1 - \al$ and thus
\[
\f{ \se_i }{ \se_i - \al \ba _i } = \f{1}{1 - \al}, \qqu i = 0, 1, \cd, k.
\]
Consequently, $\se_i = \ba_i$ for $i = 0, 1, \cd, k$.  This contradicts to the assumption that $\se_i \neq \ba_i$ for some $i \in \{0, 1, \cd,
\ka \}$. This proves our claim and thus
\[
\prod_{i = 0}^k  \li ( \f{ \se_i }{ \se_i - \al \ba _i } \ri )^{ \f{ \se_i - \al \ba_i}{1 - \al} } < \f{1}{1 - \al},
\]
which can be written as \be \la{crt8896b}
 (1 - \al) \ln ( 1 - \al) + \sum_{i = 0}^k ( \se_i - \al \ba_i ) \ln \f{ \se_i }{ \se_i - \al \ba _i }
< 0. \ee As a result of (\ref{crt8896a}) and (\ref{crt8896b}), we have that
\[
\li [ \li ( \f{\mcal{N} - n}{\mcal{N}} \ri )^{\mcal{N} - n} \prod_{i = 0}^k \li ( \f{\mcal{C}_i}{\mcal{C}_i - z_i} \ri )^{ \mcal{C}_i - z_i }
\ri ]^{\mcal{N}} \to 0.
\]
This leads to the truth of assertion (IV).  The proof of the theorem is thus completed.

\sect{Proof of Theorem  \ref{GenmulInvHyper}} \la{GenmulInvHyperapp}

We need to define some quantities. Define
\[
\bs{\mcal{C}} = [\mcal{C}_0, \mcal{C}_1, \cd, \mcal{C}_\ka]^\top, \qqu \wh{\bs{\mcal{C}}} = [\wh{\mcal{C}}_0, \wh{\mcal{C}}_1, \cd,
\wh{\mcal{C}}_\ka]^\top,
\]
\[
\mscr{X} = \li \{ [x_0, x_1, \cd, x_\ka]^\top: x_0 = \ga, \; x_i \in \bb{Z}^+, \; i = 1, \cd, \ka \; \tx{and} \; \f{1}{\mcal{N}} \sum_{i=0}^\ka
x_i  < 1  + \f{1}{\mcal{N}} \ri \}
\]
and
\[
\bs{\varTheta} = \li \{ [\vse_0, \vse_1, \cd, \vse_\ka]^\top:  \; \vse_i \in \bb{R}, \; \f{\vse_0}{\mcal{N}} > \f{\ga - 1}{\mcal{N}}, \;
\f{\vse_i}{\mcal{N}}
> 0 \; \tx{for} \; i = 0, 1, \cd, \ka \; \tx{and} \; \sum_{i=0}^\ka \vse_i = \mcal{N} \ri \}
\]
Let $g (\bs{x}, \bs{\vse})$, where $\bs{x} = [x_0, x_1, \cd, x_\ka]^\top \in \mscr{X}$ and $\bs{\vse} = [\vse_0, \vse_1, \cd, \vse_\ka]^\top \in
\bs{\varTheta}$,  denote the pmf defined by (\ref{geninvhhyper89}) with $\mcal{C}_i$ replaced by $\vse_i$ for $i = 0, 1, \cd, \ka$.  Define
$\mscr{M} (\bs{x}, \bs{\vse}) = \f{ f(\bs{x}) }{ g (\bs{x}, \bs{\vse}) }$ for $\bs{x} \in \mscr{X}$ and  $\bs{\vse}  \in \bs{\varTheta}$ such
that $g (\bs{x}, \bs{\vse})
> 0$.  Define $\mscr{S} = \{ \bs{x} \in \mscr{X}: \; f(\bs{x}) > 0 \}$ and \[
\mscr{A} = \{ \bs{\vse} \in \bs{\varTheta}: | \bs{\vse} | \bs{\prec} | \bs{\mcal{C}} |, \; g(\bs{z}, \bs{\vse}) > 0 \}, \qqu \mscr{B} = \{
\bs{\vse} \in \bs{\varTheta}: | \bs{\vse} | \bs{\succ} | \bs{\mcal{C}} |, \; g(\bs{z}, \bs{\vse}) > 0 \}.
\]

\beL \la{pre89663mulinv} Assume that $f(\bs{z}) > 0$.  Then, $f(\bs{X}) \bb{I}_{\{\bs{X} \bs{\prec} \bs{z} \}} \leq \mscr{M} (\bs{z}, \bs{\vse})
g(\bs{X}, \bs{\vse})$ holds for all $\bs{\vse} \in \mscr{A}$. Similarly, $f(\bs{X}) \bb{I}_{\{\bs{X} \bs{\succ} \bs{z} \}} \leq \mscr{M}
(\bs{z}, \bs{\vse}) g(\bs{X}, \bs{\vse})$ holds for all $\bs{\vse} \in \mscr{B}$.

\eeL

\bpf First, we shall show that  $f(\bs{X}) \bb{I}_{\{\bs{X} \bs{\prec} \bs{z} \}} \leq \mscr{M} (\bs{z}, \bs{\vse}) g(\bs{X}, \bs{\vse})$ holds
for all $\bs{\vse} \in \mscr{A}$.  Since $f(\bs{z}) = g(\bs{z}, \bs{\mcal{C}})$, the set $\mscr{A}$ is nonempty as a consequence of the
assumption that $f(\bs{z}) > 0$.  Hence, $\mscr{M} (\bs{z}, \bs{\vse})$ is well-defined for $\bs{z} \in \mscr{S}$ and $\bs{\vse} \in \mscr{A}$.
Note that $f(\bs{x}) \bb{I}_{\{\bs{x} \bs{\prec} \bs{z} \}} \leq \mscr{M} (\bs{z}, \bs{\vse}) g(\bs{x}, \bs{\vse})$ holds for $\bs{\vse} \in
\mscr{A}$ and $\bs{x} \in \mscr{X}$ such that $f(\bs{x}) = 0$. To show the first assertion,  it suffices to show that $f(\bs{x})
\bb{I}_{\{\bs{x} \bs{\prec} \bs{z} \}} \leq \mscr{M} (\bs{z}, \bs{\vse}) g(\bs{x}, \bs{\vse})$ holds for $\bs{\vse} \in \mscr{A}$ and $\bs{x}
\in \mscr{S}$.
 We need to consider two cases as follows.

Case (I): $\mcal{N} > 0$.

Case (II): $\mcal{N} < 0$.

In  Case (I), as a consequence of $\mcal{N} > 0$, we have $\vse_i > 0, \; \mcal{C}_i > 0$ for $i = 0, 1, \cd, \ka$. Consider $\bs{x} \in
\mscr{S}$ such that $\bs{x} \bs{\prec} \bs{z}$.  Note that $\prod_{i=1}^k \prod_{\ell = 1}^{x_i} (\vse_i - \ell + 1) \geq \prod_{i=1}^k
\prod_{\ell = 1}^{z_i} (\vse_i - \ell + 1) > 0 $ because of $\bs{x} \bs{\prec} \bs{z}$ and $g(\bs{z}, \bs{\vse}) > 0$.   Hence, $g(\bs{x},
\bs{\vse}) > 0$ for $\bs{\vse} \in \mscr{A}$ and $\bs{x} \in \mscr{S}$ such that $\bs{x} \bs{\prec} \bs{z}$.  Since $f(\bs{x}) > 0$, we have
that
\[
\mscr{M} (\bs{x}, \bs{\vse}) = \f{ f(\bs{x}) }{ g (\bs{x}, \bs{\vse}) } = \prod_{i=1}^\ka \f{ \bi{\mcal{C}_i}{x_i}  }{ \bi{\vse_i}{x_i} } =
\prod_{i=1}^k  \prod_{\ell = 1}^{x_i}  \f{ \mcal{C}_i - \ell  + 1  }{ \vse_i - \ell + 1 }
\]
is well-defined and positive for $\bs{\vse} \in \mscr{A}$ and $\bs{x} \in \mscr{S}$ such that $\bs{x} \bs{\prec} \bs{z}$. Hence, we have a
meaningful ratio \be \la{case899inv}
 \f{ \mscr{M}
(\bs{z}, \bs{\vse})  }{ \mscr{M} (\bs{x}, \bs{\vse})  } =  \prod_{i=1}^k  \prod_{\ell = x_i + 1}^{z_i}  \f{ \mcal{C}_i - \ell  + 1  }{ \vse_i -
\ell + 1 } \ee for $\bs{x} \in \mscr{S}$ such that $\bs{x} \bs{\prec} \bs{z}$. Since $\vse_i \leq \mcal{C}_i$ for $i = 1, \cd, \ka$, it follows
from (\ref{case899inv}) that $\mscr{M} (\bs{x}, \bs{\vse} ) \leq \mscr{M} (\bs{z}, \bs{\vse})$ holds for $\bs{\vse} \in \mscr{A}$ and $\bs{x}
\in \mscr{S}$ such that $\bs{x} \bs{\prec} \bs{z}$.

In  Case (II), as a consequence of $\mcal{N} < 0$, we have $\vse_i < 0, \; \mcal{C}_i < 0$ for $i = 0, 1, \cd, \ka$. Note that
\[
\mscr{M} (\bs{x}, \bs{\vse}) =  \prod_{i=0}^k  \prod_{\ell = 1}^{x_i} \f{ \mcal{C}_i - \ell  + 1 }{ \vse_i - \ell + 1  } = \prod_{i=0}^k
\prod_{\ell = 1}^{x_i} \f{  |\mcal{C}_i| + \ell  - 1 }{ |\vse_i|  + \ell - 1 } > 0 \qu \tx{for $\bs{x} \in \mscr{X}$}
\]
and \be \la{case899binv}
 \f{ \mscr{M}
(\bs{z}, \bs{\vse})  }{ \mscr{M} (\bs{x}, \bs{\vse})  } =  \prod_{i=1}^k  \prod_{\ell = x_i + 1}^{z_i}  \f{ | \mcal{C}_i | + \ell  - 1  }{ |
\vse_i | + \ell - 1 }  \qqu  \tx{ for $\bs{x} \in \mscr{X}$ such that} \; \bs{x} \bs{\prec} \bs{z}. \ee Since $| \vse_i | \leq | \mcal{C}_i |$
for $i = 1, \cd, \ka$, it follows from (\ref{case899binv}) that $\mscr{M} (\bs{x}, \bs{\vse} ) \leq \mscr{M} (\bs{z}, \bs{\vse})$ holds for
$\bs{x} \in \mscr{S}$ such that $\bs{\vse} \in \mscr{A}$ and $\bs{x} \bs{\prec} \bs{z}$.  Therefore, in both cases, we have that
\[
\f{ f(\bs{x})}{g(\bs{x}, \bs{\vse})} = \mscr{M} (\bs{x}, \bs{\vse} )  \leq \mscr{M} (\bs{z}, \bs{\vse} )
\]
holds for $\bs{\vse} \in \mscr{A}$ and $\bs{x} \in \mscr{S}$ such that $\bs{x} \bs{\prec} \bs{z}$. Thus, we have shown that $f(\bs{x})
\bb{I}_{\{\bs{x} \bs{\prec} \bs{z} \}} \leq \mscr{M} (\bs{z}, \bs{\vse}) g(\bs{x}, \bs{\vse})$ holds for $\bs{\vse} \in \mscr{A}$ and  $\bs{x}
\in \mscr{S}$.  This completes the proof of the first assertion of the lemma.

The second assertion of the lemma can be shown in a similar manner.  Since $f(\bs{z}) = g(\bs{z}, \bs{\mcal{C}})$, the set $\mscr{B}$ is
nonempty as a consequence of the assumption that $f(\bs{z}) > 0$.  Hence, $\mscr{M} (\bs{z}, \bs{\vse})$ is well-defined for $\bs{z} \in
\mscr{S}$ and $\bs{\vse} \in \mscr{B}$. Note that $f(\bs{x}) \bb{I}_{\{\bs{x} \bs{\succ} \bs{z} \}} \leq \mscr{M} (\bs{z}, \bs{\vse}) g(\bs{x},
\bs{\vse})$ holds for $\bs{\vse} \in \mscr{B}$ and $\bs{x} \in \mscr{X}$ such that $f(\bs{x}) = 0$. To show the second assertion,  it suffices
to show that $f(\bs{x}) \bb{I}_{\{\bs{x} \bs{\succ} \bs{z} \}} \leq \mscr{M} (\bs{z}, \bs{\vse}) g(\bs{x}, \bs{\vse})$ holds for $\bs{\vse} \in
\mscr{B}$ and $\bs{x} \in \mscr{S}$.
 We need to consider two cases as follows.

Case (I): $\mcal{N} > 0$.

Case (II): $\mcal{N} < 0$.

In  Case (I), as a consequence of $\mcal{N} > 0$, we have $\vse_i > 0, \; \mcal{C}_i > 0$ for $i = 0, 1, \cd, \ka$. Consider $\bs{x} \in
\mscr{S}$ such that $\bs{x} \bs{\succ} \bs{z}$.  Note that $\prod_{i=1}^k \prod_{\ell = 1}^{x_i} (\vse_i - \ell + 1) \geq \prod_{i=1}^k
\prod_{\ell = 1}^{x_i} (\mcal{C}_i - \ell + 1) > 0 $ because of $\bs{x} \bs{\succ} \bs{z}$ and $f (\bs{x}) > 0$.   Hence, $g(\bs{x}, \bs{\vse})
> 0$ for $\bs{\vse} \in \mscr{B}$ and $\bs{x} \in \mscr{S}$ such that $\bs{x} \bs{\succ} \bs{z}$.  Since $f(\bs{x}) > 0$, we have that
\[
\mscr{M} (\bs{x}, \bs{\vse}) = \f{ f(\bs{x}) }{ g (\bs{x}, \bs{\vse}) } = \prod_{i=0}^\ka \f{ \bi{\mcal{C}_i}{x_i}  }{ \bi{\vse_i}{x_i} } =
\prod_{i=0}^k  \prod_{\ell = 1}^{x_i}  \f{ \mcal{C}_i - \ell  + 1  }{ \vse_i - \ell + 1 }
\]
is well-defined and positive for $\bs{\vse} \in \mscr{B}$ and $\bs{x} \in \mscr{S}$ such that $\bs{x} \bs{\succ} \bs{z}$. Hence,  we have a
meaningful ratio \be \la{case899invB8}
 \f{ \mscr{M} (\bs{x}, \bs{\vse})  }{ \mscr{M} (\bs{z}, \bs{\vse})  } =  \prod_{i=1}^k  \prod_{\ell = z_i + 1}^{x_i}  \f{ \mcal{C}_i - \ell  + 1  }{ \vse_i -
\ell + 1 } \ee for $\bs{x} \in \mscr{S}$ such that $\bs{x} \bs{\succ} \bs{z}$. Since $ \vse_i \geq \mcal{C}_i$ for $i = 1, \cd, \ka$, it follows
from (\ref{case899invB8}) that $\mscr{M} (\bs{x}, \bs{\vse} ) \leq \mscr{M} (\bs{z}, \bs{\mcal{C}})$ holds for $\bs{\vse} \in \mscr{B}$ and
$\bs{x} \in \mscr{S}$ such that $\bs{x} \bs{\succ} \bs{z}$.

In  Case (II), as a consequence of $\mcal{N} < 0$, we have $\vse_i < 0, \; \mcal{C}_i < 0$ for $i = 0, 1, \cd, \ka$. Note that
\[
\mscr{M} (\bs{x}, \bs{\vse}) =  \prod_{i=0}^k  \prod_{\ell = 1}^{x_i} \f{ \mcal{C}_i - \ell  + 1 }{ \vse_i - \ell + 1  } = \prod_{i=0}^k
\prod_{\ell = 1}^{x_i} \f{  |\mcal{C}_i| + \ell  - 1 }{ |\vse_i|  + \ell - 1 } > 0 \qu \tx{for $\bs{x} \in \mscr{X}$}
\]
and \be \la{case899binvB899}
 \f{ \mscr{M}
(\bs{x}, \bs{\vse})  }{ \mscr{M} (\bs{z}, \bs{\vse})  } =  \prod_{i=1}^k  \prod_{\ell = z_i + 1}^{x_i}  \f{ | \mcal{C}_i | + \ell  - 1  }{ |
\vse_i | + \ell - 1 }  \qqu  \tx{ for $\bs{x} \in \mscr{X}$ such that} \; \bs{x} \bs{\succ} \bs{z}. \ee Since $ | \vse_i | \geq | \mcal{C}_i |$
for $i = 1, \cd, \ka$, it follows from (\ref{case899binvB899}) that $\mscr{M} (\bs{x}, \bs{\vse} ) \leq \mscr{M} (\bs{z}, \bs{\vse})$ holds for
$\bs{\vse} \in \mscr{B}$ and $\bs{x} \in \mscr{X}$ such that  $\bs{x} \bs{\succ} \bs{z}$.  Therefore, in both cases, we have that
\[
\f{ f(\bs{x})}{g(\bs{x}, \bs{\vse})} = \mscr{M} (\bs{x}, \bs{\vse} )  \leq \mscr{M} (\bs{z}, \bs{\vse} )
\]
holds for $\bs{\vse} \in \mscr{B}$ and $\bs{x} \in \mscr{S}$ such that $\bs{x} \bs{\succ} \bs{z}$. Thus, we have shown that $f(\bs{x})
\bb{I}_{\{\bs{x} \bs{\succ} \bs{z} \}} \leq \mscr{M} (\bs{z}, \bs{\vse}) g(\bs{x}, \bs{\vse})$ holds for $\bs{\vse} \in \mscr{B}$ and  $\bs{x}
\in \mscr{S}$.  This completes the proof of the second assertion of the lemma.

\epf

Applying  Theorem \ref{ThM888} and Lemma \ref{pre89663mulinv}, we have the following results.

\beL

\la{generhyp86389inv}

Assume that $f(\bs{z}) > 0$.  Then,
\[
\Pr \{ \bs{X} \bs{\prec} \bs{z} \}  \leq \prod_{i=0}^\ka \f{ \bi{ \mcal{C}_i  }{ z_i } } { \bi{\vse_i}{z_i}  } \qu \tx{for $\bs{\vse} \in
\mscr{A}$},
\]
\[
\Pr \{ \bs{X} \bs{\succ} \bs{z} \}  \leq \prod_{i=0}^\ka \f{ \bi{ \mcal{C}_i  }{ z_i } } { \bi{\vse_i}{z_i}  } \qu \tx{for $\bs{\vse} \in
\mscr{B}$}.
\]

\eeL

\beL

\la{proper68833inv}

 $g(\bs{z}, \wh{\bs{\mcal{C}}}) > 0$ for $\bs{z} \in \mscr{X}$.  \eeL

\bpf

Recall that $\wh{\bs{\mcal{C}}} = [\wh{\mcal{C}}_0, \wh{\mcal{C}}_1, \cd, \wh{\mcal{C}}_\ka]^\top$ with $\wh{\mcal{C}}_i = \mcal{N} \f{z_i}{n}$
for $i = 0, 1, \cd, \ka$.

 In the case of $\mcal{N} < 0$, we have $\wh{\mcal{C}}_i \leq 0$ for $i = 0, 1, \cd, \ka$, which implies
that $1 + \wh{\mcal{C}}_i \leq 1$ for $i = 0, 1, \cd, \ka$. Hence, $g(\bs{z}, \wh{\bs{\mcal{C}}}) > 0$ for $\bs{z} \in \mscr{X}$ if $\mcal{N} <
0$.

In the case of $\mcal{N} > 0$, as a consequence of  $\f{n - 1}{\mcal{N}} < 1$, we have $0 < n < \mcal{N} + 1$ and $\wh{\mcal{C}}_i \geq 0$ for
$i = 0, 1, \cd, \ka$. For $i = 0, 1, \cd, \ka$, we have $z_i < \wh{\mcal{C}}_i + 1 = 1$ if $z_i = 0$; and $ z_i = \f{n}{\mcal{N}}
\wh{\mcal{C}}_i < \wh{\mcal{C}}_i + \f{ \wh{\mcal{C}}_i }{\mcal{N}} = \wh{\mcal{C}}_i + \f{z_i}{n} \leq \wh{\mcal{C}}_i + 1$ if $z_i > 0$.
Hence, $g(\bs{z}, \wh{\bs{\mcal{C}}}) > 0$ for $\bs{z} \in \mscr{X}$ if $\mcal{N} > 0$.

\epf

\beL

\la{more89663inv}

$| \wh{\bs{\mcal{C}}} | \bs{\prec} | \bs{\mcal{C}} | \Leftrightarrow \bs{z} \bs{\prec} \wh{\bs{\mu}}$ and $\wh{\bs{\mcal{C}}} \bs{\succ} |
\bs{\mcal{C}} |  \Leftrightarrow \bs{z} \bs{\succ} \wh{\bs{\mu}}$.

\eeL

\bpf

In the case of $\mcal{N} > 0$, we have $\mcal{C}_i > 0$ and $\wh{\mcal{C}}_i \geq 0$ for $i = 0, 1, \cd, \ka$. Hence,
\[
| \wh{\bs{\mcal{C}}} | \bs{\prec} | \bs{\mcal{C}} | \Leftrightarrow \wh{\bs{\mcal{C}}} \bs{\prec} \bs{\mcal{C}} \Leftrightarrow \bs{z}
\bs{\prec} \wh{\bs{\mu}}, \qqu | \wh{\bs{\mcal{C}}} | \bs{\succ} | \bs{\mcal{C}} | \Leftrightarrow \wh{\bs{\mcal{C}}} \bs{\succ} \bs{\mcal{C}}
\Leftrightarrow \bs{z} \bs{\succ} \wh{\bs{\mu}}.
\]

In the case of $\mcal{N} < 0$, we have $\mcal{C}_i < 0$ and $\wh{\mcal{C}}_i \leq 0$ for $i = 0, 1, \cd, \ka$. Hence,
\[
| \wh{\bs{\mcal{C}}} | \bs{\prec} | \bs{\mcal{C}} | \Leftrightarrow \wh{\bs{\mcal{C}}} \bs{\succ} \bs{\mcal{C}} \Leftrightarrow \bs{z}
\bs{\prec} \wh{\bs{\mu}}, \qqu | \wh{\bs{\mcal{C}}} | \bs{\succ} | \bs{\mcal{C}} |  \Leftrightarrow \wh{\bs{\mcal{C}}} \bs{\prec} \bs{\mcal{C}}
\Leftrightarrow \bs{z} \bs{\succ} \wh{\bs{\mu}}.
\]

\epf

\bsk

We are now in a position to prove the theorem. Clearly, the assertions (I) and (II) of Theorem \ref{GenmulInvHyper} follow from Lemmas
\ref{generhyp86389inv}, \ref{proper68833inv} and \ref{more89663inv}.  The assertions (III) and (IV) of Theorem \ref{GenmulInvHyper} can be shown
by the same arguments as that of the counterparts of Theorem \ref{GenmulHyper}.

\sect{Proof of Theorem \ref{mulDrich}} \la{mulDrichapp}

Define $\vse_0 = \al_0$,
\[
\varTheta = \{ [\vse_0, \vse_1, \cd, \vse_\ka]^\top: 0 < \vse_i \leq \al_i, \; i = 1, \cd, \ka \}
\]
and
\[
\ro ( \bs{x}, \bs{\vse}) = \ln \f{ \mcal{B}(\bs{\vse})  }{  \mcal{B}(\bs{\al}) } +  \sum_{i=1}^\ka ( \al_i - \vse_i ) \ln x_i.
\]
It can be seen that $\ro ( \bs{x}, \bs{\vse})$ is a concave function of $\bs{x}$ such that $\ro ( \bs{x}, \bs{\vse}) \leq \ro ( \bs{z},
\bs{\vse})$ for all $\bs{\vse} \in \varTheta$ provided that $\bs{x} \bs{\prec} \bs{z}$.  Moreover, $f(\bs{X}, \bs{\al}) = \exp ( \ro ( \bs{X},
\bs{\vse}) ) f (\bs{X}, \bs{\vse})$ for all $\bs{\vse} \in \varTheta$.  Making use of these facts and Theorem \ref{con89666}, we have
\[
\Pr \{ \ovl{\bs{\mcal{X}}}_n \bs{\prec} \bs{z} \} \leq \inf_{\bs{\vse} \in \varTheta} \li [ \f{ \mcal{B}(\bs{\vse})  }{  \mcal{B}(\bs{\al}) }
\prod_{i=1}^\ka \f{ z_i^{\al_i} } { z_i^{\vse_i} } \ri ]^n.
\]
To derive an explicit bound, we shall use the method of moments by seeking $\bs{\vse}$ such that $\bs{z} = \bb{E}_{\bs{\vse}} [\bs{X}]$. From
this equation, we have
\[
z_i = \f{\vse_i}{\sum_{\ell =0}^\ka \vse_\ell}, \qqu i = 1, \cd, \ka,
\]
which implies that $\f{\vse_i}{z_i}$ are equal for all $i = 1, \cd, \ka$.  It follows that
\[
z_i = \f{1}{ \f{\al_0}{\vse_i} + \sum_{\ell = 0}^\ka  \f{\vse_\ell}{\vse_i}  } = \f{1}{ \f{\al_0}{\vse_i} + \sum_{\ell = 0}^\ka \f{z_\ell}{z_i}
},  \qqu i = 1, \cd, \ka,
\]
from which we obtain
\[
\f{\al_0}{\vse_i} = \f{1}{z_i} - \sum_{\ell = 0}^\ka \f{z_\ell}{z_i}, \qqu i = 1, \cd, \ka. \] Hence,
\[
\vse_i = \f{ \al_0  z_i}{ 1 - \sum_{\ell = 1}^\ka z_\ell } = \f{ \al_0  z_i}{ z_0} = \wh{\al}_i, \qqu i = 1, \cd, \ka.
\]
Noting that $z_0 \geq \mu_0 = \f{\al_0}{ \sum_{\ell = 0}^\ka \al_\ell }$,  we have $\f{ \al_0}{ z_0} \leq \sum_{\ell = 0}^\ka \al_\ell$ and thus
\[
\wh{\al}_i = \vse_i \leq  z_i \sum_{\ell = 0}^\ka \al_\ell \leq \mu_i \sum_{\ell = 0}^\ka \al_\ell = \al_i, \qqu i = 1, \cd, \ka.
\]
This establishes (\ref{conthm1388}).  The proof of the theorem is thus completed.

\sect{Proof of Theorem  \ref{matrixgga} } \la{matrixggaapp}

Define
\[
g (\bs{x}, \bs{\vse}) = \f{| \bs{\vse} |^{-\al}}{ \ba^{p \al}   \Ga_p (\al ) } | \bs{x} |^{\al -( p + 1) \sh 2} \exp \li (  - \f{1}{\ba} \tx{tr}
( \bs{\vse}^{-1} \bs{x} )\ri ),
\]
where $\bs{\vse} \bs{\succ} \bs{\Si}$.  Define
\[
\ro (\bs{x}, \bs{\vse}) = \al \ln \f{ | \bs{\vse} | } { | \bs{\Si} | }  - \f{1}{\ba} \tx{tr} ( [ \bs{\Si}^{-1} - \bs{\vse}^{-1} ] \bs{x} ).
\]
Then, $\ro (\bs{x}, \bs{\vse}) = \f{ f(\bs{x})  }{  g (\bs{x}, \bs{\vse})  }$ is a linear function of $\bs{x}$ such that $f(\bs{X}) = \exp(  \ro
(\bs{X}, \bs{\vse})  ) \; g(\bs{X}, \bs{\vse})$.  Clearly, $\bs{\Si}^{-1} \bs{\succ} \bs{\vse}^{-1}$.  Let $A = \sq{\bs{\Si}^{-1} -
\bs{\vse}^{-1} }$ and $B = \sq{\bs{y} - \bs{x}}$.  Note that \[ \tx{tr} \li [ ( \bs{\Si}^{-1} - \bs{\vse}^{-1} ) ( \bs{x} - \bs{z} ) \ri ]  =
\tx{tr} \li (  A A B B \ri ) = \tx{tr} ( BAA B ) =  \tx{tr} \li [  (A B)^\top (AB) \ri ] \geq  0  \] for $\bs{x} \bs{\succ} \bs{z}$.  It follows
that $\ro (\bs{x}, \bs{\vse}) \leq \ro (\bs{z}, \bs{\vse})$ for $\bs{x} \bs{\succ} \bs{z}$.  Making use of these facts and Theorem
\ref{con89666}, we have \be \la{mga89}
 \Pr \{ \ovl{\bs{\mcal{X}}}_n \bs{\succ} \bs{z} \} \leq \inf_{ \bs{\vse} \bs{\succ} \bs{\Si} } \li [  \f{
f(\bs{z}) } {  g (\bs{z}, \bs{\vse}) } \ri ]^n. \ee Now let $\bs{\vse} = \f{ \bs{z} }{\al \ba}$.  Since $\bs{z} \bs{\succ} \al \ba \bs{\Si}$, we
have $\bs{\vse} = \f{ \bs{z} }{\al \ba} \bs{\succ} \bs{\Si}$.  Note that \bee  \f{ f(\bs{z})  } {  g (\bs{z}, \f{ \bs{z} }{\al \ba} )  }  =  \li
(  \f{e}{\al \ba} \ri )^{ p \al } \li ( \f{ | \bs{z} | } { | \bs{\Si} | } \ri )^{\al} \exp \li (   - \f{1}{\ba} \tx{tr} ( \bs{\Si}^{-1} \bs{z} )
\ri ). \eee It follows from (\ref{mga89}) that
\[
\Pr \{ \ovl{\bs{\mcal{X}}}_n \bs{\succ} \bs{z} \} \leq  \li [ \li (  \f{e}{\al \ba} \ri )^{ p \al } \li (  \f{ | \bs{z} | } { | \bs{\Si} | } \ri
)^{\al} \exp \li ( - \f{1}{\ba} \tx{tr} ( \bs{\Si}^{-1} \bs{z} ) \ri ) \ri ]^n \qu \tx{for $\bs{z} \bs{\succ} \al \ba \bs{\Si}$}.
\]
 Taking $\bs{z} = \vro \al \ba \bs{\Si}$, we have \bee  \Pr \{ \bs{X} \bs{\succ} \vro \al \ba \bs{\Si} \} & \leq  &
\li [ \li ( \f{e}{\al \ba} \ri )^{ p \al } \li (  \f{ | \vro \al \ba \bs{\Si} | } { | \bs{\Si} | } \ri )^{\al} \exp \li (   - \f{1}{\ba} \tx{tr}
( \bs{\Si}^{-1} \vro \al \ba \bs{\Si} ) \ri ) \ri ]^n\\
& = & \li [  \vro \exp (1 - \vro ) \ri ]^{n p \al} \eee for $\vro  \geq 1$.  In a similar manner, we can show that
\[
\Pr \{ \ovl{\bs{\mcal{X}}}_n \bs{\prec} \bs{z} \} \leq  \li [ \li (  \f{e}{\al \ba} \ri )^{ p \al } \li (  \f{ | \bs{z} | } { | \bs{\Si} | } \ri
)^{\al} \exp \li ( - \f{1}{\ba} \tx{tr} ( \bs{\Si}^{-1} \bs{z} ) \ri ) \ri ]^n \qu \tx{for $\bs{z} \bs{\prec} \al \ba \bs{\Si}$},
\] and  $\Pr \{ \bs{X} \bs{\prec} \vro \al \ba \bs{\Si} \} \leq  \li [  \vro \exp (1 - \vro ) \ri ]^{n p \al}$ for $\vro \in (0, 1]$. This completes the proof of the theorem.

\end{document}